\documentclass[11pt]{article} 
\usepackage[preprint]{tmlr}
\usepackage{microtype}  

\usepackage{xcolor}
\usepackage{fourier}   

\usepackage{caption}
\usepackage{amsmath}
\usepackage{amsfonts} 
\usepackage{amsthm} 
\usepackage{mathtools}

\usepackage{amsmath,amsfonts,bm}

\def\eqref#1{equation~\ref{#1}}

\def\1{\bm{1}}

\def\vv{{\bm{v}}}

\def\mR{{\bm{R}}}

\DeclareMathAlphabet{\mathsfit}{\encodingdefault}{\sfdefault}{m}{sl}
\SetMathAlphabet{\mathsfit}{bold}{\encodingdefault}{\sfdefault}{bx}{n}

\newcommand{\E}{\mathbb{E}}

\usepackage[e]{esvect} 
\PassOptionsToPackage{hyphens}{url} 
\usepackage{hyperref}
\hypersetup{
    linktocpage,
    colorlinks,
    linkcolor={blue!80!black},
    citecolor={green!50!black},
    urlcolor={blue!80!black}
}
\newcommand{\Def}[1]{{\color{red!60!black} #1}}
\usepackage{url}

\usepackage{graphicx}
\usepackage{sidecap}

\newtheorem{theorem}{Theorem}[section]       
\newtheorem{proposition}[theorem]{Proposition}
\newtheorem{definition}[theorem]{Definition}
\newtheorem{corollary}[theorem]{Corollary}
\newtheorem{lemma}[theorem]{Lemma}
\newtheorem{remark}[theorem]{Remark}
\newtheorem{example}[theorem]{Example}

\usepackage{enumitem}

\usepackage{xspace}
\newcommand*{\eg}{e.g.\@\xspace}
\newcommand*{\ie}{i.e.\@\xspace}

\newcommand*{\cf}{cf.\@\xspace}

\renewcommand{\mR}{\mathbb{R}}
\newcommand{\Pup}{\overline{P}}
\newcommand{\Plow}{\underline{P}}
\newcommand{\Rup}{\overline{R}}
\newcommand{\Rlow}{\underline{R}}

\newcommand{\avgseq}{\vv{\Sigma X}}

\newcommand{\CP}{\operatorname{CP}}
\newcommand{\limsupn}{\limsup_{n \rightarrow \infty}}

\newcommand{\dynkin}{\mathcal{G}}
\newcommand{\gbr}{\operatorname{GBR}}
\newcommand{\vvomi}{\vv{\Omega}(i)}
\newcommand{\linfty}{L^\infty}
\newcommand{\linftyd}{(\linfty)^*}

\newcommand{\natext}{\operatorname{NatExt}}
\newcommand{\comp}{\mathrm{C}}
\newcommand{\eventonce}{2^{\vv{\Omega}}_{1+}}
\newcommand{\dynkingen}{\operatorname{PD}}
\newcommand{\PF}{\operatorname{PF}}

\newcommand{\naturals}{\mathbb{N}}
\newcommand{\reals}{\mathbb{R}}

\newcommand{\pnew}{{p_{\mathrm{new}}}}

\newcommand{\hpnew}{{\hat{p}_{\mathrm{new}}}}
\newcommand{\pold}{{p_{\mathrm{old}}}}
\newcommand{\relint}{\operatorname{relint}}
\newcommand{\co}{\operatorname{co}}

\newcommand{\Scal}{\mathcal{S}}
\newcommand{\Ncal}{\mathcal{N}}
\newcommand{\Scalt}{\tilde{\mathcal{S}}}
\newcommand{\Ncalt}{\tilde{\mathcal{N}}}
\usepackage{algpseudocode,algorithm,stmaryrd,caption,subcaption,graphicx}
\graphicspath{{./figures/}}

\newcommand{\coclosed}{\overline{\operatorname{co}}}

\renewcommand*{\eqref}[1]{Equation~\ref{#1}}

\title{Strictly Frequentist Imprecise Probability}  

\author{\name Christian Fröhlich \email christian.froehlich@uni-tuebingen.de \\
      \addr Department of Computer Science\\
      University of Tübingen
      \AND
      \name Rabanus Derr \email rabanus.derr@uni-tuebingen.de \\
      \addr Department of Computer Science\\
      University of Tübingen
      \AND
      \name Robert C. Williamson \email bob.williamson@uni-tuebingen.de\\
      \addr Department of Computer Science \\
      University of Tübingen\\
      and Tübingen AI Center}

\def\month{MM}  
\def\year{YYYY} 
\def\openreview{\url{https://openreview.net/forum?id=XXXX}} 

\begin{document}

\maketitle

\begin{abstract}
Strict frequentism defines probability as the limiting relative
frequency in an infinite sequence. What if the limit does 
not exist? We present a broader theory, which is applicable 
also to random phenomena that exhibit diverging 
relative frequencies. In doing so, we develop a close 
connection with the theory of imprecise probability: 
the cluster points of relative frequencies yield a coherent upper prevision. We show that a natural 
frequentist definition of conditional probability recovers 
the generalized Bayes rule. This also suggests an 
independence concept, which is related to epistemic 
irrelevance in the imprecise probability literature.
Finally, we prove constructively that, for a finite set of elementary events, there exists a sequence for 
which the cluster points of relative frequencies coincide 
with a prespecified set which demonstrates the naturalness, and arguably completeness, of our theory.

\end{abstract}


\section{Introduction}
\label{sec:introductionfreqip}

\vspace*{-7mm}
\hfill\begin{minipage}{0.55\textwidth}
    \footnotesize{\it
    Do other statistical properties, that can not be reduced
    to stochasticness, exist? This question did not attract
    any attention until the applications of the probability
    theory concerned only natural sciences. The situation is
    definitely [changing], when one studies social phenomena: 
    the stochasticness gets broken as soon as we [deal] with 
    deliberate activity of people.}
    \hfill --- Victor Ivanenko and Valery Labkovsky \citeyearpar{ivanenko1993}\\
\end{minipage}

It is now almost universally acknowledged that probability theory 
ought to be based on Kolmogorov's
\citeyearpar{kolmogorov1933grundbegriffe} mathematical axiomatization (translated in~\citep{kolmogorov2018foundations}).\footnote{An important exception is quantum probability \citep{gudder1979stochastic, khrennikov2016probability}.} However, if probability is defined in this purely measure-theoretic fashion, what warrants its application to real-world problems of decision making under uncertainty? To those in the so-called \textit{frequentist} camp, the justification is essentially due to the \textit{law of large numbers}, which comes in both an empirical and a theoretical flavour.
Our motivation for the present paper comes from questioning both of these.

By the empirical version of the law of large numbers (LLN), we mean not a ``law'' which can be proven to hold, but the following hypothesis, which seems to guide many scientific endeavours. Assume we have obtained data $x_1,..,x_n$ as the outcomes of some experiment, which has been performed $n$ times under ``statistically identical'' conditions. Of course, conditions in the real-world can never truly be identical --- otherwise the outcomes would be constant, at least under the assumption of a deterministic universe. Thus, ``identical'' in this context must be a weaker notion, that all factors which we have judged as relevant to the problem at hand have been kept constant over the repetitions.\footnote{In fact, we do not need that conditions stay exactly constant, but that they change merely in a way which is so benign that the relative frequencies converge. That is, in the limit we should obtain a stable statistical aggregate.}
The empirical ``law'' of large numbers, which \citet{gorban2017statistical} calls the \textit{hypothesis of
(perfect) statistical stability} then asserts that in the long-run, relative frequencies of events and 
sample averages converge. These limits are then conceived of as the \textit{probability} of an 
event and the \textit{expectation}, respectively. Thus, even if relative frequencies can fluctuate in the 
finite data setting, we expect that they stabilize as more and more data is acquired. 
Crucially, this hypothesis of perfect statistical stability is not amenable to falsification, 
since we can never refute it in the finite data setting. It is a matter of faith to assume convergence 
of relative frequencies. On the other hand, there is now ample experimental evidence that relative frequencies 
can fail to stabilize even under very long observation intervals \citep[Part II]{gorban2017statistical}. 
We say that such phenomena display \textit{unstable (diverging)} relative frequencies. 
Rather than refuting the stability hypothesis, which is impossible, we question its adequateness as 
an idealized modeling assumption: we view convergence as the idealization of approximate stability in the finite case, whereas divergence idealizes instability. 
Thus, if probability is understood as limiting relative frequency, then the applicability of Kolmogorov's
theory to empirical phenomena is limited to those which are statistically stable; 
the founder himself remarked:
\begin{quote}
    Generally speaking there is no ground to believe that a random phenomenon should
possess any definite probability \citep{kolmogorov1983logical}.
\end{quote}
Building on the works of \citet{mises1964mathematical}, \citet{walley1982towards} and \citet{ivanenkobook}, our goal is to establish a broader theory, which is also applicable to ``random'' phenomena which are outside of the scope of Kolmogorov's theory by exhibiting unstable relative frequencies.

One attempt to ``prove'' (or justify) the empirical law of large numbers, which in our view is doomed to fail, is to invoke the theoretical law of large numbers, which is a purely formal, mathematical statement. The strong law of large numbers states that if $X_1,X_2,..$ is a sequence of independent and identically distributed (i.i.d.) random variables with finite expectation $\Def{\E[X] \coloneqq \E[X_1]=\E[X_2]=\cdots}$, then the sample average $\Def{\bar{X}_n \coloneqq \frac{1}{n} \sum_{i=1}^n X_i}$ converges almost surely to the expectation:
\[
    P\left(\lim_{n \rightarrow \infty} \bar{X}_n = \E[X]\right) = 1,
\]
where $P$ is the underlying probability measure in the sense of Kolmogorov. To interpret this statement correctly, some care is needed. It asserts that $P$ assigns measure $1$ to the set of sequences for which the sample mean converges, but not that this happens \textit{for all} sequences. Thus one would need justification for identifying ``set of measure 0`` with ``is negligible'' (``certainly does not happen''), which in particular requires a justification for $P$. With respect to a different measure, this set might not be negligible at all \citep[p.\@\xspace 8]{schnorr2007zufalligkeit}; see also \citep{calude1999most,whenlargeisalso} for critical arguments. Moreover, the examples in \citep[Part II]{gorban2017statistical} show that sequences with seemingly non-converging relative frequencies (fluctuating substantially even for long observation intervals) are not ``rare'' in practice. In Appendix~\ref{app:pathornormal} we examine the question of how pathological or normal such sequences are in more depth.

Conceptually, the underlying problem is that the probability measure $P$, which is used to measure the event 
$\left\{\lim_{n \rightarrow \infty} \bar{X}_n = \E[X]\right\}$ has no clear meaning. Of course, in the subjectivist spirit, one could interpret it as assigning a \textit{belief} in the statement that convergence takes place. But it is unclear what a frequentist interpretation of $P$ would look like. As \citet{lacazefrequentism} observed:
\begin{quote}
    Importantly, “almost sure convergence” is also given a frequentist interpretation. Almost
sure convergence is taken to provide a justification for assuming that the relative frequency
of an attribute \textit{would} converge to the probability in actual experiments \textit{were} the experiment
to be repeated indefinitely [emphasis in original].
\end{quote}
But again, it is unclear on what ground $P$ can be given this interpretation and according to \citet{hajek2009fifteen} this leads to a regress to mysterious ``meta-probabilities''. Furthermore, the theoretical LLN requires that $P$ be countably additive, which is problematic under a frequency interpretation \citep[pp.\@\xspace 229--230]{hajek2009fifteen}.

Given these complications, we opt for a different approach, namely a \textit{strictly frequentist} one. Reaching back to Richard von Mises' \citeyearpar{mises1919grundlagen} foundational work, a strictly frequentist theory explicitly defines probability in terms of limiting relative frequencies in a sequence. Importantly, we here \textit{do not assume} that the elements of the sequence are random variables with respect to an abstract, countably additive probability measure. Instead, like von Mises, we actually take the notion of a sequence as the primitive entity in the theory. As a consequence, countable additivity does not naturally arise in this setting, and hence we do not subscribe to the frequentist interpretation of the classical strong LLN. 

The core motivation for our work is to drop the assumption 
of perfect statistical stability and instead to explicitly 
model the possibility of unstable (diverging) 
relative frequencies.
Rather than merely conceding that the ``probability'' might 
vary over time \cite[pp.\@\xspace 27ff.]{borel1950} (which begs 
the question what such ``probabilities'' mean) we follow 
the approach of  Ivanenko \citeyearpar{ivanenkobook}, reformulate his construction of a \textit{statistical regularity} of a sequence, 
and discover that it is closely connected to the subjectivist theory of \textit{imprecise probability}. 
In essence, to each sequence we can naturally associate a \emph{set} of probability measures, 
which constitute the statistical regularity that describes the cluster points of relative frequencies and consequently also those of sample averages. Since this works for \textit{any} sequence and \textit{any} event, we have thus countered a typical argument against frequentism, namely that the limit may not exist and hence probability is undefined \citep{hajek2009fifteen}. The relative frequencies induce a coherent upper probability and the sample averages induce a coherent upper prevision in the sense of \citet{walley1991statistical}. In the convergent case, this reduces to a precise, finitely additive probability and a linear prevision, respectively. 
Furthermore, we derive in a natural way a conditional upper prevision;
remarkably, this approach recovers the \textit{generalized Bayes rule}, the arguably most important updating principle in imprecise probability.

Furthermore, we demonstrate that the reverse direction works, too: given a set of probability measures, we can explicitly construct a sequence, which corresponds to this set in the sense that its relative frequencies have this set of cluster points. Thereby we establish strictly frequentist \textit{semantics} for imprecise probability: a subjective decision maker who uses a set of probability measures to represent their belief can also be understood as assuming an implicit underlying sequence and reasoning in a frequentist way thereon.


\subsection{Von Mises - The Frequentist Perspective}
Our approach is inspired by Richard von Mises \citeyearpar{mises1919grundlagen} (refined and summarized in \citep{mises1964mathematical}) axiomatization of probability theory. 
In contrast to the subjectivist camp, von Mises concern was to develop a theory for repetitive events; which gives rise to a theory of probability that is mathematical, but which can also be used to reason about the physical world.
\begin{quote}
The calculus of probability, i.e. the theory of probabilities, in so far as
they are numerically representable, is the theory of definite observable
phenomena, repetitive or mass events. Examples are found in games of
chance, population statistics, Brownian motion etc.
\citep[p.\@\xspace 102]{mises1981probability}.
\end{quote}
Hence, von Mises is not concerned with the probability of single events, which he deems meaningless, but instead always views an event as part of a larger \textit{reference class}. Such a reference class is captured by what he terms a \textit{collective}, a disorderly sequence which exhibits both global regularity and local irregularity. 

\begin{definition}
Consider a tuple $\Def{\left(\Omega,\vv{\Omega},\mathcal{A},\mathcal{S}\right)}$ with the following data: 
\begin{enumerate}[nolistsep] \item a sequence $\Def{\vv{\Omega} \colon \mathbb{N} \rightarrow \Omega}$;
\item a set of \Def{\emph{selection rules}} $\Def{\mathcal{S} \coloneqq \{\vv{S}_j \colon j \in \mathcal{J}\}}$,
where  for each $j$ in a countable index set $\mathcal{J}$,  
$\vv{S}_j \colon \mathbb{N} \rightarrow \{0,1\}$
and $\vv{S}_j(i)=1$ for infinitely many $i \in \mathbb{N}$; 
\item   a non-empty set system 
$\Def{\mathcal{A} \subseteq 2^\Omega}$, where for simplicity we assume 
$|\mathcal{A}|<\infty$.\footnote{In fact, $\mathcal{A}$ does not necessarily has to be finite. 
Since an infinite domain of probabilities does not contribute a lot to a better understanding of 
the frequentist definition at this point, we restrict ourselves to the finite case here. The 
reader can find details in \citep{mises1964mathematical}.}
\end{enumerate}

This tuple forms a collective if the following two axioms hold.
\begin{enumerate}[nolistsep,label=\textbf{\emph{vM\arabic*.}}, ref=vM\arabic*]
   \item \label{vm1} The limiting relative frequency for $A \in \mathcal{A} \subseteq 2^\Omega$ exists:
   \[
       \Def{P(A) \coloneqq \lim_{n \rightarrow \infty} \frac{1}{n} \sum_{i=1}^n 
       \chi_A\left(\vvomi\right)}.\footnote{The function $\Def{\chi_A}$ denotes the indicator gamble for a set 
       $A \subseteq \Omega$, \ie $\Def{\chi_A(\omega)\coloneqq 1}$ if $\omega \in A$ and 
       $\Def{\chi_A(\omega)\coloneqq 0}$ otherwise.}
   \]
   We call this limit the \Def{\emph{probability of $A$}}.
   \item \label{vm2} For each $j \in \mathcal{J}$, the selection rule $\vv{S}_j$ does not change limiting relative frequencies:\footnote{To be precise, a selection rule in the sense of von Mises is a map from the set of finite $\Omega$-valued strings to $\{0,1\}$, \ie a selection rule is able to ``see'' all previous elements when deciding whether or not to select the next one. Our formulation is more restrictive to avoid notational overhead, but when a sequence is fixed, the two formulations are equivalent.}
   \[
   \lim_{n \rightarrow \infty} \frac{\sum_{i=1}^n  \chi_A\left(\vvomi\right) \cdot \vv{S}_j(i)}{\sum_{i=1}^n \vv{S}_j(i)} = P(A) \ \ \ \ \forall A\in\mathcal{A}.
   \]
\end{enumerate}
\end{definition}
Here, we view $\vv{\Omega}$ as a sequence of elementary outcomes $\omega \in \Omega$, 
for some possibility space $\Omega$ on which we have a set system of events $\mathcal{A}$.
Axiom \ref{vm1} explicitly defines the probability of an event in terms of the limit of its relative frequency. Demanding that this limit exists is non-trivial, since 
this need not be the case for an arbitrary sequence. Intuitively, \ref{vm1} 
expresses the hypothesis of statistical stability, which captures a global 
regularity of the sequence.

In contrast, \ref{vm2} captures a sense of \textit{randomness} or local irregularity.  
Note it actually comprises \emph{two claims}: 1) the limit exists and 2) it is the same as 
the limit in \ref{vm1}. It is best understood by viewing a selection 
rule $\vv{S}_j \colon \mathbb{N} \rightarrow \{0,1\}$ as selecting a 
subsequence of the original sequence $\vv{\Omega}$ and then demanding that the 
limiting relative frequencies thereof coincide with those of the original sequence. 
Such a selection rule is called \Def{\emph{admissible}}, whereas a selection 
rule  which would give rise to different limiting relative frequencies for at 
least one $A \in \mathcal{A}$ would be \Def{\emph{inadmissible}}. Why do we 
need axiom \ref{vm2}? Von Mises calls this the ``law of the excluded 
gambling system'' and it is the key to capture the notion of 
randomness in his framework. Intuitively, if a selection rule is inadmissible, 
an adversary could use this knowledge to strategically offer a bet on 
the next outcome and thereby make long-run profit, at the expense of 
our fictional decision maker. A \emph{random} sequence, however, is 
one for which there does not exist such a betting strategy. It turns out,
that this statement cannot hold in its totality. A sequence cannot be 
random with respect to all selection rules except in trivial cases 
(\cf Kamke's critique of von Mises' notion of randomness, nicely 
summarized in \citep{lambalgen1987mises}). Thus, von Mises explicitly 
relativizes randomness with respect to a problem-specific set of 
selection rules \citep[p.\@\xspace 12]{mises1964mathematical}.\footnote{This class 
of selection rules necessarily must be specified in advance; confer \citep{shen2009}.
A prominent line 
of work aspires to fix the set of selection rules as all partially computable selection rules \citep{church1940concept}, 
but there is no compelling reason to elevate this to a universal choice;
\cf \citep{derr2022fairness} for an elaborated critique.} A sequence which 
forms a collective (``is random with respect to'') one set of selection 
rules, might not form a collective with respect to another set.

In our view, the role of the randomness axiom \ref{vm2} is similar to 
the role of more familiar randomness assumptions like the standard 
\textit{i.i.d.} assumption: to empower inference from finite data. 
In this work, however, we will be exclusively concerned with the 
idealized case of infinite data, since our focus is the axiom 
(or hypothesis) of statistical stability.

We are motivated by the following question. What happens to 
von Mises approach when axiom \ref{vm1} breaks down? That is, 
when relative frequencies of at least some events do not converge. 
Our answer leads to a confluence with a 
theory that is thoroughly grounded in the subjectivist camp: 
the theory of \textit{imprecise probability}. 
In summary, we establish a strictly frequentist theory of imprecise probability. 

\subsection{Imprecise Probability - The Subjectivist Perspective}
\label{sec:walleyintro}
We briefly introduce the \textit{prima facie} unrelated, subjectivist theory of 
imprecise probability, or more specifically, the theory of \textit{lower and upper 
previsions} as put forward by \citet{walley1991statistical}. Orthodox Bayesianism 
models belief via the assignment of precise probabilities to propositions, 
or equivalently, via a linear expectation functional. 
In contrast, in Walley's theory,  belief is interval-valued and the linear expectation is replaced by a pair 
of a lower and upper expectation. Hence, the theory is 
strictly more expressive than orthodox Bayesianism, which can be recovered as a special case.

We assume an underlying possibility set $\Omega$, where $\omega \in \Omega$
is an elementary event, 
which includes all relevant information. 
We call a function $X \colon \Omega \rightarrow \mR$, which 
is bounded, \ie $\sup_{\omega \in \Omega} |X(\omega)| < \infty$, 
a \Def{\emph{gamble}} and collect all 
such functions in the set $\linfty$. The set of gambles 
$\linfty$ carries a vector space structure with scalar 
multiplication $(\lambda X)(\omega) = \lambda X(\omega)$, $\lambda \in \mR$, and addition 
$(X+Y)(\omega) = X(\omega) + Y(\omega)$. For a constant gamble $c(\omega)=c \text{ } \forall \omega$ 
we write simply $c$. Note that Walley's theory in the general case does not require that a vector space 
of gambles is given, but definitions and results simplify significantly in this case.

We interpret a gamble as assigning an uncertain loss $X(\omega)$ to 
each elementary event, that is, in line with the convention 
in insurance and machine learning, we take positive values to represent 
loss and negative values to represent reward.\footnote{Unfortunately,
this introduces tedious sign flips when comparing results 
to \citet{walley1991statistical}.} We imagine a decision maker who is faced 
with the question of how to value a gamble $X$; 
the orthodox answer would be the expectation $\E[X]$ 
with respect to a subjective probability measure.

\citet{walley1991statistical} proposed a betting interpretation of imprecise probability, which is inspired by \citet{de2017theory}, who identifies probability with fair betting rates. The goal is to axiomatize a functional $\Rup \colon \linfty \rightarrow \mR$, which assigns to a gamble the smallest number $\Rup(X)$ so that $X-\Rup(X)$ is a \text{desirable} transaction to our decision maker, where she incurs the loss $X(\omega)$ but in exchange gets the reward $-\Rup(X)$. Formally:
\[
    \Def{\Rup(X) \coloneqq \inf\{\alpha \in \mR \colon X - \alpha \in \mathcal{D} \}}, 
\]
where $\mathcal{D}$ is a \textit{set of desirable gambles}. \citet[Section 2.5]{walley1991statistical} argued for a criterion of coherence, which any reasonable functional $\Rup$ should satisfy, and consequently obtained the following characterization \cite[Theorem 2.5.5]{walley1991statistical}, which we shall take here as an axiomatic \textit{definition} instead.\footnote{Here, we need the vector space assumption on the set of gambles. We also note that \citet[pp.\@\xspace 64--65]{walley1991statistical} himself made a similar definition, but then proposed the more general coherence concept.} 
\begin{definition}
A functional $\Rup \colon \linfty \rightarrow \mR$ is a 
\Def{\emph{coherent upper prevision}} if it satisfies $\forall X,Y \in \linfty$:
\begin{enumerate}[nolistsep,label=\textbf{\emph{UP\arabic*.}}, ref=UP\arabic*]
  \item \label{item:UP1} $\Rup(X) \leq \sup(X)$ \quad \:\:\:\:\:\:\:\:\:\:\:\:\:\:\:\:\:\:\:\:  (bounds)
  \item \label{item:UP2} $\Rup(\lambda X) = \lambda \Rup(X) \text{, } \forall \lambda \in \mR^+$ \quad (positive homogeneity)
  \item \label{item:UP3} $\Rup(X+Y) \leq \Rup(X) + \Rup(Y)$ \quad  \thinspace\thinspace (subadditivity)
\end{enumerate}
\end{definition}
Together, these properties also imply $\forall X,Y \in \linfty$ \citep[p.\@\xspace 76]{walley1991statistical}:
\begin{enumerate}[nolistsep,start=4,label=\textbf{UP\arabic*.}, ref=UP\arabic*]
  \item \label{item:P6} $\Rup(X + c) = \Rup(X) + c \text{, } \forall c \in \mR$  \qquad \:\:\:\:\:\:\:\:\:\:\:\: \thinspace \thinspace (translation equivariance)
  \item \label{item:P7} $X(\omega) \leq Y(\omega) \text{ } \forall \omega \in \Omega \Rightarrow \Rup(X) \leq \Rup(Y)$ \quad (monotonicity)
\end{enumerate}
To a coherent upper prevision, we can define its conjugate \Def{\emph{lower prevision}} by:
\begin{align*}
    \Def{\Rlow(X)} &\Def{\coloneqq -\Rup(-X)}\\
    &= -\inf\{\alpha \in \mR\colon -X - \alpha \in \mathcal{D} \}\\
    &= \sup\{\alpha \in \mR\colon \alpha - X \in \mathcal{D} \},
\end{align*}
which specifies the highest certain loss $\alpha$ that the decision maker is willing to shoulder in exchange for giving away the loss $X(\omega)$, \ie receiving the reward $-X(\omega)$. Due to the conjugacy, it suffices to focus on the upper prevision throughout. In general, we have that $\Rlow(X)\leq \Rup(X)$ for any $X \in \linfty$. If $\Rlow(X)=\Rup(X)$ $\forall X \in \linfty$, we say that $R\coloneqq \Rup = \Rlow$ is a \textit{linear prevision}, a definition which aligns with \citet{de2017theory}.

By applying an upper prevision to indicator gambles, 
we obtain an \Def{\emph{upper probability}} $\Def{\Pup(A) \coloneqq \Rup(\chi_A)}$, 
where $A \subseteq \Omega$. Correspondingly, the \Def{\emph{lower probability}} is 
$\Def{\Plow(A) \coloneqq 1-\Pup(A^\comp) = \Rlow(\chi_A)}$. In the precise case, there is a 
unique relationship between (finitely) additive probabilities and linear previsions; however, 
upper previsions are more expressive than upper probabilities.
Finally, we remark that via the so-called \Def{\emph{natural extension}}, 
a coherent upper probability which 
is defined on some subsets of events $\mathcal{A} \subseteq 2^\Omega$ can be extended to a 
coherent upper prevision $\Def{\natext(\Pup)}$ on $\linfty$,  which is compatible with $\Pup$ 
in the sense that $\natext(\Pup)(\chi_A) = \Pup(A)$ $\forall A \in \mathcal{A}$ (\cf \citep[Section 3.1]{walley1991statistical}). 


\section{Unstable Relative Frequencies}
Assume that we have some fixed sequence 
$\vv{\Omega} \colon \mathbb{N} \rightarrow \Omega$ on a possibility set
$\Omega$ of elementary events, but that for some events 
$A \in \mathcal{A}$, where $\mathcal{A} \subseteq 2^\Omega$,
the limiting relative frequencies do not exist. What can we do then?
In a series of papers \citep{ivanenko1986functional,ivanenkoclassofcriterion, 
ivanenkomodelofnonstochasticrussian, ivanenko1993, Ivanenko2000DecisionMI, ivanenkoonregularities,
ivanenko2017expected} and a monograph \citep{ivanenkobook}, Ivanenko and 
collaborators have developed a strictly frequentist theory of ``hyper-random phenomena'' 
based on ``statistical regularities''. In essence, they tackle mass decision 
making in the context of sequences with possibly divergent relative frequencies.
Like von Mises, they take the notion of a sequence as the primitive, 
that is, without assuming an a priori probability and then invoking the 
law of large numbers. They explicitly recognize that ``stochasticness 
%
%
gets broken as soon 
as we deal with deliberate activity of people'' \citep{ivanenko1993}\footnote{This 
not only  occurs because of non-equilibrium effects, but also from feedback loops, 
what has become known as  ``performativity'' \citep{mackenzie2007economists}
or ``reflexivity'' \citep{soros2009}. See the epigraph at the beginning of the 
present paper.}.
The presentation of Ivanenko's theory is 
obscured somewhat by the great generality with which it is presented (they work 
with general nets, rather than just sequences). We build heavily upon their 
work but entirely restrict ourselves to working with sequences. While in some sense this
is a weakening, our converse result (see Section~\ref{sec:converseresult} is actually stronger as 
we show that one can achieve any 
``statistical regularity'' by taking relative frequencies of only sequences.
For simplicity, we will dispense with integrals with respect to  finitely additive 
measures in our presentation, so that there are less mathematical 
dependencies involved;\footnote{The two well-known accounts of the 
theory of integrals with finitely additive measures are 
\citep{rao1983theory} and \citep{dunford1988linear}. 
The theory of linear previsions as laid out in \citep{walley1991statistical} appears
to be an easier approach for our purposes.} instead, we work with linear 
previsions. Moreover, we establish\footnote{\citet{ivanenkoonregularities} 
mention in passing that sets of probabilities also appear 
in \citep{walley1991statistical}.} the connection to 
imprecise probability and give a different justification for the construction.

\subsection{Ivanenko's Argument --- Informally}
We begin by providing an informal summary of Ivanenko's 
construction of \textit{statistical regularities} on sequences. 
Assume we are given a fixed sequence 
$\vv{\Omega}\colon \mathbb{N} \rightarrow \Omega$ of elementary events 
$\vv{\Omega}(1),\vv{\Omega}(2),\ldots$, where we may 
intuitively think of $\mathbb{N}$ as representing time.
In contrast to von Mises, who demands the existence of relative 
frequency limits to define probabilities, we ask for something like 
a probability for \emph{all} events $A \subseteq \Omega$,
even when the relative frequencies have no limit. 
To this end, we exploit that sequences of relative frequencies 
always have a non-empty set of cluster points, 
each of which is a finitely additive probability.
Hence, a decision maker can use this set of probabilities to represent 
the global statistical properties of the sequence. 
In fact, we will see that this induces a coherent upper probability. 
Also, our decision maker is interested in assessing a value for each 
gamble $X\colon \Omega \rightarrow \mR$, which is evaluated infinitely 
often over time. Here, the sequence of averages 
$n \mapsto \frac{1}{n} \sum_{i=1}^n X(\vvomi)$ is the 
object of interest. In the case of convergent relative frequencies,
a decision maker would use the expectation to assess 
the risk in the limit, whereas in the general case of 
possible non-convergence, a different object is needed. 
This object turns out to be a coherent upper prevision. 
We provide a justification for using this upper prevision to assess 
the value of a gamble, which links it to imprecise probability.

\subsection{Ivanenko's Argument --- Formally}
\label{sec:ivanenkoformal}

Let $\Def{\Omega}$ be an arbitrary (finite, countably infinite or uncountably infinite) set of outcomes and fix $\Def{\vv{\Omega} \colon \mathbb{N} \rightarrow \Omega}$, an $\Omega$-valued sequence. 
We define a \Def{gamble $X\colon \Omega \rightarrow \mR$} as a bounded function from $\Omega$ to $\mR$, \ie $\exists K \in \mR\colon |X(\omega)| \leq K$ $\forall \omega \in \Omega$ and collect all such gambles in the set $\linfty$. We assume the vector space structure on $\linfty$ as in Section~\ref{sec:walleyintro}.

The set $\linfty$ becomes a Banach space, \ie a complete normed vector space, under the 
supremum norm $\Def{\|X\|_{\linfty} \coloneqq \sup_{\omega \in \Omega} |X(\omega)|}$. We denote the topological 
dual space of $\linfty$ by $\Def{\linftyd}$. Recall that it consists exactly of the continuous 
linear functionals $\phi\colon \linfty \rightarrow \mR$. We endow $\linftyd$ with the \emph{weak*-topology},
which is the weakest topology (\ie with the fewest open sets) that makes all evaluation functionals 
of the form $\Def{X^* \colon \linftyd \rightarrow \mR}$, $\Def{X^*(E) \coloneqq E(X)}$ for any 
$X \in \linfty$ and $E \in \linftyd$ continuous. Consider the following subset of $\linftyd$:
\[
    \Def{\PF(\Omega) \coloneqq \left\{E \in \linftyd \colon E(X) \geq 0 \text{ whenever } X \geq 0, E(\chi_\Omega)=1\right\}.}
\]
Due to the Alaoglu-Bourbaki theorem, this set is compact under the weak* topology, see Appendix~\ref{app:weakstarcompact}.

A finitely additive probability $P\colon \mathcal{A} \rightarrow [0,1]$ on some set system $\mathcal{A}$, where $\Omega \in \mathcal{A}$, is a function such that:
\begin{enumerate}[nolistsep,label=\textbf{PF\arabic*.}, ref=PF\arabic*]
   \item \label{prop:pf1} $P(\Omega)=1$.
   \item \label{prop:pf2} $P(A \cup B)=P(A)+P(B)$ whenever $A \cap B= \emptyset$ and $A,B \in \mathcal{A}$.
\end{enumerate}

We induce a sequence of finitely additive probabilities $\vv{P}$ where
$\Def{\vv{P}(n) \coloneqq A \mapsto \frac{1}{n} \sum_{i=1}^n
\chi_A(\vv{\Omega}(i))}$
for each $n \in \mathbb{N}$. It is easy to 
check that indeed $\vv{P}(n)$ is a finitely additive probability on the whole powerset $2^\Omega$ 
for any $n \in \mathbb{N}$. We shall call $\vv{P}$ the \Def{\emph{sequence of empirical probabilities}}.
Due to \citep[Corollary 3.2.3]{walley1991statistical}, a finitely additive probability defined on 
$2^\Omega$ can be uniquely extended (via natural extension) to a linear prevision 
$E_P \colon \linfty \rightarrow \mR$, so that $E_P(\chi_A) = P(A)$ $\forall A \subseteq \Omega$. 
Furthermore, we know from \citep[Corollary 2.8.5]{walley1991statistical}, that there is 
a one-to-one correspondence between elements of $\PF(\Omega)$ and linear previsions 
$E_P \colon \linfty \rightarrow \mR$. Hence, we associate to each empirical probability 
$\vv{P}(n)$ an \Def{\emph{empirical linear prevision}} 
$\Def{\vv{E}(n) \coloneqq X \mapsto \natext(\vv{P}(n))(X)}$, where $X \in \linfty$ and we 
denote the natural extension by $\natext$. We thus obtain a sequence 
$\vv{E} \colon \mathbb{N} \rightarrow \PF(\Omega)$.

On the other hand, each gamble $X \in \linfty$ induces a sequence of evaluations 
as  $\Def{\vv{X}\colon \mathbb{N} \rightarrow \mR}$, where $\Def{\vv{X}(n) \coloneqq X\left(\vv{\Omega}(n)\right)}$.
For $X \in \linfty$, we define the sequence of averages of the gamble over time as 
$\Def{\avgseq \colon \mathbb{N} \rightarrow \mR}$, where 
$\Def{\avgseq(n) \coloneqq \frac{1}{n} \sum_{i=1}^n X\left(\vv{\Omega}(i)\right)}$. For each finite $n$, we can 
also view the average as a function in $X$, \ie $X \mapsto \frac{1}{n} \sum_{i=1}^n X\left(\vv{\Omega}(i)\right)$. 
Observe that this is a coherent linear prevision and by applying it to indicator gambles $\chi_A$, 
we obtain $\vv{P}(n)$. Hence, we know from \citep[Corollary 3.2.3]{walley1991statistical} that this linear 
prevision is in fact the natural extension of $\vv{P}(n)$, \ie 
$\vv{E}(n) = X \mapsto \frac{1}{n} \sum_{i=1}^n X\left(\vv{\Omega}(i)\right) = X \mapsto \avgseq(n)$.
This concludes the technical setup; we now begin reproducing Ivanenko's argument.

Since $\PF(\Omega)$ is a compact topological space under the subspace topology induced by the weak*-topology 
on $\linftyd$, we know that any sequence $\vv{E}\colon \mathbb{N} \rightarrow \linftyd$ has a 
non-empty closed set of cluster points. Recall that a point $z$ is a \Def{\emph{cluster point}} of 
a sequence $\vv{S} \colon \mathbb{N} \rightarrow \mathcal{T}$, where $\mathcal{T}$ is any topological space, if:
\[
    \forall N \mbox{, where $N$ is any neighbourhood of } z \text{ with respect to } \mathcal{T} \text{, } \forall n_0 \in \mathbb{N}: \exists n \geq n_0: \vv{S}(n) \in N.
\]
We remark that this \emph{does not} imply that those cluster points are limits of 
convergent subsequences.\footnote{This would hold under sequential compactness, which 
is not fulfilled here in general, but it is for finite $\Omega$.} 
We denote the \Def{\emph{set of cluster points}} as $\Def{\CP(\vv{E})}$. 
Equivalently, by applying these linear previsions to indicator gambles, we obtain the \Def{set of 
finitely additive probabilities} $\Def{\mathcal{P} \coloneqq \left\{A \mapsto E(\chi_A) \colon E \in \CP(\vv{E})\right\}}$.
Due to the one-to-one relationship, we might work with either $\CP(\vv{E})$ or $\mathcal{P}$. 
Following Ivanenko, we call $\mathcal{P}$ the \Def{\emph{statistical regularity}} of the sequence 
$\vv{\Omega}$; in the language of imprecise probability,  it is a \Def{\emph{credal set}}.

We further define 
\[
            \Def{\Rup(X) \coloneqq  \sup \left\{E(X) \colon E \in \CP(\vv{E}) \right\} =
            \sup \left\{E_P(X) \colon P \in \mathcal{P} \right\} , \quad X \in \linfty,}
 \]
where $E_P \coloneqq \natext(P)$, and
\[
\Def{\Pup(A) \coloneqq \sup \left\{E(\chi_A) \colon E \in \CP(\vv{E}) \right\} =
\sup \left\{ P(A) \colon P \in \mathcal{P} \right\}, \quad A \subseteq \Omega.}
\]
Observe that $\Rup$ is defined on all $X \in \linfty$ and $\Pup$ is defined 
on \textit{all} subsets of $\Omega$, even if $\Omega$ is uncountably infinite, 
since each $P \in \mathcal{P}$ is a finitely additive probability on $2^\Omega$.
We further observe that $\Rup$ is a coherent upper prevision on 
$\linfty$ or equivalently, a coherent risk measure in the sense of 
\citet{artzner1999coherent}.\footnote{For the close connection of coherent 
upper previsions and coherent risk measures we refer to 
\citep{frohlich2022risk}.} Correspondingly, $\Pup$ is a coherent upper 
probability on $2^\Omega$, which is obtained by applying $\Rup$
to indicator functions. This follows directly from the 
\textit{envelope theorem} in \citep[Theorem 3.3.3]{walley1991statistical}.

So far, the definition of $\Rup$ and $\Pup$ may seem unmotivated. Yet they play a special role, as we now show.

\begin{proposition}
The sequence of averages $\avgseq$ has the set of cluster points 
\[
    \operatorname{CP}\left(\avgseq\right) = \left\{E(X) \colon E \in \CP(\vv{E})\right\} =
    \left\{E_P(X) \colon {P} \in \mathcal{P}\right\},
\]
and therefore
\[
\Rup(X) = \sup \CP\left(\avgseq\right) = \limsupn \avgseq(n).
\]
\end{proposition}
\begin{proof}
First observe that
\[
    \vv{E}(n)(X) = \avgseq(n).
\]

We use the following result from \citep[Lemma 3]{ivanenko2017expected}.\footnote{A subtle point in the argument, which \citet{ivanenko2017expected} do not make visible, is the sequential compactness of $\mR$, which means that for any cluster point of an $\mR$-valued sequence we can find a subsequence converging to it. }
\begin{lemma}
Let $f\colon Y \rightarrow \mR$ be a continuous function on a compact space $Y$ and $\vv{y}$ a $Y$-valued sequence. Then $\CP\left(n \mapsto f\left(\vv{y}(n)\right)\right) 
= f\left(\CP\left(\vv{y}\right)\right)$.  \end{lemma}
On the right side, the application of $f$ is to be understood as applying $f$ to each element in the set $\CP\left(\vv{y}\right)$.
Consider now the evaluation functional $X^* \colon \PF(\Omega) \rightarrow \mR$, $X^*(E) \coloneqq E(X)$, which is continuous under the weak*-topology. The application of the lemma with $f=X^*$, $Y=\PF(\Omega)$,$\vv{y}=\vv{E}$ hence gives:
\[
    \CP\left(n \mapsto X^*\left(\vv{E}(n)\right)\right) = X^*\left(\CP\left(\vv{E}\right)\right).
\]
But since $X^*\left(\vv{E}(n)\right)=\avgseq(n)$, we obtain that $\CP\left(\avgseq\right) = 
X^*\left(\CP\left(\vv{E}\right)\right) = \left\{E(X) \colon E \in \CP\left(\vv{E}\right)\right\}$.
\end{proof}

A similar statement holds for the coherent upper probability.
\begin{corollary}
For any $A\subseteq \Omega$,  $\displaystyle\Pup(A) = \limsupn \left(\vv{P}(n)(A)\right) = \limsupn \frac{1}{n} \sum_{i=1}^n \chi_A\left(\vvomi\right)$.
\end{corollary}
\begin{proof}
Just observe that
$\vv{P}(n)(A) = \vv{\Sigma \chi_A}(n)$ and apply the previous result.
\end{proof}

Thus the limes superior of the sequence of relative frequencies induces a 
coherent upper probability on $2^\Omega$; similarly, the limes superior of 
the sequence of a gamble's averages induces a coherent upper prevision on 
$\linfty$. By conjugacy, we have that the lower prevision and lower probability are
:
\begin{align*}
    \Rlow(X) = &\inf \left\{E(X) \colon E \in \CP\left(\vv{E}\right) \right\}
           = \liminf_{n\rightarrow\infty} \avgseq(n), \quad \forall X \in \linfty, \\
    \Plow(A) = &\inf \left\{ P(A) \colon P \in \mathcal{P} \right\} 
           =\liminf_{n\rightarrow\infty} \frac{1}{n} \sum_{i=1}^n 
           \chi_A\left(\vvomi\right), \quad A \subseteq \Omega,
\end{align*}
which are obtained in a similar way using the limes inferior. 
Finally, when an event is \textit{precise} in the sense that 
$\Pup(A)=\Plow(A)$ (and thus the $\liminf$ equals the $\limsup$ 
and hence the limit exists), we denote the upper (lower) probability as 
$P(A)$ and say that the precise probability of $A$ exists.

\subsection{The Induced Upper Prevision}
\label{sec:inducedriskmeasure}
We have seen that the upper prevision $\Rup$, as we have just defined it, 
has the property that it is induced by the statistical regularity of 
the sequence, and at the same time corresponds to taking the supremum 
over the cluster points of the sequence of averages of a gamble over 
time. The set of cluster points is in general a complicated object, 
hence it is unclear why one should take the supremum to reduce it 
to a single number in a decision making context. Our goal in this 
section is to provide some intuition \textit{why} it is reasonable to 
use $\Rup$, as we  defined it, to assess the risk inherent in a sequence.
\citet{Ivanenko2000DecisionMI} argued that $\Rup$ is the unique 
object which satisfies a certain list of axioms, which are similar to 
those for an upper prevision, but including a so-called 
``principle of guaranteed result'', which appears rather mysterious to us.

Our setup is as follows. We imagine an individual decision maker, 
who is faced with a fixed sequence $\vv{\Omega} \colon \mathbb{N}
\rightarrow \Omega$ and various gambles $X\colon \Omega 
\rightarrow \mR$. The question to the decision maker is how to 
value this gamble in light of the sequence, \ie imagining that the 
gamble is evaluated at each $\vv{\Omega}(1),\vv{\Omega}(2),...$, 
infinitely often. Here, $X(\vv{\Omega}(i))$ represents a loss 
for positive values, and a gain for negative values. 
We can think of the $\vv{\Omega}(i)$ as the states of nature, 
and the sequence determines which are realized and how often.
Importantly, we view our decision maker as facing a 
\textit{mass decision}, \ie the gamble will not only be 
evaluated once, but instead infinitely often.

Which gambles are \textit{desirable} to our decision maker? Desirable gambles\footnote{We remark that throughout the paper we always mean \textit{almost desirability}, \cf \citep[Section 3.7]{walley1991statistical}.} are those which the decision maker would accept even without any reward for it. We argue that an appropriate set of desirable gambles is given by:
\[
    \Def{\mathcal{D}_{\vv{\Omega}} \coloneqq \left\{X \in \linfty \colon 
       \limsup_{n \rightarrow \infty}  \avgseq \leq 0\right\}.}
\]
Consider what $X \in \mathcal{D}_{\vv{\Omega}}$ means. If the limes superior of the gamble's average sequence, \ie the growth rate of the accumulated loss, is negative or zero, then we are guaranteed that there is no strictly positive accumulated loss which we will face infinitely often. The choice of the average as the aggregation functional is justified from the mass decision character of the setting, since we assume that our decision maker does not care about individual outcomes, but only about long-run capital. Now, given this set of desirable gambles, we seek a functional $\Rup(X) \colon \linfty \rightarrow \mR$, so that when at each time step $i$, the transaction $X(\vvomi) - \Rup(X)$ takes place, this results in a desirable gamble for our decision maker. Our decision maker shoulders the loss $X(\vvomi)$, while at the same time asking for $-\Rup(X)$ in advance. Intuitively, $\Rup(X)$ is supposed to be the certainty equivalent of the ``uncertain'' loss $X$, in the sense that $X(\vvomi)$ will vary over time. Therefore we define, in correspondence with $\mathcal{D}_{\vv{\Omega}}$, the upper and lower previsions ($\forall X \in \linfty$): 
\begin{align}
    \label{eq:rdeffromdesirable}
    \Def{\Rup(X) \coloneqq} &\Def{\inf\left\{\alpha \in \mR\colon X-\alpha \in \mathcal{D}_{\vv{\Omega}}\right\}}\\
    \Def{\Rlow(X) \coloneqq -\Rup(-X) =} &\Def{\sup\left\{\alpha \in \mR\colon \alpha - X \in 
      \mathcal{D}_{\vv{\Omega}}\right\}.}\nonumber
\end{align}
When a set of desirable gambles and an upper (lower) prevision are in this correspondence, it holds that $X-\Rup(X) \in \mathcal{D}_{\vv{\Omega}}$ and furthermore $\Rup(X)$ is the least (smallest) functional for which this holds. We can now observe that in fact $\Rup(X) = \limsupn \avgseq(n)$, since
\[
        \mathcal{D}_{\vv{\Omega}} =\left\{X \in \linfty \colon \Rup(X) \leq 0\right\}
\]
is the general correspondence for a set of desirable gambles and a coherent upper prevision. Hence, we have motivated the definition of $\Rup(X)$ in Section~\ref{sec:ivanenkoformal}. It is easy to see explicitly (\cf Appendix~\ref{app:inducedriskmeasure}) that $X-\Rup(X) \in \mathcal{D}_{\vv{\Omega}}$ and that $\Rup(X) = \limsupn \avgseq(n)$ is in fact the smallest number such that the relation in Equation~\ref{eq:rdeffromdesirable} holds.

\section{From Cluster Points to Sequence}
\label{sec:converseresult}




In the previous section, we have shown how from a given sequence we can construct a coherent upper prevision from the set of cluster points $\CP(\vv{E})$.
 In this section, we show the converse, thus ``closing the loop'': given an arbitrary coherent upper prevision, we construct a sequence $\vv{\Omega}$ such that the induced upper prevision is just the specified one. We take this to be an argument for the well-groundedness of our approach. For simplicity, we assume a finite possibility space $\Omega$.

 \begin{theorem}
 \label{theorem:converse}
 Let $|\Omega|<\infty$. Let $\Rup$ be a coherent upper prevision on $\linfty$. There exists a sequence $\vv{\Omega}$ such that we can write $\Rup$ as:
\[
            \Rup(X) = \Rup_{\vv{\Omega}}(X) = \sup \left\{E(X) \colon E \in \mathcal{E}_{\vv{\Omega}} \right\}, \quad \mathcal{E}_{\vv{\Omega}} \coloneqq \CP\left(\vv{E}_{\vv{\Omega}}\right) \quad \forall X \in \linfty,
\]
where we now make the dependence on the sequence $\vv{\Omega}$ explicit in the notation, \ie $\vv{E}_{\vv{\Omega}}(n) = X \mapsto \frac{1}{n} \sum_{i=1}^n X\left(\vv{\Omega}(i)\right)$.
 \end{theorem}
 
The significance of this result is that it establishes strictly frequentist \textit{semantics} for imprecise probability. It shows that to any decision maker who, in the subjectivist fashion, uses a coherent upper prevision, we can associate a sequence, which would yield the same upper prevision in a strictly frequentist way. We interpret this result as evidence for the naturalness, and arguably completeness, of our theory.

The key to prove this is Theorem~\ref{th:CP-r-C}, for which we introduce some convenient notation.
 For $k\in\naturals$, let $\Def{[k]\coloneqq\{1,\ldots,k\}}$ and
define the \Def{$(k-1)$-simplex} as
\[
    \Def{\Delta^k \coloneqq \left\{d=(d_1,\ldots,d_k) \in \mR^k\colon
    \sum_{i=1}^n d_i = 1,\  d_i \geq 0 \  \forall i\in[k]\right\}.}
\]

It is also helpful to have a dual notation for sequences $z\colon\naturals\rightarrow
[k]$, whereby we write either $\Def{z(i)}$ or $\Def{z_i}$ to mean the same thing.
\begin{definition}
	Suppose $k\in\naturals$ and $x\colon\naturals\rightarrow [k]$. For
	any $n\in\naturals$ define the \Def{\emph{relative frequency of 
	$x$ with respect to $i$ at $n$}}, 
	$r_i^x\colon\naturals\rightarrow [0,1]$ via
	\[
		\Def{r_i^x(n)\coloneqq \textstyle\frac{1}{n}|\{j\in[n]\colon
		x_j=i\}|}
	\]
	and the \Def{\emph{relative frequency  of $x$ at $n$}},
	$r^x\colon\naturals\rightarrow\Delta^k$ as
	\begin{equation}
		  \Def{r^x(n)\coloneqq r_{[k]}^x(n)=
		       \left(\begin{array}{c}r_1^x(n)\\ \vdots
			\\ r_k^x(n)\end{array}\right).}
		\label{eq:r-x-def}
	\end{equation}
\end{definition}

\begin{theorem}\label{th:CP-r-C}
	Suppose $k\in\naturals$ and $C$ is a rectifiable closed curve in $\Delta^k$. There
	exists $x\colon\naturals\rightarrow [k]$ such that $\CP(r^x)=C$.
\end{theorem}
The proof (which is constructive) is in Appendix~\ref{app:construction} along
with an example.  From this, we obtain the following Corollary (proven in Appendix~\ref{app:fromboundarytocurves}). Denote the topological boundary of a set $D$ as $\partial D$.

\begin{corollary}
\label{cor:CP-convexbody}
Suppose $k \in \naturals$ and $D \subseteq \Delta^k$ is a non-empty convex set. There exists $x\colon\naturals\rightarrow [k]$ such that $\CP(r^x)=\partial D$.
\end{corollary}



Since we have a finite possibility space $\Omega$, we can identify each linear prevision $\mathcal{E}$ with a point in the simplex, by assigning coordinates to its underlying finitely additive probability; in the case of $\vv{E}_{\vv{\Omega}}(n)$, this is the relative frequency $r^{\vv{\Omega}}(n)$. This is formalized in the following.
\begin{proposition}
\label{prop:simplexcorrespondence}
Let $\vv{E}(n) : \mathbb{N} \rightarrow \PF(\Omega)$ be a sequence of linear previsions with 
underlying probabilities $\vv{P}(n) \coloneqq A \mapsto \vv{E}(n)(A)$. 
Then $E \in \CP\left(\vv{E}(n)\right)$ with respect to the weak* topology if and only if 
the sequence $\vv{D} \colon \mathbb{N} \rightarrow \Delta^k$, $\vv{D}(n) \coloneqq 
\left(\vv{P}(n)(\omega_1),\ldots, \vv{P}(n)(\omega_k)\right)$ has as cluster 
point $d_E=\left(E\left(\chi_{\{\omega_1\}}\right),\ldots,E\left(\chi_{\{\omega_k\}}\right)\right)$
with respect to the topology induced by the Euclidean norm on $\mR^k$.
\end{proposition}

For the proof, see Appendix~\ref{app:euclideansufficient}. Combining Corollary~\ref{cor:CP-convexbody} and Proposition~\ref{prop:simplexcorrespondence} allows us to now prove Theorem~\ref{theorem:converse}.

\begin{proof}[Proof of Theorem~\ref{theorem:converse}]
Let $\Omega=[k]$. If $\Rup$ is a coherent upper prevision on $\linfty$, we can write it as \citep[Theorem 3.6.1]{walley1991statistical}:
\[
\Rup(X) = \sup \left\{E(X) \colon E \in \mathcal{E} \right\}, \quad \forall X \in \linfty,
\]
for some weak* compact and convex set $\mathcal{E} \subseteq \PF(\Omega)$. From \citep[Theorem 3.6.2]{walley1991statistical} we further know that
\[
\Rup(X) = \sup \left\{E(X) \colon E \in \mathcal{E} \right\} = \sup \left\{E(X) \colon E \in \operatorname{ext} \mathcal{E} \right\},\quad \forall X \in \linfty,
\]
where $\operatorname{ext}$ denotes the set of extreme points of $\mathcal{E}$.\footnote{A point $E \in \mathcal{E}$ is an extreme point of $\mathcal{E}$ if it cannot be written as a convex combination of any other elements in $\mathcal{E}$.} 
Then:
\begin{align*}
    \Rup(X) &= \sup \left\{E(X) \colon E \in \mathcal{E} \right\}\\
    &= \sup \left\{E(X) \colon E \in \operatorname{ext} \mathcal{E} \right\}\\
    &\leq \sup \left\{E(X) \colon E \in \partial \mathcal{E} \right\}\\
    &\leq \sup \left\{E(X) \colon E \in \mathcal{E} \right\} = \Rup(X),
\end{align*}
since $\operatorname{ext}\mathcal{E} \subseteq \partial \mathcal{E}$ and $\partial \mathcal{E} \subseteq \mathcal{E}$; note that $\mathcal{E}$ is closed. In summary, $\Rup(X) = \sup \left\{E(X) \colon E \in \partial \mathcal{E} \right\}$.

Now choose $D \coloneqq \left\{\left(E(\chi_{\omega_1}),\ldots,E(\chi_{\omega_k})\right) \colon E \in \mathcal{E}\right\}$, which is a non-empty convex set in $\Delta^k$. We then obtain from Corollary~\ref{cor:CP-convexbody} a sequence $\vv{\Omega} \colon \naturals \rightarrow [k]$ with $\CP(r^{\vv{\Omega}})=\partial D$. But then it follows from Proposition~\ref{prop:simplexcorrespondence} that the sequence $\vv{E}_{\vv{\Omega}}$ has cluster points $\CP\left(\vv{E}_{\vv{\Omega}}\right) = \partial \mathcal{E}$. Thus 
\[
\Rup(X) = \sup \left\{E(X) \colon E \in \CP\left(\vv{E}_{\vv{\Omega}}\right) \right\}, \quad \forall X \in \linfty,
\]
which concludes the proof.
\end{proof}

\citet{ivanenkobook} offers a somewhat similar result to Theorem~\ref{theorem:converse} by generalizing from sequences to \textit{sampling nets}. Ivanenko's \citeyearpar{ivanenkobook} main result states that ``any sampling directedness has a regularity, and any regularity is
the regularity of some sampling directedness.'' \citep[Theorem 4.2]{ivanenkobook}. We provide a brief introduction to Ivanenko's setup in Appendix~\ref{app:ivanenkonets}. Our result is more parsimonious in the sense that it relies only on sequences, which are arguably more intuitive objects than such sampling nets.

Our result should also be compared to Theorem 4.2 in  \citep{walley1982towards} and Theorem 2.2 in \citep{papamarcou1991unstable}. On the one hand, our result is stronger since it holds for upper previsions, whereas Theorem 4.2 in \citep{walley1982towards} and Theorem 2.2 in \citep{papamarcou1991unstable} hold for upper probabilities only; note that upper previsions are more expressive than upper probabilities.\footnote{Indeed, the proof of Theorem 4.2 in \citep{walley1982towards} exploits this simplification by assuming that the credal set has a finite number of extreme points.} On the other hand, Theorem 2.2 in \citep{papamarcou1991unstable} is stronger in the sense that it guarantees that the same upper probability is induced when applying selection rules.

We observe that two sequences $\vv{\Omega}_1$, $\vv{\Omega}_2$, might have different sets of cluster points $\CP\left(\vv{E}_{\vv{\Omega}_1}\right)$, $\CP\left(\vv{E}_{\vv{\Omega}_2}\right)$, but when their convex hull coincides, the same upper probability and prevision is induced.\footnote{Assume $\Rup(X) \coloneqq \sup\left\{E(X) \colon E \in \mathcal{E}\right\}$. Then indeed $\Rup(X) = \sup\left\{E(X) \colon E \in \coclosed \mathcal{E}\right\}$, where $\coclosed$ denotes the weak* closure of the convex hull; \cf \citep[Section 3.6]{walley1991statistical}.} Thus, in light of the argument in Section~\ref{sec:inducedriskmeasure}, \textit{for the purpose of mass decision making}, we may consider these sequences equivalent. While in the classical case, relative frequencies are the relevant description of a sequence, the statistical regularity provides an analogous description in the general case; moreover, we differentiate only ``up to the same convex hull'' for decision making.


\section{Unstable Conditional Probability}
An interesting aspect of the strictly frequentist approach is that there is a natural way of introducing conditional probability for events $A,B \subseteq \Omega$, which is the same for the case of converging or diverging relative frequencies. Furthermore, this approach generalizes directly to gambles. We will observe that this, perhaps surprisingly, yields the \textit{generalized Bayes rule}. In the precise case, the standard Bayes rule is recovered.

Recall that for a countably or finitely additive probability $Q$, we can define conditional probability as:
\begin{equation}
\label{eq:defcondmeasure}
\Def{Q(A|B) \coloneqq \frac{Q(A\cap B)}{Q(B)}, \quad A,B \subseteq \Omega, \text{ if } Q(B)>0.}
\end{equation}
Important here is the condition that $Q(B)>0$. Conditioning on events of measure zero may create trouble. Kolmogorov then allows the conditional probability to be arbitrary. This is rather unfortunate, as there arguably are settings where one would like to condition on events of measure zero.

As a prerequisite, given a linear prevision $E \in \PF(\Omega)$, we define the conditional linear prevision as:
\begin{equation}
\label{def:condlinearprev}
    \Def{E(X|B) \coloneqq \frac{E\left(\chi_B X\right)}{E\left(\chi_B\right)}, \quad \text{ if } E\left(\chi_B\right)>0, B \subseteq \Omega.}
\end{equation}
The application to indicator gambles then recovers conditional probability. 
As long as $E\left(\chi_B\right)>0$, it is insignificant whether we condition 
the linear prevision, or instead condition on the level of its 
underlying probability and then naturally extend it;
confer  \citep[Corollary 3.2.3]{walley1991statistical}.

Nearly in line with Kolmogorov's conditional probability, von Mises started from the following intuitive, frequentist view: the probability of an event $A$ conditioned on an event $B$ is the frequency of the occurence of the event $A$ given that $B$ happens. In what follows, we build upon this idea, which von Mises called ``partition operation'' \citep[p.\@\xspace 22]{mises1964mathematical}.
\citet[Section 4.3]{walley1982towards} have extended this definition to the divergent case of conditional probability on a finite possibility space; we further extend it to conditional upper previsions on arbitrary possibility spaces and link them to the generalized Bayes rule.
As a technical preliminary, we define a wrapper function
$\Def{\Psi \colon \PF(\Omega)\cup\{\bot\} \rightarrow \PF(\Omega)}$ as:
\begin{equation*}
    \Def{\Psi(P) \coloneqq \begin{cases}
    P_0 & \text{ if } P=\bot,\\
    P &  \text{ otherwise},\\
    \end{cases} }
\end{equation*}
where $P_0$ is an arbitrary finitely additive probability on $2^\Omega$.

\subsection{Conditional Probability}
\label{sec:unstablecondprob}
Recall our sequence of unconditional finitely additive probabilities $\vv{P}(n) \coloneqq A \mapsto \frac{1}{n} \sum_{i=1}^n \chi_A\left(\vv{\Omega}(i)\right)$. We want to define a similar sequence of \textit{conditional} finitely additive probabilities. A very natural approach is the following:
let $A,B \subseteq \Omega$ be such that $\vvomi \in B$ for at least one $i \in \mathbb{N}$. 
We write $\Def{\eventonce}$ for the set of such events, \ie events which occur at least once in the sequence.
Define a sequence of conditional probabilities $\Def{\vv{P}(\cdot|B)\colon \mathbb{N} \rightarrow \PF(\Omega)}$ by 
\begin{equation}
\label{eq:defcondprobseq}
    \Def{\vv{P}(\cdot|B)(n) \coloneqq \Psi\left(A \mapsto \frac{\sum_{i=1}^n (\chi_{A} \cdot 
     \chi_B)\left(\vv{\Omega}(i)\right)}{\sum_{i=1}^n \chi_B\left(\vv{\Omega}(i)\right)}\right),}
\end{equation}
where we consider only those $\vvomi$ which lie in $B$, and hence we adapt the relative frequencies to the 
occurrence of $B$. Informally, this is simply counting $|\mbox{$A$ and $B$ occured}|/|\mbox{$B$ occured}|$. Until $B$ occurs for the first time, the denominator will be $0$ and thus the mapping undefined 
 (returning the falsum $\bot$).
Throughout, we demand that the event $B$ on which we condition is in $\eventonce$, \ie 
occurs \emph{at least once} in the sequence. Note that this is a much weaker condition than demanding 
that $P(B)>0$, if $B$ is precise. Denote by $\Def{n_B}$ the smallest index so that 
$\vv{\Omega}(n_B) \in B$. Note that $\vv{P}(A|B)(n) = \vv{P}(n)(A \cap B) / \vv{P}(n)(B)$ for $n \geq n_B$.

\begin{proposition}
$\vv{P}(\cdot|B)$ is a sequence of finitely additive probabilities.
\end{proposition}
\begin{proof}
For $n < n_B$, this is clear due to $\Psi$. Now let $n \geq n_B$.

\ref{prop:pf1}: $\vv{P}(\Omega|B)(n)=1$: obvious.

\ref{prop:pf2}: If $A,C \subseteq \Omega$, $A \cap C = \emptyset$, then we show that $\vv{P}(A \cup C|B)(n) = \vv{P}(A|B)(n) + \vv{P}(C|B)(n)$. 
\begin{align*}
\vv{P}(A \cup C|B)(n) &= \frac{\sum_{i=1}^n \left(\chi_{A\cup C}
\cdot \chi_B\right)\left(\vv{\Omega}(i)\right)}{\sum_{i=1}^n \chi_B\left(\vv{\Omega}(i)\right)} \\
&= \frac{\sum_{i=1}^n \left(\chi_{A} \cdot \chi_B\right)\left(\vv{\Omega}(i)\right) + 
   \sum_{i=1}^n \left(\chi_{C} \cdot \chi_B\right)\left(\vv{\Omega}(i)\right)}{\sum_{i=1}^n 
   \chi_B\left(\vv{\Omega}(i)\right)}\\
&= \vv{P}(n)(A\cap B) / \vv{P}(n)(B) + \vv{P}(n)(C\cap B) / \vv{P}(n)(B)\\
&= \vv{P}(A|B)(n) + \vv{P}(C|B)(n).
\end{align*}
Noting that since $A$ and $C$ are disjoint, $\vv{\Omega}(i)$ cannot lie in both at the same time for any $i$.
\end{proof}
Even though the probability is conditional, we deal with a sequence of finitely additive probabilities again.
Hence, we can now essentially repeat the argument from Section~\ref{sec:ivanenkoformal}.
To each $\vv{P}(\cdot|B)(n)$, associate its uniquely corresponding linear prevision $\Def{\vv{E}(\cdot|B)(n)}$,
which is of course given by ($\forall X \in \linfty$, $n \geq n_B$):
\[
    \Def{\vv{E}(\cdot|B)(n) = \avgseq|B(n) \coloneqq X \mapsto 
    \frac{\sum_{i=1}^n \left(X \cdot \chi_B\right)\left(\vv{\Omega}(i)\right)}{\sum_{i=1}^n \chi_B\left(\vv{\Omega}(i)\right)}.}
\]
It is easy to check that $\vv{E}(\cdot|B)(n)$ is coherent. For $n<n_B$, set 
$\Def{\vv{E}(\cdot|B)(n) = \natext(P_0)}$.
From the weak* compactness of $\PF(\Omega)$, we obtain a non-empty closed set of cluster points 
$\CP(\vv{E}(\cdot|B))$.


\begin{definition}
If $B \in \eventonce$, we define the conditional upper prevision and the conditional upper probability as:
\[
\Def{ \Rup(X|B) \coloneqq \sup\left\{\tilde{E}(X) \colon \tilde{E} \in \CP\left(\vv{E}(\cdot|B)\right)\right\};
\quad \Pup(A|B) \coloneqq \sup\left\{Q(A) \colon Q \in \CP\left(\vv{P}(\cdot|B)\right)\right\}, \quad A \subseteq \Omega.}
\]
\end{definition}
Since they are expressed via an envelope representation,\footnote{An envelope representation expresses a coherent upper prevision as a supremum over a set of linear previsions.} $\Rup$ and $\Pup$ are automatically coherent \citep[Theorem 3.3.3]{walley1991statistical}. By similar reasoning as in Section~\ref{sec:ivanenkoformal}, we get the following representation.

\begin{proposition}
The conditional upper prevision (probability) can be represented as:
\[
\Rup(X|B) = \limsupn \avgseq|B(n), \quad X \in \linfty; \quad \Pup(A|B) = \limsupn \vv{P}(A|B)(n), \quad A \subseteq \Omega.
\]
\end{proposition}
Also, we obtain the corresponding lower quantities $\Rup(X|B) = \liminf_{n \rightarrow \infty}  \avgseq|B(n)$ and $\Plow(A|B) = \liminf_{n \rightarrow \infty}  \vv{P}(A|B)(n)$. 
Note that these definitions also have reasonable frequentist semantics even when $B$ 
occurs only finitely often; then the sequence $\vv{P}(\cdot|B)$ is eventually constant and 
we have $\vv{P}(A|B) = |\mbox{$A$ and $B$ occured}|/|\mbox{$B$ occured}|$. 
For instance, if $A$ and $B$ occur just once, but simultaneously, then $\Pup(A|B)=\Plow(A|B)=1$.
This is an advantage over Kolmogorov's approach, where conditioning on events of measure zero 
is not meaningfully defined.

We now further analyze the conditional upper probability and the conditional upper prevision. 
As a warm-up, we consider the case of precise probabilities. 
If for some event $A \subseteq \Omega$, we have $\Pup(A|B)=\Plow(A|B)$, 
we write $\Def{\tilde{P}(A|B) \coloneqq \lim_{n \rightarrow \infty} \vv{P}(A|B)(n)}$.

\begin{proposition}
\label{prop:condprecise}
    Assume $P(B),P(A\cap B)$ exist for some $A,B \subseteq \Omega$ and $P(B)>0$. Then it holds that $\tilde{P}(A|B)=P(A|B)$, where $P(\cdot|B)$ is the conditional probability in the sense of Equation~\ref{eq:defcondmeasure}.
\end{proposition}
\begin{proof}
\begin{align*}
    \tilde{P}(A|B) &= \lim_{n \rightarrow \infty} \vv{P}(A|B)(n)\\
    &= \lim_{n \rightarrow \infty} \frac{\frac{1}{n} \sum_{i=1}^n 
    (\chi_{A} \cdot \chi_B)\left(\vv{\Omega}(i)\right)}{\frac{1}{n} \sum_{i=1}^n \chi_B\left(\vv{\Omega}(i)\right)}\\
    &\overset{\ref{proof:prop34step1}}{=} \frac{\lim_{n \rightarrow \infty} \frac{1}{n} 
        \sum_{i=1}^n (\chi_{A} \cdot \chi_B)\left(\vv{\Omega}(i)\right)}{\lim_{n \rightarrow \infty}
          \frac{1}{n} \sum_{i=1}^n \chi_B\left(\vv{\Omega}(i)\right)}\\
    &\overset{}{=}\frac{P(A\cap B)}{P(B)}\\
    &\overset{\ref{proof:prop34step2}}{=} P(A|B).
\end{align*}
\begin{enumerate}[nolistsep,label=(\arabic*), ref=(\arabic*)]
   \item \label{proof:prop34step1} The limits exist by assumption and the denominator is $>0$.
   \item \label{proof:prop34step2} In the sense of Equation~\ref{eq:defcondmeasure}.
\end{enumerate}
\end{proof}
Thus, when the relative frequencies of $B$ and $A \cap B$ converge, we reproduce the classical definition of conditional probability. Now what happens under non-convergence?

\subsection{The Generalized Bayes Rule}
We now relax the assumptions of Proposition~\ref{prop:condprecise} and only demand that $\Plow(B)>0$.\footnote{This condition is indispensable in order to make the connection to the generalized Bayes rule.} Then we observe that the conditional upper prevision coincides with the \textit{generalized Bayes rule}, which is an important updating principle in imprecise probability (see \eg \citep{introtoiplowerprev}).
The unconditional set of desirable gambles is:
\[
    \Def{\mathcal{D}_{\vv{\Omega}} \coloneqq \left\{X \in \linfty \colon \limsupn \avgseq \leq 0\right\} = 
    \left\{X \in \linfty\colon \Rup(X) \leq 0\right\}.}
\]
\begin{definition}
\label{def:gbrdef}
For $\Plow(B)>0$, we define the \Def{\emph{conditional set of desirable gambles}} as:
\[
\Def{\mathcal{D}_{\vv{\Omega}|B} \coloneqq \left\{X \in \linfty \colon X 
      \chi_B \in \mathcal{D}_{\vv{\Omega}}\right\} = \left\{X \in \linfty\colon \limsupn 
      \vv{\Sigma (X\chi_B)} \leq 0\right\},}
\]
and a corresponding upper prevision, which we call the \Def{\emph{generalized Bayes rule}}, as:
\begin{align}
\label{eq:ourgbrdef}
\Def{\gbr(X|B)} &\Def{\coloneqq \inf\left\{\alpha \in \mR \colon X - 
      \alpha \in \mathcal{D}_{\vv{\Omega}|B}\right\}}\\
&= \inf\left\{\alpha \in \mR \colon \chi_B (X - \alpha) \in \mathcal{D}_{\vv{\Omega}}\right\}\nonumber\\
&= \inf\left\{\alpha \in \mR \colon \Rup\left(\chi_B (X - \alpha)\right) \leq 0\right\}.\nonumber
\end{align}
\end{definition}

\begin{remark}
\label{remark:gbrdef}
\normalfont
In fact, \citet[Section 6.4]{walley1991statistical} defines the generalized Bayes rule as the solution of $\Rup(\chi_B(X-\alpha))=0$ for $\alpha$. It can be checked that this solution coincides with Definition~\ref{def:gbrdef},\footnote{The conditional set of desirable gambles is considered for instance in \citep{augustin2014introduction} and \citep{wheeler2021gentle}, but there the link to the generalized Bayes rule is not made technically clear.} see Appendix~\ref{app:gbrdefcoincidence}.
\end{remark}

\begin{proposition}
\label{prop:gbrcoincidence}
Let $\Plow(B)>0$. It holds that $\Rup(X|B)=\gbr(X|B)$.
\end{proposition}
\begin{proof}
It is not hard to check that $\Rup(\cdot|B)$ is a coherent upper prevision on $\linfty$, hence we can represent it as \citep[Theorem 3.8.1]{walley1991statistical}:
\[
\Rup(X|B) = \inf\left\{\alpha \in \mR\colon X-\alpha \in \mathcal{D}_{\Rup(\cdot|B)}\right\}, 
\quad \mbox{\ where\ \ } \Def{\mathcal{D}_{\Rup(\cdot|B)} \coloneqq \left\{X \in \linfty\colon 
  \Rup(X|B) \leq 0\right\}}.
\]
We show that $\Rup(X|B)=\gbr(X|B)$ by showing that $\mathcal{D}_{\Rup(\cdot|B)} = \mathcal{D}_{\vv{\Omega}|B}$.

Let $X \in L^\infty$. On the one hand, we know
\begin{align}
    X \in \mathcal{D}_{\vv{\Omega}|B}  \Longleftrightarrow X\chi_B \in \mathcal{D}_{\vv{\Omega}}
    &\Longleftrightarrow  \Rup\left(X\chi_B\right) \leq 0\nonumber\\
    \label{eq:generalizedbayesrulelimsup}&\Longleftrightarrow \limsupn \sum_{i=1}^n \frac{(X\chi_B)
    \left(\vv{\Omega}(i)\right)}{n} \leq 0.
\end{align}
On the other hand,
\begin{align}
    \label{eq:frequencyconditionallimsup}
    X \in \mathcal{D}_{\Rup(\cdot|B)}  \Longleftrightarrow \Rup(X|B) \leq 0 \Longleftrightarrow 
    \limsupn \frac{\frac{1}{n} \sum_{i=1}^n (X\chi_B)\left(\vv{\Omega}(i)\right)}{\frac{1}{n} 
    \sum_{i=1}^n \chi_B\left(\vv{\Omega}(i)\right)} \leq 0.
\end{align}
It remains to show that the two limit statements (Equation~\ref{eq:generalizedbayesrulelimsup} 
and Equation~\ref{eq:frequencyconditionallimsup}) are equivalent. Due to the limit operation we 
can neglect the terms $n=1,\ldots,n_B-1$. Furthermore, we know that $\vv{b}(n) \in (0,1]$, $n \geq n_B$, 
and also $0 < \lim \inf_{n \rightarrow \infty} \vv{b}(n)$. 
Thus, defining $\Def{\vv{a}(n) \coloneqq \frac{1}{n}\sum_{i=1}^n (X\chi_B)\left(\vv{\Omega}(i)\right)}$ 
and $\Def{\vv{b}(n) \coloneqq \frac{1}{n} \sum_{i=1}^n \chi_B\left(\vv{\Omega}(i)\right)}$, we can 
leverage Lemma~\ref{lemma:gbrlemma},\footnote{To rigorously apply the Lemma, we would again introduce 
a wrapper for the sequence $\vv{b}(n)$ to ensure strict positivity, since finitely many terms 
$i=1,\ldots, n_B-1$ might be zero.} included in Appendix~\ref{app:gbrlemma} to show:
\begin{align*}
\limsupn \vv{a}(n) \leq 0 &\Longleftrightarrow \limsupn \frac{\vv{a}(n)}{\vv{b}(n)} \leq 0.
\end{align*}

\end{proof}

\begin{remark} \normalfont
\label{remark:wronggbr}
Note that $X \mapsto \limsupn \vv{\Sigma (X\chi_B)} = \Rup\left(X \chi_B\right)$ is \emph{not} in 
general a coherent upper prevision on $\linfty$, as it can violate \ref{item:UP1};
see Appendix~\ref{app:wronggbr}. In general, 
we have $\gbr(X|B)\neq \Rup\left(X \chi_B\right)$.
\end{remark}

As a consequence, we can apply the classical representation result for the generalized Bayes rule.

\begin{corollary}
If $\Plow(B)>0$, the conditional upper prevision can be obtained by updating each linear prevision in the set of cluster points, that is:
\[
    \Rup(X|B) = \sup\left\{E(X|B) \colon E \in \CP\left(\vv{E}\right)\right\},
\]
where conditioning of the linear previsions is in the sense of Definition \ref{def:condlinearprev}.
\end{corollary}
This follows from \citep[Theorem 6.4.2]{walley1991statistical}. 
Intuitively, it makes no difference whether we consider the cluster points of the 
sequence of conditional probabilities or whether we condition all probabilities in the 
set of cluster points in the classical sense.


\section{Unstable Independence}
Closely related to conditional probability is the concept of \textit{statistical independence}.
Independence plays a central role not only in Kolmogorov's \citep[p.\@\xspace 37]{Durrett:2019tt}, but more generally in most probability theories (\citet{levin1980concept}; \citet[Section IIF, IIIG and VH]{fine2014theories}). Already \citet[Introduction, p.\@\xspace 6]{DeMoivre1738} nicely summarized a pre-theoretical, probabilistic notion of real-world independence:
\begin{quote}
Two events are independent, when they have no connexion one with the other,
and that the happening of one neither forwards nor obstructs the happening of
the other.
\end{quote}
This intuitive conception was then formalized by \cite{kolmogorov1933grundbegriffe} (translated in~\citep{kolmogorov2018foundations})  in the following classical definition.
\begin{definition} Let $(\Omega,\mathcal{F},P)$ be a probability space.\footnote{Here, $\Omega$ is the possibility space, $\mathcal{F}$ is a $\sigma$-algebra and $P$ is a countably additive probability measure.}
We call events $A \in \mathcal{F}$ and $B \in \mathcal{F}$ \Def{\emph{classically independent}}
if $P(A\cap B) = P(A)P(B)$.
\end{definition}
If $P(B)>0$, then we can equivalently express this condition as $P(A|B)=P(A)$ by using the definition of conditional probability.

Kolmogorov's definition is formal and it has been questioned whether it is an adequate expression of what 
we mean by independence in a statistical context \citep{von2006note}. As it is stated in purely 
measure-theoretic terms, it is unclear whether it has reasonable frequentist semantics. 
In our framework, we construct an intuitive definition of independence, where the independence of 
events is based on an independence notion of \textit{processes} 
(\cf \citep[p.\@\xspace35-39]{mises1964mathematical}). Therefore, our definition is 
thoroughly grounded in the frequentist setting. Furthermore, we shall generalize the 
independence concept to the case of possible divergence, where new subtleties come into play.
We will then consider how our definitions relate to the classical case when relative 
frequencies converge. Assume that a sequence $\vv{\Omega}$ is given and we have 
constructed an upper probability $\Pup$ as in Section~\ref{sec:ivanenkoformal}.

\begin{definition}
We call an event $B \in \eventonce$ \Def{\emph{irrelevant}} to another event $A \subseteq \Omega$ if:
\begin{align*}
\Pup(A|B)&=\Pup(A).
\end{align*}
\end{definition}

This definition captures the concept of \textit{epistemic irrelevance} in the imprecise probability literature \citep{miranda2008survey}. Why does this definition possess reasonable frequentist semantics? Consider what $\Pup(A|B)$ means (see Section~\ref{sec:unstablecondprob}): we are considering a subsequence, induced by the indicator gamble $\chi_B$, that is, we condition (in an intuitive sense) on the occurence of $B$; and on this subsequence, we then consider an unconditional upper probability. If this then coincides with the orginal upper probability, our decision maker values $A$ just the same whether $B$ occurs or not. Thus $B$ is irrelevant for putting a value on $A$.
In contrast to the classical, precise case, irrelevance is not necessarily symmetric. Hence, we define independence as follows.
\begin{definition}
\label{def:indevents}
Let $A,B \in \eventonce$.
We call $A$ and $B$ \Def{\emph{independent}} if $\Pup(A|B)=\Pup(A)$ and $\Pup(B|A)=\Pup(B)$.
\end{definition}

Thus, we have obtained a grounded concept of independence for events. We note that Definition~\ref{def:indevents} is similar to a condition proposed by \citet{walley1982towards} for independence of joint experiments\footnote{\citet{walley1982towards} considered outcomes of ``joint experiments'' in $\Omega \times \Omega$. They furthermore demanded that lower limits of relative frequencies factorize; to us it is not clear from a strictly frequentist perspective why this condition should be introduced.}; they did not propose an independence concept for gambles.

How can we extend this to an irrelevance and independence concept for gambles? First, we briefly recall how this is done in the classical case.
\begin{definition}
\label{def:independencervsclassical}
Let $(\Omega,\mathcal{F},P)$ be a probability space and fix the Borel $\sigma$-algebra $\mathcal{B}$ on $\mR$. 
Given two gambles $X,Y \colon \Omega \rightarrow \mR$, we say that they are \Def{\emph{classically independent}} 
if:
\[
P(A \cap B) = P(A)P(B) \quad \forall A \in \sigma(X), B \in \sigma(Y),
\]
where the $\sigma$-algebra generated from a gamble $X$, $\sigma(X)$, is defined as the smallest $\sigma$-algebra which $X$ is measurable with respect to:
\[
    \sigma(X) \coloneqq \sigma\left(X^{-1}(\mathcal{B})\right),
\]
and $\sigma(\mathcal{H})$ is the smallest $\sigma$-algebra containing all sets $H \in \mathcal{H}$, $\mathcal{H} \subseteq 2^\Omega$.

\end{definition}
Thus independence of gambles is reduced to independence of events. But note that this definition inherently depends on the choice of the Borel $\sigma$-algebra on $\mR$. 
In our case, this is similar: to define irrelevance and independence on gambles, we need to fix a set system on $\mR$, but we leave the choice open in general.

\begin{definition}
\label{def:irrelevanceofrvs}
Assume a set system $\mathcal{H} \subseteq 2^\mR$ and two gambles $X,Y\colon \Omega \rightarrow \mR$ are given. We call $Y$ irrelevant to $X$ with respect to $\mathcal{H}$ if
    \[
    \Pup(X^{-1}(A)|Y^{-1}(B))=\Pup(X^{-1}(A)) \quad \forall A, B \in \mathcal{H} \text{ if } Y^{-1}(B) \in \eventonce.
    \]
    Similarly, we call them independent when both directions hold.
\end{definition}

Observe that if $\mathcal{H} = \mathcal{B}$ and $\Pup$ was actually a precise $P$ on $\sigma(X)$ and $\sigma(Y)$, this definition would be equivalent to Definition~\ref{def:independencervsclassical} (modulo the subtlety regarding conditioning on measure zero events), due to the following.
\begin{lemma}
Given set systems $\mathcal{A},\mathcal{B} \subseteq \Omega$, in the precise case, the following statements are equivalent.
\begin{enumerate}[nolistsep,label=\textbf{\emph{PI\arabic*.}}, ref=PI\arabic*]
   \item \label{prop:pi1} $P(A|B)=P(A)$ $\forall A \in \mathcal{A}, B \in \mathcal{B}$ and $P(B)>0$.
   \item \label{prop:pi2} $P(A \cap B)=P(A)P(B)$ $\forall A \in \mathcal{A}, B \in \mathcal{B}$.
\end{enumerate}
\end{lemma}
\begin{proof}
Obviously \ref{prop:pi2} implies \ref{prop:pi1} by the definition of conditional probability.
One only has to check that when \ref{prop:pi1} holds, that \ref{prop:pi2} holds even if $P(B)=0$. But if $P(B)=0$, then also $P(A \cap B)=0$ due to monotonicity of $P$ in the sense of a capacity.
\end{proof}

\begin{example} \normalfont
Choose $\mathcal{H} \coloneqq \{(-\infty,a] \colon a \in \mR\}$ in Definition~\ref{def:irrelevanceofrvs}. Such an $\mathcal{H}$ is called a $\Pi$-system, which is a non-empty set system that is closed under finite intersections. This particular $\Pi$-system can in fact be used to define independence in the classical case, which is done in terms of the joint cumulative distribution function. In Appendix~\ref{app:independenceviapi}, we investigate this approach to defining independence and discuss subtle differences to the classical, countably additive case.
\end{example}




\section{Related Work}
\label{sec:relatedwork} 

We examine previous research at the intersection of 
frequentism and imprecise probability. While divergence of relative 
frequencies has been linked to imprecise probability before, 
this has almost exclusively been done in settings which 
are not \textit{strictly} frequentist.
\citet{fine1970apparent} was one of the first authors to critically 
evaluate the hypothesis of statistical stability. 
\citet{fine1970apparent} observed that this widespread hypothesis 
is regarded as a ``striking instance of order in chaos'' in 
the statistics community, and sought to challenge its nature as 
an empirical ``fact''. In contrast to our approach, 
\citet{fine1970apparent} was concerned with finite sequences and 
the question what it means for such a sequence to be random.
While Fine did mention von Mises, \citet{fine1970apparent} opted for 
a randomness definition based on computational complexity. 
Intuitively, one can consider a sequence random if it 
cannot be generated by a short computer program (\ie 
universal Turing machine). Fine then showed that statistical 
stability (``apparent convergence'') occurs \textit{because of},
and not in spite of, high randomness of the sequence. 
In contrast, a sequence for which relative frequencies 
diverge has low computational complexity. We consider these 
findings surprising, and believe that an interesting avenue 
for future research with respect to statistical stability lies 
in the comparison of the computational complexity approach to 
von Mises randomness notion based on selection rules. 
We agree with \citet{fine1970apparent} that apparent 
convergence is not some law of nature, but rather a 
consequence of data handling.

The previously mentioned paper may be seen as a predecessor 
to a long line of work by Terrence  Fine and collaborators, 
\citep{fine1976computational, walley1982towards, kumar1985stationary, grize1987continuous, 
fine1988lower, papamarcou1991stationarity, papamarcou1991unstable,
sadrolhefazi1994finite, fierens2009frequentist}; see also 
\citep{finehandbook} for an introduction. 
A central motivation behind this work was to develop a 
frequentist model for the puzzling case of stationary, 
unstable phenomena with bounded time averages. 
What differentiates this work from ours is that we 
take a \textit{strictly frequentist} approach: 
we explicitly define the upper probability and upper 
prevision from a given sequence. In contrast, 
the above works (with the exceptions of Section 4.3 in \citep{walley1982towards}, \citep{papamarcou1991unstable}
and \citep{fierens2009frequentist}) use an imprecise probability 
to represent a single trial in a sequence of unlinked 
repetitions of an experiment, and then induce an imprecise 
probability via an infinite product space. This is in the 
spirit of, and can be understood as a generalization of, 
the standard frequentist approach, where one would 
assume that $X_1,X_2,\ldots$ form an i.i.d. sequence of 
random variables; here, there is both an ``individual $P$,'' 
as well as an induced ``aggregate $P$'' on the infinite 
product space, which can be used to measure an event such 
as convergence or divergence of relative frequencies.

When a single trial is assumed to be governed by an imprecise 
probability, how can this be interpreted? And what is the 
interpretation of the mysterious ``aggregate imprecise probability''?
This model falls prey to similar criticisms as we outlined 
in the Introduction (Section~\ref{sec:introductionfreqip}) 
concerning the theoretical law of large numbers. In fact, 
\citet{walley1982towards} subscribed to a frequency-propensity 
interpretation (specifically, they were inspired by 
\citet{giere1973objective}), where the imprecise probability 
of a single trial represents its propensity, that is, 
its tendency or disposition to produce a certain outcome. 
Consequently, one obtains a propensity for compound trials 
in terms of an imprecise probability and thus one can 
ascribe a lower and upper probability to events such as 
divergence of relative frequencies. To us, the meaning of 
such a propensity is unclear.
While we are not against a propensity interpretation as such, 
our motivation was to work with a parsimonious 
set of assumptions. To this end, we took the sequence as the 
primitive entity, without relying on an underlying 
``individual'' (imprecise) probability.

Closely related to our work is \citep{papamarcou1991unstable}, 
who were also inspired by von Mises. The authors proved that,
for any set of probability measures $\mathcal{P}$ on 
$(\Omega,2^\Omega)$, $|\Omega|<\infty$, and any countable set 
of place selection rules $\mathcal{S}$,\footnote{See 
\citep{papamarcou1991unstable} for the definition. 
Intuitively, a place selection rule is causal, 
\ie depends only on past values.} the existence of a 
sequence $\vv{\Omega} \colon \mathbb{N} \rightarrow 
\Omega$ with the following property can be guaranteed 
\citep[Theorem 2.2]{papamarcou1991unstable}:
\[
\forall \vv{S}_j \in \mathcal{S} : \forall A \subseteq \Omega : \limsupn \frac{\sum_{i=1}^n  
\chi_A\left(\vvomi\right) \cdot 
\vv{S}_{\!\!j}(i)}{\sum_{i=1}^n \vv{S}_{\!\!j}(i)}  
= \sup\{P(A) \colon P \in \mathcal{P}\}.
\]
That is, the sequence has the specified upper probability 
(take $\vv{S}_{\!\!j}(i)=1$ $\forall i \in \mathbb{N}$)
and this property is stable under subselection.
Note that this claim is in one sense weaker than our 
Proposition~\ref{theorem:converse}, where we construct 
a sequence for which the set of cluster points is 
exactly a prespecified one - coherent upper previsions are more expressive than coherent upper probabilities; on the other hand it is stronger, since the property holds also when applying selection rules.

Within the setup of \citep{walley1982towards}, 
\citet[Theorem 1, Theorem 2]{cozmanconvex} proposed an estimator 
for the underlying imprecise probability of the sequence. 
Specifically, they computed relative frequencies along a 
set of selection rules (however without referring to von Mises)
and then took their minimum to obtain a lower probability;
in a specific technical sense, this estimation succeeds.
What motivated the authors to do this is an assumption on 
the data-generating process: at each trial, ``nature'' 
may select a different distribution from a set of probability 
measures; the trials are then independent but not identically 
distributed.
This viewpoint also motivated \citet{fierens2009frequentist}, 
who restricted themselves to finite sequences. They offered the metaphor 
of an \textit{analytical microscope}. With more and more complex selection 
rules (``powerful lenses''), along which relative frequencies are computed, more and more structure of the set of probabilities comes to light. The authors also proposed a way to simulate data from a set of probability measures.

\citet{cattaneo2017empirical} investigated an empirical, frequentist interpretation of imprecise probability in a similar setting, where $X_1,X_2,..$ is a sequence of precise Bernoulli random variables, but $p_i \coloneqq P(X_i = 1)$ is chosen by nature and may differ from trial to trial, hence $p_i \in [\underline{p_i},\overline{p_i}]$. The author drew the sobering conclusion that ``imprecise probabilities do not have a generally valid,
clear empirical meaning, in the sense discussed in this paper''.

Works which proposed extensions (modifications) of the law of large numbers to imprecise probabilities include \citep{marinacci1999limit,maccheroni2005strong,de2008weak,chen2013strong,peng2019nonlinear}. 

Separate from the imprecise probability literature, \citet{gorban2017statistical} studied the phenomenon of statistical stability and its violations in depth, including theory and experimental studies.
Similarly, the work of \citet{ivanenkobook}, a major motivation for our work, does not appear to be known in the imprecise probability literature.

\section{Conclusion}
In this work, we have extended strict frequentism to the case of possibly divergent relative frequencies and sample averages, tying together threads from \citep{mises1919grundlagen}, \citep{ivanenkobook} and \citep{walley1991statistical}. In particular, we have recovered the generalized Bayes rule from a strictly frequentist perspective. Furthermore, we have established strictly frequentist \textit{semantics} for imprecise probability, by demonstrating that (under the mild assumption that $|\Omega|<\infty$) we can explicitly construct a sequence for which the relative frequencies have a prespecified set of cluster points, corresponding to the coherent upper prevision.

The hypothesis of perfect statistical stability is typically taken for 
granted by practitioners of statistics, without recognizing that 
it is just that --- a \textit{hypothesis}; 
see Appendix \ref{sec:pathologies-or-norm} for an elaboration of this point.
Importantly, when one blindly \textit{assumes} convergence of relative 
frequencies, one will not notice when it is violated --- in the practical case, when only a finite sequence is given, such a violation amounts to instability of relative frequencies even for long observation intervals \citep{gorban2017statistical}. 
In this work, we have rejected the assumption of stability; furthermore, in contrast to other related work, we have aimed to weaken the set of assumptions by taking the concept of a sequence as the primitive. 
However, this gives rise to the critique that no finite part of a sequence has any bearing on what the limit is, as has been pointed out by other authors whose studies attempted a frequentist understanding of imprecise probability (\eg \citep{cattaneo2017empirical}). So what is the empirical content of our theory, what are its practical implications? 

The reader may wonder why we have introduced von Mises frequentist account but not further used selection rules afterwards. In von Mises' framework, the set of selection rules expresses randomness assumptions about the sequence, similar to what the i.i.d.\@\xspace assumption achieves in the standard picture. In our view, randomness assumptions are \textit{the} key to empower generalization in the finite data setting. Hence, to supplement our theory with empirical content, the introduction of selection rules is needed. However, multiple directions can be pursued here. For instance, \citet{papamarcou1991unstable} have defined the concept of an \textit{unstable collective}, where divergence remains unchanged when applying selection rules. By contrast, we could introduce a set of selection rules and assume that relative frequencies converge \textit{within each} selection rule, but to potentially different limits across selection rules.\footnote{As demonstrated by Examples 2 and 3 in \citep{cozmanconvex}, converging relative frequencies within selection rules can lead to both overall convergence or divergence (on the whole sequence).} 
Hence, we view this paper as only the first step of a larger research agenda. The next step is to incorporate randomness assumptions into the picture and explore the connections between various possible approaches, specifically how different ways of relaxing \ref{vm1} and \ref{vm2} are related.



Finally, we remark that an interesting avenue for future research may investigate the use of \textit{nets}, which generalize the concept of a sequence. Indeed, \textit{fraction-of-time probability} \citep{gardnerbook,leskow2006foundations,napolitano2022fraction,gardner2022transitioning} is a theory of probability with remarkable parallels to von Mises' \citeyearpar{mises1919grundlagen}. Instead of sequences, this theory is based on continuous time, hence a net $\vv{\Omega} \colon \mR \rightarrow \Omega$.\footnote{We remark that the notion of a sampling net in \citep{ivanenkobook} is a different one, that is, fraction-of-time probability does not fit this concept, although it is based on the usage of a net.} Sample averages are then given by integration instead of summation. In essence, this amounts to using a different \textit{relative measure} than the counting measure, which is implicit in the work of \citet{mises1919grundlagen}. However, fraction-of-time probability was so far developed only for the convergent case; we expect that a similar construction as in Section~\ref{sec:ivanenkoformal} could be used to extend it to the case of divergence.

\subsubsection*{Acknowledgments}
This work was funded by the Deutsche Forschungsgemeinschaft (DFG, 
German Research Foundation) under Germany’s Excellence Strategy –- 
EXC number 2064/1 –- Project number 390727645.   
The authors thank the International Max Planck Research 
School for Intelligent Systems (IMPRS-IS) for 
supporting Christian Fröhlich and Rabanus Derr.  
Robert Williamson thanks Jingni Yang for a series 
of discussions over several years about Ivanenko's 
work  which provided much inspiration for the present work.

\bibliography{tmlr}

\begin{thebibliography}{128}
\providecommand{\natexlab}[1]{#1}
\providecommand{\url}[1]{\texttt{#1}}
\expandafter\ifx\csname urlstyle\endcsname\relax
  \providecommand{\doi}[1]{doi: #1}\else
  \providecommand{\doi}{doi: \begingroup \urlstyle{rm}\Url}\fi

\bibitem[Albeverio et~al.(2005)Albeverio, Pratsiovytyi, and
  Torbin]{albeverio2005topological}
Sergio Albeverio, Mykola Pratsiovytyi, and Grygoriy Torbin.
\newblock Topological and fractal properties of real numbers which are not
  normal.
\newblock \emph{Bulletin des Sciences Math{\'e}matiques}, 129\penalty0
  (8):\penalty0 615--630, 2005.

\bibitem[Albeverio et~al.(2017)Albeverio, Garko, Ibragim, and
  Torbin]{albeverio2017non}
Sergio Albeverio, Iryna Garko, Muslem Ibragim, and Grygoriy Torbin.
\newblock Non-normal numbers: Full {H}ausdorff dimensionality vs zero
  dimensionality.
\newblock \emph{Bulletin des Sciences Math{\'e}matiques}, 141\penalty0
  (2):\penalty0 1--19, 2017.

\bibitem[Artzner et~al.(1999)Artzner, Delbaen, Eber, and
  Heath]{artzner1999coherent}
Philippe Artzner, Freddy Delbaen, Jean-Marc Eber, and David Heath.
\newblock Coherent measures of risk.
\newblock \emph{Mathematical Finance}, 9\penalty0 (3):\penalty0 203--228, 1999.

\bibitem[A{\v{s}}i{\'c} \& Adamovi{\'c}(1970)A{\v{s}}i{\'c} and
  Adamovi{\'c}]{avsic1970limit}
M.D. A{\v{s}}i{\'c} and D.D. Adamovi{\'c}.
\newblock Limit points of sequences in metric spaces.
\newblock \emph{The American {M}athematical {M}onthly}, 77\penalty0
  (6):\penalty0 613--616, 1970.

\bibitem[Augustin et~al.(2014)Augustin, Coolen, De~Cooman, and
  Troffaes]{augustin2014introduction}
Thomas Augustin, Frank~P.A. Coolen, Gert De~Cooman, and Matthias~C.M. Troffaes.
\newblock \emph{Introduction to imprecise probabilities}.
\newblock John Wiley \& Sons, 2014.

\bibitem[Aveni \& Leonetti(2022)Aveni and Leonetti]{aveni2022most}
Andrea Aveni and Paolo Leonetti.
\newblock Most numbers are not normal.
\newblock In \emph{Mathematical Proceedings of the Cambridge Philosophical
  Society}, pp.\  1--11. Cambridge University Press, 2022.

\bibitem[Bagnold(1983)]{bagnold1983nature}
Ralph~Alger Bagnold.
\newblock The nature and correlation of random distributions.
\newblock \emph{Proceedings of the Royal Society of London. A. Mathematical and
  Physical Sciences}, 388\penalty0 (1795):\penalty0 273--291, 1983.

\bibitem[Bauschke et~al.(2015)Bauschke, Dao, and Moursi]{bauschke2015}
Heinz~H. Bauschke, Minh~N. Dao, and Walaa~M. Moursi.
\newblock On {F}ej{\'{e}}r monotone sequences and nonexpansive mappings.
\newblock \emph{Linear and Nonlinear Analysis}, 1\penalty0 (2):\penalty0
  287--295, 2015.

\bibitem[Bishop \& Peres(2017)Bishop and Peres]{bishop2017fractals}
Christopher~J. Bishop and Yuval Peres.
\newblock \emph{Fractals in probability and analysis}.
\newblock Cambridge University Press, 2017.

\bibitem[Borel(1963)]{borel1950}
{\'{E}}mile Borel.
\newblock \emph{Probability and Certainty}.
\newblock Walker and Company, 1963.
\newblock Translation of \emph{Probabilit\'{e} et Certitude}, Presses
  Universitaires de France, 1950.

\bibitem[Bronshteyn \& Ivanov(1975)Bronshteyn and
  Ivanov]{bronshteyn1975approximation}
Efim~M. Bronshteyn and L.D. Ivanov.
\newblock The approximation of convex sets by polyhedra.
\newblock \emph{Siberian Mathematical Journal}, 16\penalty0 (5):\penalty0
  852--853, 1975.

\bibitem[Calude \& Zamfirescu(1999)Calude and Zamfirescu]{calude1999most}
Cristian Calude and Tudor Zamfirescu.
\newblock Most numbers obey no probability laws.
\newblock \emph{Publicationes Mathematicae Debrecen}, 54\penalty0
  (Supplement):\penalty0 619--623, 1999.

\bibitem[Calude(2002)]{calude2002information}
Cristian~S. Calude.
\newblock \emph{Information and randomness: an algorithmic perspective}.
\newblock Springer Science \& Business Media, 2002.

\bibitem[Calude et~al.(2003)Calude, Marcus, and Staiger]{calude2003topological}
Cristian~S. Calude, Solomon Marcus, and Ludwig Staiger.
\newblock A topological characterization of random sequences.
\newblock \emph{Information Processing Letters}, 88\penalty0 (5):\penalty0
  245--250, 2003.

\bibitem[Canguilhem(1978)]{canguilhem1978normal}
Georges Canguilhem.
\newblock \emph{On the Normal and the Pathological}.
\newblock D. Riedel, 1978.

\bibitem[Cattaneo(2017)]{cattaneo2017empirical}
Marco E. G.~V. Cattaneo.
\newblock Empirical interpretation of imprecise probabilities.
\newblock In Alessandro Antonucci, Giorgio Corani, Inés Couso, and Sébastien
  Destercke (eds.), \emph{Proceedings of the Tenth International Symposium on
  Imprecise Probability: Theories and Applications}, volume~62 of
  \emph{Proceedings of Machine Learning Research}, pp.\  61--72. PMLR, 2017.

\bibitem[Chen(1991)]{chen1991nonequilibrium}
Ping Chen.
\newblock Nonequilibrium and nonlinearity --- a bridge between the two
  cultures.
\newblock In George~P. Scott (ed.), \emph{Time, rhythms, and chaos in the new
  dialogue with nature}, pp.\  67--85. Iowa State University Press, 1991.

\bibitem[Chen et~al.(2013)Chen, Wu, and Li]{chen2013strong}
Zengjing Chen, Panyu Wu, and Baoming Li.
\newblock A strong law of large numbers for non-additive probabilities.
\newblock \emph{International Journal of Approximate Reasoning}, 54\penalty0
  (3):\penalty0 365--377, 2013.

\bibitem[Chow \& Teicher(1988)Chow and Teicher]{chowprobability}
Yuan~Shih Chow and Henry Teicher.
\newblock \emph{Probability theory: independence, interchangeability,
  martingales}.
\newblock Springer, 2nd edition, 1988.

\bibitem[Church(1940)]{church1940concept}
Alonzo Church.
\newblock On the concept of a random sequence.
\newblock \emph{Bulletin of the American Mathematical Society}, 46\penalty0
  (2):\penalty0 130--135, 1940.

\bibitem[Cisewski et~al.(2018)Cisewski, Kadane, Schervish, Seidenfeld, and
  Stern]{cisewski2018standards}
Jessi Cisewski, Joseph~B. Kadane, Mark~J. Schervish, Teddy Seidenfeld, and
  Rafael Stern.
\newblock Standards for modest {B}ayesian credences.
\newblock \emph{Philosophy of Science}, 85\penalty0 (1):\penalty0 53--78, 2018.

\bibitem[Cozman \& Chrisman(1997)Cozman and Chrisman]{cozmanconvex}
Fabio Cozman and Lonnie Chrisman.
\newblock Learning convex sets of probability from data.
\newblock Technical report, Carnegie Mellon University, 1997.
\newblock CMU-RI-TR 97-25.

\bibitem[De~Cooman \& Miranda(2008)De~Cooman and Miranda]{de2008weak}
Gert De~Cooman and Enrique Miranda.
\newblock Weak and strong laws of large numbers for coherent lower previsions.
\newblock \emph{Journal of Statistical Planning and Inference}, 138\penalty0
  (8):\penalty0 2409--2432, 2008.

\bibitem[de~Finetti(1974/2017)]{de2017theory}
Bruno de~Finetti.
\newblock \emph{Theory of probability: A critical introductory treatment}.
\newblock John Wiley \& Sons, 1974/2017.

\bibitem[de~Moivre(1738/1967)]{DeMoivre1738}
Abraham de~Moivre.
\newblock \emph{The Doctrine of Chances: A Method of Calculating the
  Probabilities of Events in Play}.
\newblock Frank Cass and Company Limited, 2nd edition, 1738/1967.

\bibitem[Deitmar(2016)]{deitmar2016analysis}
Anton Deitmar.
\newblock \emph{Analysis}.
\newblock Springer, 2016.

\bibitem[Derr \& Williamson(2022)Derr and Williamson]{derr2022fairness}
Rabanus Derr and Robert~C. Williamson.
\newblock Fairness and randomness in machine learning: Statistical independence
  and relativization.
\newblock \emph{arXiv preprint arXiv:2207.13596}, 2022.

\bibitem[{Desrosi{\`e}res}(1998)]{desrosieres1998politics}
Alain {Desrosi{\`e}res}.
\newblock \emph{The politics of large numbers: A history of statistical
  reasoning}.
\newblock Harvard University Press, 1998.

\bibitem[Dunford \& Schwartz(1988)Dunford and Schwartz]{dunford1988linear}
Nelson Dunford and Jacob~T. Schwartz.
\newblock \emph{Linear operators, part 1: general theory}.
\newblock John Wiley \& Sons, 1988.

\bibitem[{D\"{u}rr}(2001)]{durr2001bohmian}
Detlef {D\"{u}rr}.
\newblock Bohmian mechanics.
\newblock In Jean Bricmont, Detlef D{\"u}rr, Maria~Carla Galavotti, Giancarlo
  Ghirardi, Francesco Petruccione, and Nino Zanghi (eds.), \emph{Chance in
  Physics: Foundations and Perspectives}, pp.\  115--131. Springer, 2001.

\bibitem[D{\"u}rr \& Struyve(2021)D{\"u}rr and Struyve]{durr2021typicality}
Detlef D{\"u}rr and Ward Struyve.
\newblock Typicality in the foundations of statistical physics and {B}orn’s
  rule.
\newblock \emph{Do Wave Functions Jump? Perspectives of the Work of {GianCarlo
  Ghirardi}}, pp.\  35--43, 2021.

\bibitem[{D\"{u}rr} \& Teufel(2009){D\"{u}rr} and Teufel]{durr2009}
Detlef {D\"{u}rr} and Stefan Teufel.
\newblock \emph{Bohmian Mechanics: The Physics and Mathematics of Quantum
  Theory}.
\newblock Springer, 2009.

\bibitem[D{\"u}rr et~al.(1992)D{\"u}rr, Goldstein, and Zanghi]{durr1992quantum}
Detlef D{\"u}rr, Sheldon Goldstein, and Nino Zanghi.
\newblock Quantum equilibrium and the origin of absolute uncertainty.
\newblock \emph{Journal of Statistical Physics}, 67:\penalty0 843--907, 1992.

\bibitem[Durrett(2019)]{Durrett:2019tt}
Rick Durrett.
\newblock \emph{Probability: Theory and Examples}.
\newblock Cambridge University Press, 5th edition, 2019.

\bibitem[Eggleston(1949)]{eggleston1949fractional}
Harold~Gordon Eggleston.
\newblock The fractional dimension of a set defined by decimal properties.
\newblock \emph{The Quarterly Journal of Mathematics}, 20:\penalty0 31--36,
  1949.

\bibitem[Fierens et~al.(2009)Fierens, R{\^e}go, and
  Fine]{fierens2009frequentist}
Pablo~Ignacio Fierens, Leonardo~Chaves R{\^e}go, and Terrence~L. Fine.
\newblock A frequentist understanding of sets of measures.
\newblock \emph{Journal of Statistical Planning and Inference}, 139\penalty0
  (6):\penalty0 1879--1892, 2009.

\bibitem[Fine(1970)]{fine1970apparent}
Terrence~L. Fine.
\newblock On the apparent convergence of relative frequency and its
  implications.
\newblock \emph{IEEE Transactions on Information Theory}, 16\penalty0
  (3):\penalty0 251--257, 1970.

\bibitem[Fine(1973)]{fine2014theories}
Terrence~L. Fine.
\newblock \emph{Theories of probability: An examination of foundations}.
\newblock Academic Press, 1973.

\bibitem[Fine(1976)]{fine1976computational}
Terrence~L. Fine.
\newblock A computational complexity viewpoint on the stability of relative
  frequency and on stochastic independence.
\newblock In William~L. Harper and Clifford~Alan Hooker (eds.),
  \emph{Foundations of Probability Theory, Statistical Inference, and
  Statistical Theories of Science, Volume 1}, pp.\  29--40. Springer, 1976.

\bibitem[Fine(1988)]{fine1988lower}
Terrence~L. Fine.
\newblock Lower probability models for uncertainty and nondeterministic
  processes.
\newblock \emph{Journal of Statistical Planning and Inference}, 20\penalty0
  (3):\penalty0 389--411, 1988.

\bibitem[Fine(2016)]{finehandbook}
Terrence~L. Fine.
\newblock Mathematical alternatives to standard probability that provide
  selectable degrees of precision.
\newblock In Alan Hájek and Christopher Hitchcock (eds.), \emph{{The Oxford
  Handbook of Probability and Philosophy}}. Oxford University Press, 2016.

\bibitem[Fr{\"o}hlich \& Williamson(2022)Fr{\"o}hlich and
  Williamson]{frohlich2022risk}
Christian Fr{\"o}hlich and Robert~C. Williamson.
\newblock Risk measures and upper probabilities: Coherence and stratification.
\newblock \emph{arXiv preprint arXiv:2206.03183}, 2022.

\bibitem[Galton(1889)]{galton1889natural}
Francis Galton.
\newblock \emph{Natural inheritance}.
\newblock Macmillan and Company, 1889.

\bibitem[Galvan(2006)]{galvan2006bohmian}
Bruno Galvan.
\newblock Bohmian mechanics and typicality without probability.
\newblock \emph{arXiv preprint quant-ph/0605162}, 2006.

\bibitem[Galvan(2007)]{galvan2007typicality}
Bruno Galvan.
\newblock Typicality vs. probability in trajectory-based formulations of
  quantum mechanics.
\newblock \emph{Foundations of Physics}, 37\penalty0 (11):\penalty0 1540--1562,
  2007.

\bibitem[Galvan(2022)]{galvan2022non}
Bruno Galvan.
\newblock Non-probabilistic typicality, with application to quantum mechanics.
\newblock \emph{arXiv preprint arXiv:2209.14985}, 2022.

\bibitem[Gardner(1986)]{gardnerbook}
William~A. Gardner.
\newblock \emph{Statistical Spectral Analysis: A Nonprobabilistic Theory}.
\newblock Prentice-Hall, Inc., USA, 1986.

\bibitem[Gardner(2022)]{gardner2022transitioning}
William~A. Gardner.
\newblock Transitioning away from stochastic process models, 2022.
\newblock URL
  \url{https://cyclostationarity.com/wp-content/uploads/2022/07/Transitioning-Away-from-Stochastic-Process-Models.pdf}.
\newblock Accessed: 2023-01-27.

\bibitem[Giere(1973)]{giere1973objective}
Ronald~N. Giere.
\newblock Objective single-case probabilities and the foundations of
  statistics.
\newblock In \emph{Studies in Logic and the Foundations of Mathematics},
  volume~74, pp.\  467--483. Elsevier, 1973.

\bibitem[Goldstein(2012)]{goldstein2012typicality}
Sheldon Goldstein.
\newblock Typicality and notions of probability in physics.
\newblock In \emph{Probability in physics}, pp.\  59--71. Springer, 2012.

\bibitem[Gorban(2011)]{gorban2011statistical}
Igor~I. Gorban.
\newblock Statistical instability of physical processes.
\newblock \emph{Radioelectronics and Communications Systems}, 54\penalty0
  (9):\penalty0 499, 2011.

\bibitem[Gorban(2017)]{gorban2017statistical}
Igor~I. Gorban.
\newblock \emph{The Statistical Stability Phenomenon}.
\newblock Springer, 2017.

\bibitem[Gorban(2018)]{gorban2018randomness}
Igor~I. Gorban.
\newblock \emph{Randomness and Hyper-randomness}.
\newblock Springer, 2018.

\bibitem[Grize \& Fine(1987)Grize and Fine]{grize1987continuous}
Yves~L. Grize and Terrence~L. Fine.
\newblock Continuous lower probability-based models for stationary processes
  with bounded and divergent time averages.
\newblock \emph{The Annals of Probability}, 15\penalty0 (2):\penalty0 783--803,
  1987.

\bibitem[Gu \& Lutz(2011)Gu and Lutz]{gu2011effective}
Xiaoyang Gu and Jack~H. Lutz.
\newblock Effective dimensions and relative frequencies.
\newblock \emph{Theoretical computer science}, 412\penalty0 (48):\penalty0
  6696--6711, 2011.

\bibitem[Gudder(1979)]{gudder1979stochastic}
Stanley~P. Gudder.
\newblock \emph{Stochastic methods in quantum mechanics}.
\newblock North Holland, 1979.

\bibitem[Guttmann(1999)]{guttmann1999concept}
Yair~M. Guttmann.
\newblock \emph{The concept of probability in statistical physics}.
\newblock Cambridge University Press, 1999.

\bibitem[H{\'a}jek(2009)]{hajek2009fifteen}
Alan H{\'a}jek.
\newblock Fifteen arguments against hypothetical frequentism.
\newblock \emph{Erkenntnis}, 70\penalty0 (2):\penalty0 211--235, 2009.

\bibitem[Holmes(1975)]{holmes1975geometric}
Richard Holmes.
\newblock \emph{Geometric functional analysis and its applications.}
\newblock Springer, 1975.

\bibitem[Ivanenko(2010)]{ivanenkobook}
Victor~I. Ivanenko.
\newblock \emph{Decision Systems and Nonstochastic Randomness}.
\newblock Springer, 2010.

\bibitem[Ivanenko \& Labkovskii(1986{\natexlab{a}})Ivanenko and
  Labkovskii]{ivanenko1986functional}
Victor~I. Ivanenko and Valery~A. Labkovskii.
\newblock On the functional dependence between the available information and
  the chosen optimality principle.
\newblock In Vadim~I. Arkin, A.~Shiraev, and R.~Wets (eds.), \emph{Stochastic
  Optimization}, pp.\  388--392. Springer, 1986{\natexlab{a}}.

\bibitem[Ivanenko \& Labkovskii(1986{\natexlab{b}})Ivanenko and
  Labkovskii]{ivanenkoclassofcriterion}
Victor~I. Ivanenko and Valery~A. Labkovskii.
\newblock A class of criterion-choosing rules.
\newblock \emph{Doklady Akademii Nauk SSSR}, 31\penalty0 (3):\penalty0
  204--205, 1986{\natexlab{b}}.

\bibitem[Ivanenko \& Labkovskii(1990)Ivanenko and
  Labkovskii]{ivanenkomodelofnonstochasticrussian}
Victor~I. Ivanenko and Valery~A. Labkovskii.
\newblock A model of nonstochastic randomness.
\newblock \emph{Doklady Akademii Nauk SSSR}, 310\penalty0 (5):\penalty0
  1059–1062, 1990.

\bibitem[Ivanenko \& Labkovskii(1993)Ivanenko and Labkovskii]{ivanenko1993}
Victor~I. Ivanenko and Valery~A. Labkovskii.
\newblock On the studying of mass phenonmena.
\newblock In \emph{Proceedings of the 12th IFAC Triennial World Congress,
  Sydney, Australia}, pp.\  693--696, 1993.

\bibitem[Ivanenko \& Labkovskii(2015)Ivanenko and
  Labkovskii]{ivanenkoonregularities}
Victor~I. Ivanenko and Valery~A. Labkovskii.
\newblock On regularities of mass phenomena.
\newblock \emph{Sankhya A}, 77:\penalty0 237--248, 2015.

\bibitem[Ivanenko \& Munier(2000)Ivanenko and Munier]{Ivanenko2000DecisionMI}
Victor~I. Ivanenko and Bertrand Munier.
\newblock Decision making in `random in a broad sense' environments.
\newblock \emph{Theory and Decision}, 49:\penalty0 127--150, 2000.

\bibitem[Ivanenko \& Pasichnichenko(2017)Ivanenko and
  Pasichnichenko]{ivanenko2017expected}
Victor~I. Ivanenko and Illia Pasichnichenko.
\newblock Expected utility for nonstochastic risk.
\newblock \emph{Mathematical Social Sciences}, 86:\penalty0 18--22, 2017.

\bibitem[Katsikopoulos et~al.(2021)Katsikopoulos, {\c{S}im\c{s}ek}, Buckmann,
  and Gigerenzer]{katsikopoulos2021classification}
Konstantinos~V. Katsikopoulos, {\"{O}zg\"{u}r}~{\c{S}im\c{s}ek}, Marcus
  Buckmann, and Gerd Gigerenzer.
\newblock \emph{Classification in the wild: The science and art of transparent
  decision making}.
\newblock MIT Press, 2021.

\bibitem[Khrennikov(2009)]{khrennikov2009interpretations}
Andrei~Y. Khrennikov.
\newblock \emph{Interpretations of Probability -- Second revised and extended
  edition}.
\newblock Walter de Gruyter, 2009.

\bibitem[Khrennikov(2013)]{khrennikov2013p}
Andrei~Y. Khrennikov.
\newblock \emph{$p$-Adic valued distributions in mathematical physics}.
\newblock Springer, 2013.

\bibitem[Khrennikov(2016)]{khrennikov2016probability}
Andrei~Y. Khrennikov.
\newblock \emph{Probability and randomness: quantum versus classical}.
\newblock Imperial College Press, 2016.

\bibitem[Kolmogorov(1933)]{kolmogorov1933grundbegriffe}
Andre{\u{\i}}~Nikolaevich Kolmogorov.
\newblock \emph{{Grundbegriffe der Wahrscheinlichkeitsrechnung}}.
\newblock Springer, 1933.

\bibitem[Kolmogorov(1956)]{kolmogorov2018foundations}
Andre{\u{\i}}~Nikolaevich Kolmogorov.
\newblock \emph{Foundations of the theory of probability: Second English
  Edition}.
\newblock Chelsea Publishing Company, 1956.

\bibitem[Kolmogorov(1983)]{kolmogorov1983logical}
Andre{\u{\i}}~Nikolaevich Kolmogorov.
\newblock On logical foundations of probability theory.
\newblock In K.~Ito and J.V. Prokhorov (eds.), \emph{Probability theory and
  mathematical statistics: Proceedings of the Fourth USSR -- Japan Symposium,
  held at Tbilisi, USSR, August 23-29, 1982}, pp.\  1--5. Springer, 1983.

\bibitem[Kumar \& Fine(1985)Kumar and Fine]{kumar1985stationary}
Anurag Kumar and Terrence~L. Fine.
\newblock Stationary lower probabilities and unstable averages.
\newblock \emph{Zeitschrift f{\"u}r Wahrscheinlichkeitstheorie und verwandte
  Gebiete}, 69\penalty0 (1):\penalty0 1--17, 1985.

\bibitem[La~Caze(2016)]{lacazefrequentism}
Adam La~Caze.
\newblock {Frequentism}.
\newblock In Alan Hájek and Christopher Hitchcock (eds.), \emph{{The Oxford
  Handbook of Probability and Philosophy}}. Oxford University Press, 2016.

\bibitem[Lawrence(1972)]{lawrence1972}
J.~Dennis Lawrence.
\newblock \emph{A Catalog of Special Plane Curves}.
\newblock Dover, 1972.

\bibitem[Le{\'s}kow \& Napolitano(2006)Le{\'s}kow and
  Napolitano]{leskow2006foundations}
Jacek Le{\'s}kow and Antonio Napolitano.
\newblock Foundations of the functional approach for signal analysis.
\newblock \emph{Signal processing}, 86\penalty0 (12):\penalty0 3796--3825,
  2006.

\bibitem[Levin(1980)]{levin1980concept}
Leonid~A. Levin.
\newblock A concept of independence with applications in various fields of
  mathematics.
\newblock Technical Report MIT/LCS/TR-235, MIT, Laboratory for Computer
  Science, 1980.

\bibitem[Lockwood(1961)]{lockwood1961}
Edward~H. Lockwood.
\newblock \emph{A Book of Curves}.
\newblock Cambridge University Press, 1961.

\bibitem[Maccheroni \& Marinacci(2005)Maccheroni and
  Marinacci]{maccheroni2005strong}
Fabio Maccheroni and Massimo Marinacci.
\newblock A strong law of large numbers for capacities.
\newblock \emph{The Annals of Probability}, 33\penalty0 (3):\penalty0
  1171--1178, 2005.

\bibitem[MacKenzie et~al.(2007)MacKenzie, Muniesa, and
  Siu]{mackenzie2007economists}
Donald~A. MacKenzie, Fabian Muniesa, and Lucia Siu (eds.).
\newblock \emph{Do economists make markets?: On the performativity of
  economics}.
\newblock Princeton University Press, 2007.

\bibitem[Marinacci(1999)]{marinacci1999limit}
Massimo Marinacci.
\newblock Limit laws for non-additive probabilities and their frequentist
  interpretation.
\newblock \emph{Journal of Economic Theory}, 84\penalty0 (2):\penalty0
  145--195, 1999.

\bibitem[M{\'e}ndez(1981)]{mendez1981law}
C.G. M{\'e}ndez.
\newblock On the law of large numbers, infinite games, and category.
\newblock \emph{The American Mathematical Monthly}, 88\penalty0 (1):\penalty0
  40--42, 1981.

\bibitem[Midgley(2013)]{midgley2013science}
Mary Midgley.
\newblock \emph{Science as salvation: A modern myth and its meaning}.
\newblock Routledge, 2013.

\bibitem[Milton(1981)]{milton1981origin}
John~R. Milton.
\newblock The origin and development of the concept of the `laws of nature'.
\newblock \emph{European Journal of Sociology/Archives Europ{\'e}ennes de
  Sociologie}, 22\penalty0 (2):\penalty0 173--195, 1981.

\bibitem[Miranda(2008)]{miranda2008survey}
Enrique Miranda.
\newblock A survey of the theory of coherent lower previsions.
\newblock \emph{International Journal of Approximate Reasoning}, 48\penalty0
  (2):\penalty0 628--658, 2008.

\bibitem[Miranda \& Cooman(2014)Miranda and Cooman]{introtoiplowerprev}
Enrique Miranda and Gert~de Cooman.
\newblock Lower previsions.
\newblock In \emph{Introduction to Imprecise Probabilities}, chapter~2, pp.\
  28--55. John Wiley \& Sons, Ltd, 2014.

\bibitem[Nagy(2007)]{tvs3notes}
Gabriel Nagy.
\newblock Topological vector spaces {III}: Finite dimensional spaces -- {N}otes
  from the {F}unctional {A}nalysis {C}ourse ({F}all 07 - {S}pring 08), 2007.
\newblock URL
  \url{https://www.math.ksu.edu/~nagy/func-an-2007-2008/top-vs-3.pdf}.
\newblock Accessed: 2023-01-13.

\bibitem[Napolitano \& Gardner(2022)Napolitano and
  Gardner]{napolitano2022fraction}
Antonio Napolitano and William~A. Gardner.
\newblock Fraction-of-time probability: Advancing beyond the need for
  stationarity and ergodicity assumptions.
\newblock \emph{IEEE Access}, 10:\penalty0 34591--34612, 2022.

\bibitem[Nicolis \& Prigogine(1977)Nicolis and Prigogine]{Nicolos1977}
G.~Nicolis and Ilya Prigogine.
\newblock \emph{Self-organizaiton in nonequilibrium systems}.
\newblock John Wiley \& Sons, 1977.

\bibitem[Olsen(2004)]{olsen2004extremely}
L.~Olsen.
\newblock Extremely non-normal numbers.
\newblock \emph{Mathematical Proceedings of the Cambridge Philosophical
  Society}, 137\penalty0 (1):\penalty0 43--53, 2004.

\bibitem[Oxtoby(1980)]{oxtoby1980measure}
John~C. Oxtoby.
\newblock \emph{Measure and category: A survey of the analogies between
  topological and measure spaces (2nd edition)}.
\newblock Springer, 1980.

\bibitem[Papamarcou \& Fine(1991{\natexlab{a}})Papamarcou and
  Fine]{papamarcou1991stationarity}
Adrian Papamarcou and Terrence~L. Fine.
\newblock Stationarity and almost sure divergence of time averages in
  interval-valued probability.
\newblock \emph{Journal of Theoretical Probability}, 4\penalty0 (2):\penalty0
  239--260, 1991{\natexlab{a}}.

\bibitem[Papamarcou \& Fine(1991{\natexlab{b}})Papamarcou and
  Fine]{papamarcou1991unstable}
Adrian Papamarcou and Terrence~L. Fine.
\newblock Unstable collectives and envelopes of probability measures.
\newblock \emph{The Annals of Probability}, 19\penalty0 (2):\penalty0 893--906,
  1991{\natexlab{b}}.

\bibitem[Peng(2019)]{peng2019nonlinear}
Shige Peng.
\newblock \emph{Nonlinear expectations and stochastic calculus under
  uncertainty: with robust CLT and G-Brownian motion}.
\newblock Springer Nature, 2019.

\bibitem[Philip \& Watson(1987)Philip and Watson]{philip1987some}
G.M. Philip and D.F. Watson.
\newblock Some speculations on the randomness of nature.
\newblock \emph{Mathematical Geology}, 19\penalty0 (6):\penalty0 571--573,
  1987.

\bibitem[Pichler(2013)]{pichler2013natural}
Alois Pichler.
\newblock The natural {B}anach space for version independent risk measures.
\newblock \emph{Insurance: Mathematics and Economics}, 53\penalty0
  (2):\penalty0 405--415, 2013.

\bibitem[Pitowsky(2012)]{pitowsky2012typicality}
Itamar Pitowsky.
\newblock Typicality and the role of the lebesgue measure in statistical
  mechanics.
\newblock In \emph{Probability in physics}, pp.\  41--58. Springer, 2012.

\bibitem[Prigogine(1978)]{prigogine1978time}
Ilya Prigogine.
\newblock Time, structure, and fluctuations ({N}obel prize lecture).
\newblock \emph{Science}, 201\penalty0 (4358):\penalty0 777--785, 1978.

\bibitem[Prigogine(1980)]{prigogine1980}
Ilya Prigogine.
\newblock \emph{From Being to Becoming: Time and Complexity in the Physical
  Sciences}.
\newblock W.H. Freeman and Company, 1980.

\bibitem[Prigogine \& Stengers(1985)Prigogine and Stengers]{prigogine1984order}
Ilya Prigogine and Isabelle Stengers.
\newblock \emph{Order out of chaos: Man's new dialogue with nature}.
\newblock Bantam Books, 1985.

\bibitem[Rao \& Rao(1983)Rao and Rao]{rao1983theory}
K.P.S.~Bhaskara Rao and M.~Bhaskara Rao.
\newblock \emph{Theory of charges: a study of finitely additive measures}.
\newblock Academic Press, 1983.

\bibitem[Rivas(2019)]{rivas2019role}
{\'A}ngel Rivas.
\newblock On the role of joint probability distributions of incompatible
  observables in {B}ell and {K}ochen--{S}pecker theorems.
\newblock \emph{Annals of Physics}, 411:\penalty0 167939, 2019.

\bibitem[Rose(2016)]{rose2016end}
Todd Rose.
\newblock \emph{The end of average: How to succeed in a world that values
  sameness}.
\newblock Penguin UK, 2016.

\bibitem[Ruby(1986)]{ruby1986origins}
Jane~E. Ruby.
\newblock The origins of scientific``law''.
\newblock \emph{Journal of the History of Ideas}, 47\penalty0 (3):\penalty0
  341--359, 1986.

\bibitem[Sadrolhefazi \& Fine(1994)Sadrolhefazi and
  Fine]{sadrolhefazi1994finite}
Amir Sadrolhefazi and Terrence~L. Fine.
\newblock Finite-dimensional distributions and tail behavior in stationary
  interval-valued probability models.
\newblock \emph{The Annals of Statistics}, 22\penalty0 (4):\penalty0
  1840--1870, 1994.

\bibitem[Schechter(1997)]{schechter1997handbook}
Eric Schechter.
\newblock \emph{Handbook of Analysis and its Foundations}.
\newblock Academic Press, 1997.

\bibitem[Schnorr(2007)]{schnorr2007zufalligkeit}
Claus~P. Schnorr.
\newblock \emph{Zuf{\"a}lligkeit und Wahrscheinlichkeit: eine algorithmische
  Begr{\"u}ndung der Wahrscheinlichkeitstheorie}.
\newblock Springer, 2007.

\bibitem[Seidenfeld et~al.(2017)Seidenfeld, Cisewski, Kadane, Schervish, and
  Stern]{whenlargeisalso}
Teddy Seidenfeld, Jessica Cisewski, Jay Kadane, Mark Schervish, and Rafael
  Stern.
\newblock When large also is small conflicts between measure theoretic and
  topological senses of a negligible set.
\newblock In \emph{Presented at Pitt Workshop 3/17}, 2017.

\bibitem[Shen(2009)]{shen2009}
Alexander Shen.
\newblock Algorithmic information theory and foundations of probability.
\newblock \emph{arXiv preprint arXiv:0906.4411v1}, 2009.

\bibitem[Sklar(1995)]{sklar1995physics}
Lawrence Sklar.
\newblock \emph{Physics and chance: Philosophical issues in the foundations of
  statistical mechanics}.
\newblock Cambridge University Press, 1995.

\bibitem[Sklar(2000)]{sklar2000topology}
Lawrence Sklar.
\newblock Topology versus measure in statistical mechanics.
\newblock \emph{The Monist}, 83\penalty0 (2):\penalty0 258--273, 2000.

\bibitem[Soros(2009)]{soros2009}
George Soros.
\newblock General theory of reflexivity.
\newblock \emph{The Financial Times}, 27 October 2009.

\bibitem[Spreij(2022)]{measuretheoreticprobnotes}
Peter Spreij.
\newblock Measure theoretic probability, 2022.
\newblock URL
  \url{https://staff.fnwi.uva.nl/p.j.c.spreij/onderwijs/master/mtp.pdf}.
\newblock Lecture notes.

\bibitem[Stylianou(2020)]{stylianou2020typical}
Anastasios Stylianou.
\newblock A typical number is extremely non-normal.
\newblock \emph{arXiv preprint arXiv:2006.02202}, 2020.

\bibitem[Tao(2009)]{taoweakstar}
Terence Tao.
\newblock 245{B}, notes 11: The strong and weak topologies.
\newblock
  \url{https://terrytao.wordpress.com/2009/02/21/245b-notes-11-the-strong-and-weak-topologies/},
  2009.
\newblock Accessed: 2023-01-13.

\bibitem[van Lambalgen(1987)]{lambalgen1987mises}
Michiel van Lambalgen.
\newblock {Von Mises}' definition of random sequences reconsidered.
\newblock \emph{The Journal of Symbolic Logic}, 52\penalty0 (3):\penalty0
  725--755, 1987.

\bibitem[Von~Collani(2006)]{von2006note}
Elart Von~Collani.
\newblock A note on the concept of independence.
\newblock \emph{Economic Quality Control}, 21\penalty0 (1):\penalty0 155--164,
  2006.

\bibitem[von Mises(1919)]{mises1919grundlagen}
Richard von Mises.
\newblock {Grundlagen der Wahrscheinlichkeitsrechnung}.
\newblock \emph{Mathematische Zeitschrift}, 5\penalty0 (1):\penalty0 52--99,
  1919.

\bibitem[von Mises(1981)]{mises1981probability}
Richard von Mises.
\newblock \emph{Probability, Statistics, and Truth}.
\newblock Dover, 1981.

\bibitem[von Mises \& Geiringer(1964)von Mises and
  Geiringer]{mises1964mathematical}
Richard von Mises and Hilda Geiringer.
\newblock \emph{Mathematical theory of probability and statistics}.
\newblock Academic press, 1964.

\bibitem[Walley(1991)]{walley1991statistical}
Peter Walley.
\newblock \emph{Statistical reasoning with imprecise probabilities}.
\newblock Chapman-Hall, 1991.

\bibitem[Walley \& Fine(1982)Walley and Fine]{walley1982towards}
Peter Walley and Terrence~L. Fine.
\newblock Towards a frequentist theory of upper and lower probability.
\newblock \emph{The Annals of Statistics}, 10\penalty0 (3):\penalty0 741--761,
  1982.

\bibitem[Weinert(1995)]{weinert1995laws}
Friedel Weinert (ed.).
\newblock \emph{Laws of nature: Essays on the philosophical, scientific and
  historical dimensions}, volume~8.
\newblock Walter de Gruyter, 1995.

\bibitem[Wheeler(2021)]{wheeler2021gentle}
Gregory Wheeler.
\newblock A gentle approach to imprecise probability.
\newblock In Thomas Augustin, Fabio Cozman, and Gregory Wheeler (eds.),
  \emph{Reflections on the Foundations of Statistics: Essays in Honor of Teddy
  Seidenfeld}. Springer, 2021.

\bibitem[Williams(1991)]{williams1991probability}
David Williams.
\newblock \emph{Probability with martingales}.
\newblock Cambridge {U}niversity {P}ress, 1991.

\bibitem[Zilsel(1942)]{zilsel1942genesis}
Edgar Zilsel.
\newblock The genesis of the concept of physical law.
\newblock \emph{The Philosophical Review}, 51\penalty0 (3):\penalty0 245--279,
  1942.

\end{thebibliography}
\bibliographystyle{tmlr}

\appendix
\section{Proofs}

\subsection{Weak* compactness of \texorpdfstring{$\PF(\Omega)$}{PF(Omega)}}
\label{app:weakstarcompact}
\citet[Appendix D4]{walley1991statistical} states that the set of linear previsions, $\PF(\Omega)$, 
is compact due to the Alaoglu-Bourbaki theorem, but does not explain how 
this follows. For completeness, we provide an argument.
First, we can observe, like \citet{walley1991statistical}, that the set is weak* closed. 
We use the following well known Lemma, see \eg \citep[Lemma 12.3.4.]{deitmar2016analysis}.
\begin{lemma}
Let $\mathcal{X}$ be a topological space,  $K \subseteq \mathcal{X}$ be compact and $L \subseteq K$ be closed. Then $L$ is compact.
\end{lemma}
Thus we will show that $\PF(\Omega)$ is a subset of some weak* compact set in $\linftyd$,
hence it is in fact weak* compact. From the Alaoglu-Bourbaki theorem (see \eg 
\citep[p.\@\xspace 70]{holmes1975geometric}), we know that the unit ball of the dual norm is weak* compact in $\linftyd$. By definition of the dual norm $\|\cdot\|^*$ of the $\linfty$ norm $\|X\| \coloneqq \sup_{\omega \in \Omega}\{|X(\omega)|\}$, this is the following set:
\begin{align*}
\Def{B} &\Def{\coloneqq \left\{X^* \in \linftyd \colon \|X^*\|^* \leq 1\right\}}\\
&= \left\{X^* \in \linftyd \colon \sup\left\{|X^*(X)| \colon  X \in \linfty \text{ for which } 
    \|X\| \leq 1\right\} \leq 1\right\}\\
&= \left\{X^* \in \linftyd \colon \sup\left\{|X^*(X)| \colon  X \in \linfty \text{ for which }
\sup_{\omega \in \Omega}\left\{|X(\omega)|\right\} \leq 1\right\} \leq 1\right\}.
\end{align*}
Thus, to show that $\PF(\Omega)$ is weak* compact, it suffices to show that $\PF(\Omega) \subseteq B$. That is, given some $E \in \PF(\Omega)$, we show that if $X \in \linfty$ is such that $\sup_{\omega \in \Omega}\{|X(\omega)|\} \leq 1$, then $|E(X)|\leq 1$, since then the supremum over all such $X$ is also $\leq 1$. Thus it suffices to show that $|E(X)| \leq \sup_{\omega \in \Omega} \{|X(\omega)|\}$. We know that $E(|X|) \leq \sup_{\omega \in \Omega}\{|X(\omega)|\}$, since $E$ is a coherent upper prevision, \cf \citet[2.6.1a]{walley1991statistical}. But we also have that $|E(X)| \leq E(|X|)$ from monotonicity and $E(0)=0$, see \eg \citet[Proposition 5]{pichler2013natural}, and thus $|E(X)| \leq E(|X|) \leq \sup_{\omega \in \Omega}\{|X(\omega)|\} \leq 1$, which concludes the proof. 

\subsection{Properties of the Induced Upper Prevision}
\label{app:inducedriskmeasure}
To see that $X-\Rup(X) \in \mathcal{D}_{\vv{\Omega}}$:
\begin{align*}
        X-\Rup(X) \in \mathcal{D}_{\vv{\Omega}} \Longleftrightarrow &\lim_{n\rightarrow \infty} \sup \frac{1}{n} \sum_{i=1}^n \left((X-\Rup(X))\left(\vv{\Omega}(i)\right)\right) \leq 0\\
        \Longleftrightarrow&\limsupn \avgseq - \Rup(X) \leq 0, \quad \text{since } \Rup(X) \text{ is constant}\\
        \Longleftrightarrow&\limsupn \avgseq(n) - \limsupn \avgseq(n) = 0 \leq 0.
    \end{align*}
    Now, show that $\Rup(X) = \limsupn \avgseq(n)$ is in fact the smallest number such that the relation in \eqref{eq:rdeffromdesirable} holds. Suppose there exists $\varepsilon>0$ such that $\Rup(X)-\varepsilon$ (with $\Rup$ defined as before) makes $X-(\Rup(X)-\varepsilon)$ desirable, that is
    \[
        \limsupn \avgseq - \Rup(X) + \varepsilon \leq 0,
    \]
    which is a contradiction due to our choice of $\Rup$. 

\subsection{Proof of Remark~\ref{remark:gbrdef}}
\label{app:gbrdefcoincidence}
In the literature, the generalized Bayes rule is defined as the solution $\alpha^*$ of 
$\Rup\left(\chi_B(X-\alpha)\right)=0$. We show that $\alpha^*=
\inf\left\{\alpha \in \mR \colon \Rup\left(\chi_B (X - \alpha)\right) \leq 0\right\}$. 
Of course, we get for $\alpha \coloneqq \alpha^*$ that equality holds $(=0)$. 
We just have to exclude the possibility that there exists $\tilde{\alpha} < \alpha^*$
so that $\Rup\left(\chi_B (X - \tilde{\alpha})\right) \leq 0$.

Assume such an $\tilde{\alpha}$ exists, so $\Rup\left(\chi_B X - \chi_B \tilde{\alpha}\right) \leq 0$. Write $\tilde{\alpha} + \varepsilon \leq \alpha^*$ for some $\varepsilon>0$. Since $\chi_B X - \chi_B \tilde{\alpha} - \chi_B \varepsilon \geq \chi_B X - \chi_B \alpha^*$, it follows from monotonicity that $\Rup(\chi_B X - \chi_B \tilde{\alpha} - \chi_B \varepsilon) \geq \Rup(\chi_B X - \chi_B \alpha^*)=0$. But from subadditivity, $\Rup(\chi_B X - \chi_B \tilde{\alpha} - \chi_B \varepsilon) \leq \Rup(\chi_B X - \chi_B \tilde{\alpha}) + \Rup(-\chi_B \varepsilon)$ and since $\varepsilon \Rup(-\chi_B)<0$ due to coherence and $\Plow(B)>0$, we have $\Rup(\chi_B X - \chi_B \tilde{\alpha} - \chi_B \varepsilon) < \Rup(\chi_B X - \chi_B \tilde{\alpha})$. Taking this together, we obtain $\Rup(\chi_B X - \chi_B \tilde{\alpha}) > 0$, a contradiction to the assumption.


Thus we have shown that $\alpha^*=\inf\left\{\alpha \in \mR \colon \Rup\left(\chi_B 
(X - \alpha)\right) \leq 0\right\}$. The other expressions in Definition~\ref{def:gbrdef} follow by simple manipulations.

\subsection{Supplement for Proof of Proposition~\ref{prop:gbrcoincidence}}
\label{app:gbrlemma}

\begin{lemma}
\label{lemma:gbrlemma}
Let $\vv{a}\colon \mathbb{N} \rightarrow \mR$ be a sequence and $\vv{b} \colon \mathbb{N} \rightarrow (0,1]$ be a nonnegative sequence such that $\lim \inf_{n \rightarrow \infty} \vv{b}(n) > 0$. Then:
\begin{align*}
\limsupn \vv{a}(n) \leq 0 &\Longleftrightarrow \limsupn \frac{\vv{a}(n)}{\vv{b}(n)} \leq 0.
\end{align*}
\end{lemma}
\begin{proof}
For brevity we simply write $a_n$ and $b_n$ for the sequences $\vv{a}(n)$ and $\vv{b}(n)$.
\begin{align*}
\limsupn a_n \leq 0 &\Longleftrightarrow \limsupn \frac{a_n}{b_n} \leq 0\\
\lim_{n \rightarrow \infty} \left(\sup_{k \geq n} a_k\right) \leq 0 &\Longleftrightarrow \lim_{n \rightarrow \infty} \left(\sup_{k \geq n} \frac{a_k}{b_k} \right) \leq 0.
\end{align*}
where we know that $b_n \in (0,1]$ and furthermore $0 < \lim \inf_{n \rightarrow \infty} b_n \leq \limsupn b_n$. If the sequence $b_n$ would actually converge, then the statement is clearly true, since we can then pull out the limit of $b_n$ (this is allowed).

We begin by showing that $LHS \leq 0 \implies RHS \leq 0$.
Our assumption is that
\begin{equation}
\label{eq:assumptionlhs}
\forall \varepsilon>0 : \exists n_0 \in \mathbb{N} : \forall n \geq n_0: \sup_{k \geq n} a_k < \varepsilon.
\end{equation}
Our aim is to show that 
\[
    \forall \varepsilon'>0 : \exists n_0' \in \mathbb{N} : \forall n \geq n_0': \sup_{k \geq n} \frac{a_k}{b_k} < \varepsilon' .
\]
So let some $\varepsilon'>0$ be given and fixed. We have to exhibit some $n_0'$ such that the above statement holds. Choose $\varepsilon \coloneqq \varepsilon' \cdot \lim_{n\rightarrow \infty}\inf_{k \geq n} b_k \cdot \frac{1}{\kappa}$, for an arbitrary $\kappa > 1$. Then $\varepsilon > 0 $ by our assumption that $\lim_{n\rightarrow \infty}\inf_{k \geq n} b_k > 0$, i.e. that $\Plow(B)>0$. Note that $\varepsilon \leq \varepsilon'$. Then, we know that $\exists n_0(\varepsilon)$ such that $\forall n \geq n_0(\varepsilon)$ $\sup_{k \geq n} a_k < \varepsilon$.

Also, we know that $\forall \kappa' > 1 : \exists n_0'' \in \mathbb{N} : \forall n \geq n_0''$: 
\begin{equation}
\label{eqref:kapparatio}
    \frac{\lim_{n\rightarrow \infty}\inf_{k \geq n} b_k}{\inf_{k \geq n} b_k} \leq \kappa'.
\end{equation}
Since the numerator is the limit of the denominator (which exists) and furthermore $\inf_{k \geq n} b_k$ is monotone increasing in $n$, that is, $\forall n \in \mathbb{N}$: $\lim_{n\rightarrow \infty}\inf_{k \geq n} b_k \geq \inf_{k \geq n} b_k$. Thus, the ratio approaches $1$ from above for large $n$. 

Now choose $\kappa' \coloneqq \kappa$ and $n_0' \coloneqq \max(n_0(\varepsilon), n_0''(\kappa'))$.
That is, we know that then both \eqref{eqref:kapparatio} and \eqref{eq:assumptionlhs} hold. Then we want to show:

\[
\sup_{k \geq n} \frac{a_k}{b_k} = \max\left(\sup_{k \geq n, a_k \geq 0} \frac{a_k}{b_k}, \sup_{k \geq n, a_k < 0} \frac{a_k}{b_k}\right) \stackrel{!}{<} \varepsilon',
\]
which is a legitimate decomposition of the supremum into the ``negative'' and ``nonnegative'' subsequences. But look at the second term ($a_k < 0$) and observe that since $b_k > 0$, clearly $\sup_{k \geq n, a_k < 0} \frac{a_k}{b_k} \leq 0 < \varepsilon'$. 
Thus we only have to consider the first term. Further observe that
\[
    \sup_{k \geq n, a_k \geq 0} \frac{a_k}{b_k} \leq \sup_{k \geq n, a_k \geq 0} a_k \cdot \sup_{k \geq n, a_k \geq 0} \frac{1}{b_k} = \sup_{k \geq n, a_k \geq 0} a_k \cdot \frac{1}{\inf_{k \geq n, a_k \geq 0} b_k},
\]
due to nonnegativity of the $a_k \geq 0$ and a general rule for the supremum/infimum, which applies since $b_k$ is strictly positive. Now by assumption,
\[
 \sup_{k \geq n, a_k \geq 0} a_k \cdot \frac{1}{\inf_{k \geq n, a_k \geq 0} b_k} < \varepsilon' \lim_{n\rightarrow \infty}\inf_{k \geq n} b_k \frac{1}{\kappa} \frac{1}{\inf_{k \geq n, a_k \geq 0} b_k} = \varepsilon' \cdot \underbrace{\frac{\lim_{n\rightarrow \infty}\inf_{k \geq n} b_k}{\inf_{k \geq n, a_k \geq 0} b_k}}_{\leq \kappa} \cdot \frac{1}{\kappa} \leq \varepsilon'.
\]
Noting that $\inf_{k \geq n} b_k \leq \inf_{k \geq n, a_k \geq 0} b_k$ and therefore
\[
    \frac{\lim_{n\rightarrow \infty}\inf_{k \geq n} b_k}{\inf_{k \geq n,  a_k \geq 0} b_k} \leq \frac{\lim_{n\rightarrow \infty}\inf_{k \geq n} b_k}{\inf_{k \geq n} b_k} \leq \kappa'.
\]
Altogether, we have shown that
\[
\sup_{k \geq n} \frac{a_k}{b_k} < \varepsilon',
\]
and therefore $LHS \leq 0 \implies RHS \leq 0$.

It remains to show that $RHS \leq 0 \implies LHS \leq 0$. Our assumption is that 
\begin{equation}
\label{eq:assrhsimplieslhs}
    \forall \varepsilon'>0 : \exists n_0' \in \mathbb{N} : \forall n \geq n_0': \sup_{k \geq n} \frac{a_k}{b_k} < \varepsilon.
\end{equation}
and our aim is to show that then
\[
\forall \varepsilon>0 : \exists n_0 \in \mathbb{N} : \forall n \geq n_0: \sup_{k \geq n} a_k < \varepsilon.
\]
So let $\varepsilon>0$ be fixed. Choose $\varepsilon' \coloneqq \varepsilon$ and set $n_0 \coloneqq n_0'$.
Then we want to show that $\forall n \geq n_0$:
\[
\sup_{k \geq n} a_k = \max\left(\sup_{k \geq n, a_k \geq 0} a_k, \sup_{k \geq n, a_k < 0} a_k\right) \stackrel{!}{<} \varepsilon.
\]
As to the second term, it is obviously negative, in particular $\sup_{k \geq n, a_k < 0} a_k < \varepsilon$. For the first term, where the $a_k$ are nonnegative, observe that then $a_k \leq \frac{a_k}{b_k}$ since $b_k \in (0,1]$, consequently we have $\forall n \geq n_0$:
\[
\sup_{k \geq n, a_k \geq 0} a_k \leq \sup_{k \geq n, a_k \geq 0} \frac{a_k}{b_k} \leq \sup_{k \geq n} \frac{a_k}{b_k} < \varepsilon = \varepsilon'.
\]
by our assumption \eqref{eq:assrhsimplieslhs}. And thus we have shown that $RHS \leq 0  \implies LHS \leq 0$.
\end{proof}

\subsection{Proof of Remark~\ref{remark:wronggbr}}
\label{app:wronggbr}
Take for example $X(\omega)=-1$ for a $B \subseteq \Omega$ where $\Plow(B)<1$. Then $\sup X = -1$, but
\[
\Rup\left(X \chi_B\right) = \limsupn \frac{1}{n} \sum_{i=1}^n \left(X \chi_B\right)\left(\vvomi\right)
= \limsupn - \frac{1}{n} \sum_{i=1}^n \chi_B\left(\vvomi\right) = 
- \liminf_{n \rightarrow \infty} \frac{1}{n} \sum_{i=1}^n \chi_B\left(\vvomi\right) = 
- \Plow(B),
\]
and $\sup X = -1 < - \Plow(B)$, hence \ref{item:UP1} does not hold. Thus 
$X \mapsto \Rup\left(X \chi_B\right)$ is not a coherent upper prevision on $\linfty$ 
in general (it is of course for $B=\Omega$).

\section{Independence via \texorpdfstring{$\Pi$}{Pi}-Systems}
\label{app:independenceviapi}
In this section we discuss a useful special case of defining independence via $\Pi$-systems. This is particularly insightful as it illuminates subtle differences between the definition of irrelevance (independence) in the precise, countable additive and the precise, finitely additive case.

Consider the choice of $\mathcal{H} \coloneqq \{(-\infty,a] \colon a \in \mR\}$ in Definition~\ref{def:irrelevanceofrvs}. It seems like this naturally achieves the goal of expressing independence, but we would like to leave this choice open in general.
Then irrelevance of $Y$ to $X$ means:
\[
    \Pup(\{X \leq a\}|\{Y \leq b\})=\Pup(\{X \leq a\}), \quad \forall a, b \in \mR, \text{ if } \{Y \leq b\} \in \eventonce.
\]
Compare this to the classical, precise setting where independence can also be defined as:
\begin{align}
\label{def:classicalindependence}
 P(\{X \leq a\}|\{Y \leq b\})&=P(\{X \leq a\}), \quad \forall a, b \in \mR, \text{ if } P(\{Y \leq b\})>0.\\
    \Longleftrightarrow\ \  P(\{X \leq a\}\cap\{Y \leq b\})&=P(\{X \leq a\})P(\{Y \leq b\}), \quad \forall a, b \in \mR.\nonumber\\
\Longleftrightarrow\ \ F_{X,Y}(a,b) & = F_X(a)F_Y(b), \quad \forall a, b \in \mR.\nonumber    
\end{align}
That is, in the classical, precise case it suffices to have the joint distribution function factorize. This is formalized in the following.

\begin{proposition}
Let $(\Omega,\mathcal{F},P)$ be a probability space. Assume classical independence as in Equation~\ref{def:classicalindependence} holds. Then they are also independent in the following sense:
\[
    P\left(X^{-1}(A) \cap Y^{-1}(B)\right) = P\left(X^{-1}(A)\right)P\left(Y^{-1}(A)\right), 
    \quad \forall A,B \subseteq \mathcal{B}(\mR).
\]
where $\mathcal{B}(\mR)$ is the Borel $\sigma$-algebra on $\mR$, thereby constituting a 
precise special case of Definition~\ref{def:irrelevanceofrvs}\footnote{Modulo the issue of conditioning on 
measure zero sets.}. Furthermore, they are independent in the sense 
of Definition~\ref{def:independencervsclassical}.
\end{proposition}
Essentially, it suffices to define independence based on the set systems 
$\left\{X^{-1}((-\infty,a]) \colon a \in \mR\right\}$ and $\left\{Y^{-1}((-\infty,a]) \colon 
a \in \mR\right\}$, which are the pre-images of $\mathcal{H} \coloneqq \{(-\infty,a] 
\colon a \in \mR\}$, to get independence on the whole generated $\sigma$-algebras. 
This is based on the famous $\Pi-\lambda$ Theorem.

To investigate this result in our framework, where slight differences will arise 
due to finite additivity, we need to talk about set systems. First, 
we consider the \textit{system of precision}, on which we have precise probabilities.

We define the \Def{\emph{system of precision}} $\Def{\dynkin}$ as the induced subset of $2^\Omega$ on which the relative frequencies converge:
\[
\Def{\dynkin \coloneqq \left\{A \subseteq \Omega\colon \Pup(A)=\Plow(A) =
P(A) = \lim_{n \rightarrow \infty} \frac{1}{n} \sum_{i=1}^n \chi_A(\vvomi)\right\}.}
\]
We show that $\dynkin$ always constitutes a \textit{pre-Dynkin system}, but not in general a \textit{Dynkin system}.

\begin{definition}
A set system $\mathcal{A} \subseteq 2^\Omega$ is called a \Def{\emph{pre-Dynkin system}} if the following conditions hold:
\begin{enumerate}[nolistsep,label=\textbf{\emph{PD\arabic*.}}, ref=PD\arabic*]
\item \label{pd:omega} $\Omega \in \mathcal{A}$.
\item \label{pd:complement} $A \in \mathcal{A} \implies A^\complement \in \mathcal{A}$.
\item \label{pd:closurefinite} If $A,B \in \mathcal{A}$ and $A \cap B =\emptyset$, then $A \cup B \in \mathcal{A}$.
\end{enumerate}
\end{definition}
Thus a pre-Dynkin system is closed under complements and (by induction) under finite union of disjoint sets. If condition \ref{pd:closurefinite} holds also for a countable collection of disjoint sets, \ie if $A_1,.. \in \mathcal{A}$ and $A_i \cap A_j = \emptyset$ for all $i \neq j$ implies $\bigcup_{i=1}^\infty A_i \in \mathcal{A}$, then we speak of a \Def{\emph{Dynkin system}}.

We write $\Def{\dynkingen(\mathcal{H})}$ for the \Def{intersection of all pre-Dynkin systems containing 
a set system $\mathcal{H} \subseteq 2^\Omega$}.

\begin{proposition}
The system of precision $\dynkin$ is a pre-Dynkin system, but not in general a Dynkin system.
\end{proposition}
\begin{proof}
That condition \ref{pd:omega} holds is obvious. Suppose $A \in \dynkin$, \ie 
$\lim_{n \rightarrow \infty} \frac{1}{n} \sum_{i=1}^n \chi_A\left(\vvomi\right)$ exists.
Then also $\lim_{n \rightarrow \infty} \frac{1}{n} \sum_{i=1}^n \left(1-\chi_A\right)\left(\vvomi\right)
= \lim_{n \rightarrow \infty} \frac{1}{n} \sum_{i=1}^n 1 - \chi_A\left(\vvomi\right)$ exists, 
and hence $A^\complement \in \dynkin$, \ie \ref{pd:complement} holds. 
Now suppose $A$ and $B$ are in $\dynkin$ and $A \cap B = \emptyset$. Then 
$\lim_{n \rightarrow \infty} \frac{1}{n} \sum_{i=1}^n \chi_{A \cup B}\left(\vvomi\right) =
\lim_{n \rightarrow \infty} \frac{1}{n} \sum_{i=1}^n \chi_{A}\left(\vvomi\right) + \chi_B\left(\vvomi\right) = \lim_{n \rightarrow \infty} \frac{1}{n} \sum_{i=1}^n \chi_A\left(\vvomi\right) + 
\lim_{n \rightarrow \infty} \frac{1}{n} \sum_{i=1}^n \chi_B\left(\vvomi\right)$ since both 
limits exist by assumption, hence $A \cup B \in \dynkin$.

It remains to give a counterexample to show that closure under countable disjoint union can fail. For this we simply set $\Omega = \mathbb{N}$. Then, we construct a sequence $\vv{\Omega}$ such that there exist a countable set of pairwise disjoint elements in the corresponding system of precision $\dynkin$ such that their union is not an element of $\dynkin$.
Let $\vvomi = \left\langle 1^{[1]}\, 2^{[1]}\, 3^{[2]}\, 4^{[2]}\, 5^{[4]}\, 6^{[4]}\, 
\dots \, i^{[2^{\lceil i/2\rceil}-1]} 
\, \ldots \right\rangle$, for $i\rightarrow\infty$, where $\langle \cdot \rangle$ forms a sequence 
from its inputs. The notation $i^{[j]}$ here means $j$ repetitions of $i$.

For every even natural number $2k, k \in \mathbb{N}$ we have 
$\lim_{n \rightarrow \infty} \frac{1}{n} \sum_{i=1}^n \chi_{\{ 2k\}}\left(\vvomi\right) 
\le \lim_{n \rightarrow \infty} \frac{1}{n} 2^{k} - 1 = 0$. Thus, $\{ 2k\} \in \dynkin$ for every 
$k \in \mathbb{N}$. The disjoint union of all such sets, namely 
$\{ 2k \colon k \in \mathbb{N}\}$, however, is not in $\dynkin$. To see this, we consider the 
sequence $h_i \coloneqq \chi_{\{ 2k \colon k \in \mathbb{N}\}}\left(\vvomi\right) =
\left\langle 0^{[1]}\, 1^{[1]}\, 0^{[2]}\, 1^{[2]}\, 0^{[4]}\, 1^{[4]}\, \ldots 0^{[2^i]}\, 1^{[2^i]}\,
\ldots \right\rangle$. As noticed by \citet[p. 11]{mises1964mathematical}, this sequence has no 
unique frequency limit, i.e. $\lim_{n \rightarrow \infty} \frac{1}{n} \sum_{i=1}^n 
\chi_{\{ 2k \colon k \in \mathbb{N}\}}\left(\vvomi\right)$ does not exist. 
This concludes the proof that the system of precision $\dynkin$ is a 
pre-Dynkin-system, but not generally a Dynkin-system.
\end{proof}
Note that the system of precision need not be closed under intersection \citep{rivas2019role}.

A pre-Dynkin system is closed under finite disjoint union, but as opposed to a Dynkin system not in general under countable disjoint union. A similar relation holds between a \textit{field} (also called algebra) and a $\sigma$-algebra. To avoid confusion, we stick to the name \textit{field}.

\begin{definition}
A set system $\mathcal{A}$ is a \Def{\emph{field}} if the following conditions hold.
\begin{enumerate}[nolistsep,label=\textbf{\emph{FLD\arabic*.}}, ref=FLD\arabic*]
\item \label{field:omega} $\Omega \in \mathcal{A}$.
\item \label{field:complement} $A \in \mathcal{A} \implies A^\complement \in \mathcal{A}$.
\item \label{field:closurefinite} If $A,B \in \mathcal{A}$, then $A \cup B \in \mathcal{A}$.
\end{enumerate}
Then it it is also closed under finite intersections. We write 
$\Def{\operatorname{field}(\mathcal{H})}$ for the \Def{intersection of all fields containing the set system $\mathcal{H} \subseteq 2^\Omega$}.
\end{definition}
As opposed to a $\sigma$-algebra, a field is in general closed only under finite union. Finally, we need the concept of a $\Pi$-system, which is closed under finite intersections.

\begin{definition}
A set system $\mathcal{H}$ is called a \Def{\emph{$\Pi$-system}} if it is non-empty and $A,B \in \mathcal{H} \implies A\cap B \in \mathcal{H}$.
\end{definition}
\begin{example}
\normalfont
A prominent $\Pi$-system is given by $\mathcal{H}=\{(-\infty,a] \colon a \in \mR\}$.
\end{example}

We now begin (with slight modifications) reproducing a series of results which are 
stated in the literature for the interplay of Dynkin systems and $\sigma$-algebras.
In our case, we will restate them for the interplay of pre-Dynkin system with fields.

\begin{proposition}
A set system $\mathcal{A}$ is a field if and only if it is both a pre-Dynkin system and a $\Pi$-system.
\end{proposition}
\begin{proof}
First, we show ``only if''. This is clear, since the field-condition ``closed under finite union'' is equivalent to ``closed under binary intersections'' (using closure under complement and $\Omega \in \mathcal{A}$). Hence it is necessary that $\mathcal{A}$ is a $\Pi$-system. Also, it is necessary that it be a pre-Dynkin system, since a field is closed under arbitrary finite union, so it also must be closed under finite disjoint union; also it is closed under complement.

Next, we show the ``if'', \ie that $\mathcal{A}$ being a pre-Dynkin and a $\Pi$-system imply that $\mathcal{A}$ is a field. We only have to check that it is closed under arbitrary finite union. Consider:
\[
    A \cup B = \Omega \setminus \left(A^\complement \cap B^\complement\right), 
\]
which is in $\mathcal{A}$ due to it being a pre-Dynkin and $\Pi$-system.
\end{proof}

\begin{proposition}
Pre-Dynkin-$\Pi$-Theorem: Let $\mathcal{H}$ be a $\Pi$ system. Then the generated field coincides with the generated pre-Dynkin system.
\end{proposition}
\begin{proof}
Just use the proof in \citep[p.\@\xspace 193]{williams1991probability} and replace the one occurence of ``countable union'' with ``finite union''; consequently, we must replace ``Dynkin'' with ``pre-Dynkin'' everywhere.
\end{proof}

\begin{proposition}
    Assume precise probabilities exist on some $\Pi$-system $\mathcal{H}$. Then the generated pre-Dynkin system $\dynkin(\Pi)$ is contained in the whole precise pre-Dynkin system $\dynkin$. In particular, $\dynkin$ contains a field.
\end{proposition}
\begin{proof}
Obvious.
\end{proof}

\begin{proposition}
    Uniqueness Lemma: Let $\mathcal{A}$ be a $\Pi$-system on $\Omega$, of which the generated field (or equivalently, pre-Dynkin system) is the field $\mathcal{F}$, and assume we have two finitely additive measures $P$, $Q$ on the field $\mathcal{F} \subseteq 2^\Omega$, so that $P(\Omega)=Q(\Omega)=1$ and $P(A)=Q(A)$ $\forall A \in \mathcal{A}$. Then actually $P(A)=Q(A)$ $\forall A \in \mathcal{F}$.
\end{proposition}
\begin{proof}
We show that $\mathcal{H} \coloneqq \{A \in \mathcal{F}\colon P(A)=Q(A)\}$ is a pre-Dynkin system. Clearly, $\Omega\in \mathcal{H}$ and we have closure under complement. We have to show that if $A,B \in \mathcal{H} \subseteq \mathcal{F}$ and $A\cap B =\emptyset$, then $A \cup B \in \mathcal{H}$, \ie $P(A\cup B)=Q(A \cup B)$. But this follows obviously by assumption that $P(A)=Q(A)$, $P(B)=Q(B)$ and $P,Q$ are finitely additive measures on $\mathcal{F}$.

Also, the $\Pi$-system $\mathcal{A}$ is contained in $\mathcal{H}$ by assumption. Then we get from the Pre-Dynkin-$\Pi$-Theorem that actually the field generated by $\mathcal{A}$ is in $\mathcal{H}$, which concludes the proof.
\end{proof}

\begin{example}\normalfont
    Let $\Pi=\{(-\infty,a] \colon a \in \mR\}$. Consider the induced pre-Dynkin system $\dynkingen(\Pi)$. 
    We call this the Borel field by analogy, since the induced Dynkin system is the Borel $\sigma$-algebra.
\end{example}

Perhaps the name ``Borel field'' is unfortunate, as there is no connection to topology anymore, unlike for the Borel $\sigma$-algebra. However, the name serves to emphasize the close relation to the latter.

\begin{proposition}
    Assume we have precise probabilities on $\{\{X \leq x\} \colon x \in \mR\}$ and $\{\{Y \leq y\} \colon y \in \mR\}$ and assume that the irrelevance condition~\ref{def:irrelevanceofrvs} holds on $\mathcal{H} \coloneqq \{(-\infty,a] \colon a \in \mR\}$:
    \[
        P(\{X \leq x\}|\{Y \leq y\}) = P(\{X \leq x\}), \quad \text{ if } \{Y \leq y\} \in \eventonce.
    \]
    Then the independence condition also holds for the Borel field, \ie
    \[
        P\left(X^{-1}(A)|Y^{-1}(B)\right)=P\left(X^{-1}(A)\right), \quad \forall A,B \in \operatorname{field}(\mathcal{H}), Y^{-1}(B) \in \eventonce.
    \]
\end{proposition}
\begin{proof}
Define $\mathcal{A} \coloneqq \{\{X \leq x\} \colon x \in \mR\}$ and $\mathcal{B} \coloneqq \{\{Y \leq y\} \colon y \in \mR\}$. 
With similar reasoning as in \citet[Proposition 3.12]{measuretheoreticprobnotes}, we get that:
\[
        P(A|B) = P(A), \quad \forall A \in \operatorname{field}(\mathcal{A}), B \in \operatorname{field}(\mathcal{B}).
\]
To obtain the statement, it remains to show:
\[
    \operatorname{field}\left(X^{-1}(\mathcal{H})\right) = X^{-1}(\operatorname{field}(\mathcal{H})).
\]
We can follow similar reasoning as in \citep[p.\@\xspace 12, Lemma 1]{chowprobability}, since 
nothing in the argument depends on the $\sigma$-algebra vs. field distinction. 
This concludes the argument.
\end{proof}

This gives a good justification for the \textit{precise case} to define independence via 
the $\Pi$-system $\mathcal{H} \coloneqq \{(-\infty,a] \colon a \in \mR\}$. But 
for imprecise probabilities, we have no such justification and thus should better 
directly use the whole Borel field to define independence on.

\section{Existence of Sequences with Prespecified Relative Frequency
Cluster Points}
\label{app:construction}
In this section we prove  Theorem \ref{th:CP-r-C} and thus demonstate 
the existence of sequences
$x\colon\naturals\rightarrow[k]$ whose corresponding relative frequencies
$r^x\colon\naturals\rightarrow\Delta^k$ have the property that their set of
cluster points $\CP(r^x)=C$, where $C$ is an arbitrary closed rectifiable curve
in $\Delta^k$.  We do so constructively 
by providing an explicit procedure which takes a chosen
$C$ and  constructs  a suitable $x$. we illustrate our method with
some examples.  

\citet[p.\@\xspace 11]{mises1964mathematical}
considered a binary sequence given by 
\[
x=\left\langle 0^{[1]}, 1^{[1]}, 0^{[2]},
1^{[2]}, 0^{[4]}, 1^{[4]}, \ldots 0^{[2^i]}, 1^{[2^i]},
\ldots \right\rangle
\] 
for $i\rightarrow\infty$. (The notation
$\Def{i^{[j]}}$ here means $j$ repetitions of $i$.)   It is a
straightforward calculation to show that the induced relative frequencies
have all elements of $\left[\frac{1}{3},\frac{2}{3}\right]$ as cluster points. When
$x\colon\naturals\rightarrow [k]$ with
$k>2$ a much richer set of behaviors of $\CP(r^x)$ is possible.  The question of how
common sequences with non-convergent relative frequencies are is addressed in Section~\ref{sec:pathologies-or-norm}.

\subsection{Sufficient to Work With Topology Induced by Euclidean Norm on the Simplex}
\label{app:euclideansufficient}
Let $|\Omega|=k<\infty$. A linear prevision $E$ on $\linfty$ is in a one-to-one correspondence with a finitely additive probability $P$, which we can represent as a point in the $(k-1)$-simplex (see Lemma~\ref{lemma:pfsimplex} below). It is convenient to then consider the cluster points of a sequence of such probabilities in the $(k-1)$-simplex, with respect to the topology induced by the Euclidean norm on $\mR^k$.
We show that the notion of such a cluster point coincides with a cluster point in $\linftyd$ with respect to the weak* topology. This is the goal of this section; in particular, we prove Proposition~\ref{prop:simplexcorrespondence}.

Since $|\Omega|=k<\infty$, we can represent any $X \in \linfty$ as:
\[
X(\omega) = c_1 \chi_{\{\omega_1\}} + \cdots + c_k \chi_{\{\omega_k\}},
\]
where $c_i=X(\omega_i)$.

Similarly, any $Z \in \linftyd$ can be represented as:
\begin{align*}
Z(X) &= Z\left(c_1 \chi_{\{\omega_1\}} + \cdots+ c_k \chi_{\{\omega_k\}}\right)\\
&= c_1 Z\left(\chi_{\{\omega_1\}}\right) + \cdots + c_k Z\left(\chi_{\{\omega_k\}}\right),
\end{align*}
since the $Z \in \linftyd$ are linear functionals and the coefficients $c_i$ depend on $X$. 
Intuitively, $Z\left(\chi_{\{\omega_i\}}\right) = P(\omega_i)$ if $Z \in \PF(\Omega)$. 
Define $\Def{d_i \coloneqq Z\left(\chi_{\{\omega_i\}}\right)}$ $\forall i \in 1,\ldots,n$ 
and consequently define
\[
\Def{\|Z\| \coloneqq \sqrt{d_1^2 + \cdots + d_k^2}.}
\]
For a given $Z \in \linftyd$, call $\Def{d_Z\coloneqq(d_1,..,d_k)} \in \mR^k$ the \Def{coordinate representation of $Z$}.

\begin{lemma}
$\|\cdot\|$ is a norm on $\linftyd$.
\end{lemma}
\begin{proof}
\emph{Point-separating property:} if and only if $Z=0$, where $0 \in \linftyd$ is given by $0(X)=0$, 
$\forall X \in \linfty$. But this is easily observed, due to the similar property 
holding for the Euclidean norm: $\|Z\|=0$ if and only if $d_i=0$ $\forall i=1,\ldots,n$.
Since $d_i =Z\left(\chi_{\{\omega_i\}}\right)$, this is the case exactly if $Z=0$.

\emph{Subadditivity:} $\|Y+Z\| \leq \|Y\| + \|Z\|$. For the input $Y+Z$ we get 
$d_i = (Y+Z)\left(\chi_{\{\omega_i\}}\right) = Y\left(\chi_{\{\omega_i\}}\right) +
Z\left(\chi_{\{\omega_i\}}\right)$ and then subadditivity follows from the similar 
property for the Euclidean norm.

\emph{Absolute homogeneity:} $\|\lambda Z\| = |\lambda| \|Z\|$, $\forall \lambda \in \mR$. This follows easily.

With these properties, $\|\cdot\|$ is a valid norm on $\linftyd$.
\end{proof}


We now show that the $(k-1)$-simplex is in a one-to-one correspondence with the set of linear previsions $\PF(\Omega)$ via the coordinate representation.

\begin{lemma}
\label{lemma:pfsimplex}
Let $d \in \Delta^k$. Then $Z(X) \coloneqq c_1 d_1 + ... + c_k d_k \in \PF(\Omega)$. Conversely, let $Z \in \PF(\Omega)$. Then the corresponding $d_Z \in \Delta^k$.
\end{lemma}
\begin{proof}
Let $d \in \Delta^k$ and $Z(X) \coloneqq c_1 d_1 + \cdots + c_k d_k$. Since $\sum_{i=1}^k d_i = 1$
we have immediately that $Z\left(\chi_\Omega\right)=1$, noting that $\chi_\Omega= 
1 \chi_{\{\omega_1\}} + \cdots+ 1 \chi_{\{\omega_k\}}$. Also, if $X \geq 0$, 
\ie $c_i \geq 0$, $\forall i \in [k]$, then $Z(X) \geq 0$ since $d_i \geq 0$ $\forall i$.
Thus, $Z \in \PF(\Omega)$.

Conversely, let $Z$ be a linear prevision, \ie $Z\left(\chi_\Omega\right)=1$ and 
$Z(X) \geq 0$ if $X \geq 0$. From $Z\left(\chi_\Omega\right)=1$ we can deduce that 
$\sum_{i=1}^k d_i = 1$.  If $X \geq 0$, we know that $c_i \geq 0$ $\forall i \in [k]$, 
hence $Z(X) \geq 0$ can only be true if all $d_i \geq 0$. Thus $d_Z \in \Delta^k$.
\end{proof}

We here restate Proposition~\ref{prop:simplexcorrespondence} for convenience.
\begin{proposition}
\label{prop:equivofclusterpoints}
Let $\vv{E}(n) : \mathbb{N} \rightarrow \PF(\Omega)$ be a sequence of linear previsions with 
underlying probabilities $\vv{P}(n) \coloneqq A \mapsto \vv{E}(n)(A)$. 
Then $E \in \CP\left(\vv{E}(n)\right)$ with respect to the weak* topology if and only if 
the sequence $\vv{D} \colon \mathbb{N} \rightarrow \Delta^k$, $\vv{D}(n) \coloneqq 
\left(\vv{P}(n)(\omega_1),\ldots, \vv{P}(n)(\omega_k)\right)$ has as cluster 
point $d_E=\left(E\left(\chi_{\{\omega_1\}}\right),\ldots,E\left(\chi_{\{\omega_k\}}\right)\right)$
with respect to the topology induced by the Euclidean norm on $\mR^k$.
\end{proposition}
First note that if $E \in \PF(\Omega)$, then $d_E = \left(E\left(\chi_{\{\omega_1\}}\right), \ldots,
E\left(\chi_{\{\omega_k\}}\right)\right) = \left(P(\omega_1),\ldots,P(\omega_k)\right)$, where 
$P$ is the underlying probability of $E$, and hence $\|E\| = \sqrt{P(\omega_1)^2 + \cdots 
+ P(\omega_k)^2}$. To complete the proof, we need some further statements first.

\begin{definition}
    A vector space $\mathcal{X}$ is called topological vector space if the topology on $\mathcal{X}$ is such that $(x,y) \mapsto x+y$ is continuous with respect to the product topology on $\mathcal{X} \times \mathcal{X}$ and $(\lambda,x) \mapsto \lambda x$ is continuous with respect to the product topology on $\mR \times \mathcal{X}$. We call a topology which makes $\mathcal{X}$ a topological vector space a \emph{linear topology}.
\end{definition}

\begin{remark}
The weak* topology makes $\linftyd$ a topological vector space with a Hausdorff topology,\footnote{See for instance Exercise 13 and 14 in \citet{taoweakstar}.} where vector addition and scalar product are defined pointwise: $Y+Z \coloneqq X \mapsto Y(X) + Z(X)$ $\forall X \in \linfty$, $\lambda Z \coloneqq X \mapsto \lambda Z(X)$ $\forall X \in \linfty$.
\end{remark}

\begin{remark}[well-known]
\label{remark:normtvs}
A vector space whose topology is induced by a norm is a topological vector space.
\end{remark}

\begin{proposition}
\label{prop:tvsonetopology}
On every finite dimensional vector space X there is a unique topological vector space structure. In other words, any two Hausdorff linear topologies on X coincide \citep{tvs3notes}.
\end{proposition}

Now Proposition~\ref{prop:equivofclusterpoints} directly follows.
\begin{proof}[Proof of Proposition~\ref{prop:equivofclusterpoints}]
Our norm $\|\cdot\|$ makes $\linftyd$ a topological vector space due to Remark~\ref{remark:normtvs}, and any topology induced by a norm is Hausdorff; but the weak* topology also makes $\linftyd$ a topological vector space, and the weak* topology is Hausdorff. Hence we can conclude from Proposition~\ref{prop:tvsonetopology} that they coincide. But then the two notions of what a cluster point is of course also coincide, since this depends just on the topology.
\end{proof}

Thus, for the proof of Theorem~\ref{theorem:converse}, we will work exclusively with the topology induced by the Euclidean metric
restricted to $\Delta^k$. For $z\in\Delta^k$ and $\epsilon>0$ define 
the \Def{$\epsilon$-neighbourhood (ball) }
\[ \Def{N_\epsilon(z)\coloneqq\left\{p\in\Delta^k\colon \|p-z\|<\epsilon\right\}} \] where
$\|\cdot\|$ is the Euclidean norm (restricted to the simplex).  Then  from
\citep[p.\@\xspace 430]{schechter1997handbook} we have an equivalent definition of 
a cluster point:
\begin{definition}
	Say that $z\in\Delta^k$ is a \Def{\emph{cluster point}} of a sequence
	$r\colon\naturals\rightarrow\Delta^k$ (and  denote by
    $\Def{\CP(r)}$ the \Def{\emph{set of all  cluster points of $r$}})
     if for all $\epsilon>0$, 
	$\left|\{n\in\naturals\colon r(n)\in N_\epsilon(z)\right|=\aleph_0$, where
	$\Def{\aleph_0\coloneqq|\naturals|}$ is the cardinality of the natural numbers. 
\end{definition}

\subsection{Further Notation}
We will work solely with the topology on $\Delta^k$ induced by the 
Euclidean metric;  by  the argument in the previous subsection 
the cluster points w.r.t. this topology coincide with 
those w.r.t. the weak* topology. 

We introduce further notation to assist in stating our algorithm. The $i$th
canonical unit vector in $\Delta^k$ is denoted 
$\Def{e_i\coloneqq(0,\ldots,1,\dots,0)}$,
where the 1 is in the $i$th position. The boundary of the simplex is 
\[ \Def{\partial\Delta^k\coloneqq\{(z_1,\ldots,z_k)
\colon z_1,\ldots,z_k\ge 0,\  z_1+\cdots+z_k=1\}.} \]
If $p_1,p_2\in\Delta^k$ then 
$\Def{l(p_1,p_2)\coloneqq\{\lambda p_1+(1-\lambda) p_2\colon \lambda\in [0,1]\}}$
is the \Def{\emph{line segment connecting $p_1$ and $p_2$}}.  If
$C\subset\Delta^k$ is a rectifiable closed curve parametrised by
$c\colon[0,1]\rightarrow\Delta^k$, its length is
$\Def{\operatorname{length}(C)=\int_0^1 |c'(t)| dt}$.  For $y\in\reals$, 
$\Def{\lfloor y\rceil}$ is the nearest integer to $y$: 
$\Def{\lfloor y\rceil\coloneqq\lfloor
y+\textstyle\frac{1}{2}\rfloor}$.  We apply certain operations elementwise.
For example, if $z=\langle z_1,\ldots,z_k\rangle\in\Delta^k$ and
$\iota\in\naturals^k$, then $\Def{\lfloor Tz\rceil\coloneqq\langle\lfloor
Tz_1\rceil,\ldots, \lfloor T z_k\rceil\rangle}$ and for $T\in\naturals$,
$\Def{\iota/T}\in\reals^n$ is simply $\langle\iota_1/T,\ldots,\iota_k/T\rangle$.
If $i<j\in \naturals$ the ``interval'' $\Def{[i,j]\coloneqq \{m\in
\naturals\colon i\le m \le j\}}$.

To avoid confusion, we will reserve ``sequence'' for the infinitely long 
$x\colon\naturals\rightarrow [k]$ and use ``\Def{segment}'' to denote finite 
length strings  $z\colon[n]\rightarrow [k]$ which we will write explicitly
as $\langle z_1,\ldots, z_n\rangle$.
We construct the sequence $x$ attaining the desired behavior of $r^x$ 
by iteratively appending a series of segments.
We denote the empty segment as $\langle\,\rangle$.
If $x^1$ and $x^2$ are two finite segments of lengths $\ell_1$ and $\ell_2$ then
is their \Def{concatenation} is the length $\ell_1+\ell_2$ segment
$\Def{x^1 x^2\coloneqq \left\langle x_1^1,\ldots,x_{\ell_1}^1,
x_1^2,\ldots, x_{\ell_2}^2\right\rangle}$. 
We extend the  $i^{[j]}$ notation to segments: if $z=\langle
z_1,\ldots,z_\ell\rangle$, then $\Def{z^{[\iota]}\coloneqq\langle{z,z \ldots,
z}\rangle} $ is the length $\ell\iota$ segment formed by concatenating $\iota$
copies of $z$.  Given $n\in\naturals$ and 
a sequence $x\colon\naturals\rightarrow [k]$, the \Def{shifted 
sequence $x^{+n}\colon\naturals\rightarrow [k]$} is defined via
	$\Def{x^{+n}(i)\coloneqq x(i+n)}$ for $i\in\naturals$.
 
\subsection{Properties of Relative Frequency Sequences}
Our construction of $x$ relies upon the following elementary
property of relative frequency sequences.
\begin{lemma}
	\label{lemma:r-decomposition}
	Suppose $k,n,m\in\naturals$, $x\colon\naturals\rightarrow [k]$. 
    Then 
	\begin{equation}
		r^x(n+m) = \displaystyle\frac{n}{n+m} r^x(n) +\frac{m}{n+m}
		r^{x^{+n}}(m).
		\label{eq:r-decomposition}
	\end{equation}
\end{lemma}
\begin{proof}
	For any $i\in[k]$ we have
	\begin{align*}
		r_i^x(n+m) &= \textstyle\frac{1}{n+m}|\{j\in[n+m]\colon
			x(j)=i\}|\\
		&= \textstyle\frac{1}{n+m}\left(|\{j\in[n]\colon x(j)=i\}|
			+ |\{j\in[n+m]\setminus[n]\colon x(j)=i\}|\right)\\
			&= \textstyle\frac{1}{n+m}\frac{n}{n} 
			    |\{j\in[n]\colon x(j)=i\}|
			    + \frac{1}{n+m}\frac{m}{m}
			    |\{t\in[m]\colon x(n+t)=i\}|\\
		&= \textstyle\frac{n}{n+m} r_i^x(n) +
			    \frac{m}{n+m} \frac{1}{m}|\{t\in[m]\colon
			    x^{+n}(t)=i\}|\\
			    & = \textstyle\frac{n}{n+m} r_i^x(n) +
			    \frac{m}{n+m} r_i^{x^{+n}}(m).
	\end{align*}
	Since this holds for all $i\in[k]$ we obtain 
	\eqref{eq:r-decomposition}.
\end{proof}
Observe that (\ref{eq:r-decomposition}) also holds when 
$x\colon[n+m]\rightarrow [k]$ is a segment, in which case 
$x^{+n}=\langle x_{n+1},\ldots,x_{n+m}\rangle$.
Furthermore note  that (\ref{eq:r-decomposition}) is a convex combination of the two
points $r^x(n)$ and $r^{x^{+n}}(m)$ in $\Delta^k$ since
$\frac{n}{n+m}+\frac{m}{n+m}=1$ and both coefficients are positive. These
two points are (respectively) the relative frequency of $x$ at $n$, and the relative
frequency of $x^{+n}$ at $m$. This latter sequence will be the piece ``added
on'' at each stage of our construction and forms the basis of our piecewise
linear construction of $r^x$ such that its cluster points are a given
$C\subset\Delta^k$.

The set of cluster points of any sequence is closed. In addition, we have

\begin{lemma} \label{lem:connected}
	For any $k\in\naturals$ and $x\colon\naturals\rightarrow [k]$,
	$\CP(r^x)$ is a connected set.
\end{lemma}
This follows immediately from \citep[Lemma 2.6]{bauschke2015}  upon observing that 
$\lim_{n\rightarrow\infty} \|r^x(n)-r^x(n+1)\|=0$ since $\|r^x(n)-r^x(n+1)\|
=\|r^x(n)-\frac{n}{n+1}r^x(n)-\frac{1}{n+1} e_{(x(n+1)}\|=
\frac{1}{n+1}\|r^x(n)-e_{x(n+1)}\|\le \frac{2}{n+1}$.
The boundedness  of $r^x$ 
is essential for this to hold ---  for unbounded sequences the set 
of cluster points need not be connected \citep{avsic1970limit}.
\subsection{Logic of the Construction}
\label{app:logicofconstruction}
\begin{figure}
\vspace*{-2.7cm}
\begin{center}
\def\svgwidth{0.5\textwidth}
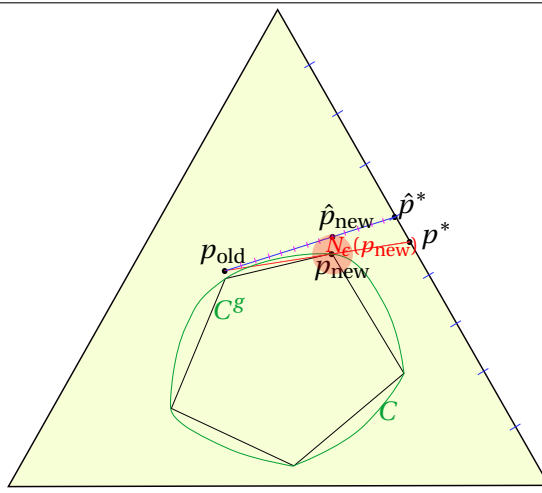
\end{center}
\vspace*{-4cm}
\caption{Illustration of construction of the sequence $x$. 
The figure shows how a single segment is 
created. We have the desired curve $C$ (in dark green) and a polygonal approximation 
$C^g$ using 5 segments (thus $V^g=5)$. We assume that we have already constructed
the first $n$ elements of $x$ and thus we can compute $r^x(n)$. 
Denote this by $p_{\mathrm{old}}$, an element of the simplex, itself shown in pale green.
We hope to append a segment $y$ of length $m$ to $x$ such $r^{xy}(n+m)=p_\mathrm{new}$.
However we can only choose elements of $y$ from $[k]$ and that means that at each step 
we move towards one of the vertices of the simplex. We note that ideally we would move 
from $p_\mathrm{old}$ towards $p^*\in\partial\Delta^k$.  In order to deal with the 
restrictions on directions we can head, we quantize $p^*$ as $\hat{p}^*=\lfloor T
p^*\rceil/T$ where in this example we have chosen $T=9$. As argued in the main text
this means we will thus be restricted to heading towards one of a fixed set of 
points on the boundary (marked in blue). Observe $\hat{p}^*$ is located at one 
of these points in the diagram, but in general it might not even be on the boundary 
of the simplex. 
We then construct $y$ to move $r^x$ towards $\hat{p}^*$ and take sufficient steps to 
move to $\hat{p}_{\mathrm{new}}$. This too is done in repeated steps by setting 
$y=z^{[\tilde{\ell}]}$ where $z$ is a shorter segment (marked by purple ticks
on the line segment $l(p_\mathrm{old},\hat{p}^*)$) which will move towards 
$\hat{p}^*$ by a small amount. The end result is that we get 
$r^{xy}(n+m)=\hat{p}_\mathrm{new}$ which 
is contained within $N_\epsilon(p_\mathrm{new})$, an $\epsilon$-ball centered
at $p_\mathrm{new}$.
\label{fig:construction}
}
\end{figure}

The idea of our construction is as follows (see Figure \ref{fig:construction}
below for a visual aid). In order to satisfy the
definition of cluster points, we need to return to each neighbourhood of
each point in $C$ infinitely often. To that end we iterate through an
infinite sequence of generations indexed by $g$.  For each $g$, we
approximate $C$ by a polygonal approximation $C^g$ comprising $V^g$
seperate segments. We
choose the sequence $(V^g)_{g\in\naturals}$ so that $C^g$ approaches $C$ in
an appropriate sense.  Then for generation $g$ we append elements to $x$ to
ensure the sequence of relative frequencies makes another cycle
approximately following $C^g$.  We control the approximation error of this
process and ensure its error is of a size that also decreases with
increasing $g$. 

We now describe the construction of a single generation. Thus
suppose $g$ is now fixed and suppose the current partial sequence (segment)
$x$ has length $n$. We suppose (and will argue this is ok later)
that $r^x(n)$ is  close to one of the vertices of $C^{g-1}$. We
then choose a finer approximation $C^g$ of $C$ (since
$V^g>V^{g-1}$). 

For each vertex $p_v^g$, $v\in[V^g]$ we append elements to $x$ resulting in
a segment of length $n'$. We do this in a manner such that we move the
relative frequency from $r^x(n)$ to $r^x(n')\approx p_v^g$.  We do so by
appending multiple copies of a vector $z$ to $x$ where $r^z(m)$ points in
the same direction as the direction one needs to go from $p_{\mathrm{old}}$
to $p_{\mathrm{new}}$.  This can only be done approximately because with a
finite length segment, the set of directions one can move the relative
frequencies is quantized. We choose the fineness of the quantization to be
fine enough to achieve the accuracy we need.  That is governed by the
parameter $T\in\naturals$. We then append $\tilde{\ell}$ copies of $z$ to $x$ where
$\tilde{\ell}$ is the integer closest to the real number $\ell$ that would
be the ideal number of steps needed to get to the desired point
$p_{\mathrm{new}}$. We also control the error incurred by approximating
$\ell$ by $\tilde{\ell}$. The upshot of this is that with the resulting
extension to $x$ we have $r^x(n')$ is sufficiently close to
$p_{\mathrm{new}}$. We then repeat this operation for all the vertices
$p_v^g$ for $v\in[V^g]$. This completes generation $g$. We show below that
for each generation $g$, \emph{all} the points in the relative frequency
sequence are adequately close to $C^g$, where ``adequately close'' is
quantified and increases in accuracy as $g$ increases.

\begin{algorithm}[t]
	\caption{Construction of $x$ such that $\CP(r^x)=C$\label{alg:C}}
\begin{algorithmic}[1]
	\Require $C\subset\Delta^k$, a rectifiable closed curve 
	 parametrized as $c\colon[0,1]\rightarrow\Delta^k$
	    \Require $V\colon\naturals\rightarrow\naturals $ 
		    \Comment{Number of segments at generation $g$;
		    as a function of $n$}
	    \Require $T\colon\naturals\rightarrow\naturals $
		    \Comment{Controls quantization of angle; needs to be increasing}
	    \State  $x\gets \langle 1\rangle$
	    \Comment{Arbitrary initialization $x_1=1$}
	    \State $p_{\mathrm{old}}\gets e_1$  
	    \Comment{$\pold=r^{\langle  1\rangle}(1)=e_1$}
	\State $n\gets 1$ \Comment{$n$ is always updated to correspond to
		the current length of $x$}
	\State $g \gets 1$
	\While{{true}} \Comment{Iterate over repeated generations
		$g$; $V$ is chosen at start of generation}
        	\State $V\gets V(g)$  
		\Comment{Choose $V$ for generation $g$}
	     \State $p_v \gets c(v/V) \ \mbox{for\ } v=0,\ldots,V$  
		     \Comment{Vertices of $C^g\coloneqq \bigcup_{v\in[V]}
		     l(p_{v-1},p_{v})$ }
	    \State $v \gets 0$  
	    \While{$v\le V$}  \Comment{For all vertices of $C^g$}
		    \State $T \gets T(n)$  \Comment{Quantization of angle;
			    chosen per segment}
		    \State $p_{\mathrm{new}} \gets
		    p_{v+1}$ \Comment{The next vertex of $C^g$}
		    \State $\gamma_i\gets
		    {p_{\mathrm{old},i}}/(p_{\mathrm{old},i}- p_{\mathrm{new},i}) 
			\mbox{\ for\ } i\in[k]$
			\Comment{Will have 
				$p_{\mathrm{old}}\approx p_v$}
		    \State $\gamma\gets\min\{\gamma_i\colon i\in[k],
			\gamma_i>0\}$ \Comment{See (\ref{eq:gamma-def})}
			\State $p^*\gets
			\gamma(p_{\mathrm{new}}-p_{\mathrm{old}})
			+ p_{\mathrm{old}}$
			    \Comment{Determine $p^*\in\partial\Delta^k$}
		    \State $\iota\gets \lfloor Tp^*\rceil$
			    \Comment{Elementwise;
			    $\iota=(\iota_1,\ldots,\iota_k)$}
		    \State $\tilde{p}^*\gets \iota/T $
			    \Comment{Elementwise; quantized version
			    of $p^*$}
		    \State $\tilde{T}\gets \sum_{i=1}^k \iota_i$ 
			    \Comment{Will have $\tilde{T}\approx T$}
		    \State $y\gets
			    \langle 1^{[\iota_1]},\ldots,k^{[\iota_k]}\rangle$
			    \Comment{The string $y$ is thus of length $\tilde{T}$}
		    \State $\tilde{\ell}\gets {\lceil\frac{n}{T(\gamma-1)}\rceil}$
		       \Comment{Integer number of repetitions of
		       $y$ needed}
		    \State $x\gets x\, y^{[\tilde{\ell}]}$
		    \Comment{Construct new $x$ by appending $z$, comprising
			    $\tilde{\ell}$ copies of $y$}
			    \State $n\gets n+\tilde{\ell}\tilde{T}$
			    \Comment{Length of $x$ now}
		    \State $p_{\mathrm{old}}\gets r^x(n)$
			    \Comment{Relative frequency at current $n$}
		\State $v\gets v+1$  \Comment{Move onto next vertex of $C^g$}
	    \EndWhile
	    \State $g\gets g+1$
	    \Comment{Move onto next generation of the construction}
	\EndWhile
	\Comment{Procedure never terminates}
\end{algorithmic}
\end{algorithm}

We consistently use the following terminology in describing our algorithm:
\begin{description}[nolistsep]
	\item[generation] These are indexed by $g$ and entail an entire
		pass around the curve $C$, or more precisely its polygonal
		approximation $C^g\coloneqq \bigcup_{v\in[V^g]}
		     l(p_{v-1},p_{v})$ 
	\item[segment]  Corresponds to a single line segment $
		l(p_{v-1},p_{v})$
		of the $g$th polygonal approximation.
	\item[piece] Corresponds to appending $z=\langle
		1^{[\iota_1]},\ldots,k^{[\iota_k]}\rangle $ to $x$, which 
		results in moving  $r^x(n)$ in the
		direction $\hat{p}^*$.
	\item[step] The appending of a single element of $z$, which will
		always move $r^x(n)$ towards one of the vertices of the
		simplex $e_i$ ($i\in[k]$).
\end{description}

The end result is that we have constructed a procedure (Algorithm
\ref{alg:C}) which runs
indefinitely ($g$ increases without bound), and which has the property that for any 
choice of $\epsilon>0$,
if one waits long enough, there will be a sufficiently large $g$ such that
all the relative frequencies associated with generated $g$ are within
$\epsilon_g$ of $C$, and $(\epsilon_g)_{g\in\naturals}$  is a null
sequence. We will thus conclude that $\CP(r)\supseteq C$.  We will 
also argue that $\CP(r)\subseteq C$ completing the proof.

\subsection{Construction of \texorpdfstring{$p^*$}{p-star} and its Approximation \texorpdfstring{$\tilde{p}^*$}{p-star tilde}}

The basic idea of the construction is to exploit Lemma
\ref{lemma:r-decomposition}. Suppose $n\in\naturals$ (and suppose it is
``large'') and fix $m=1$ in (\ref{eq:r-decomposition}) to obtain
\begin{equation}
\label{eq:r-one-step}
	r^x(n+1) = \textstyle\frac{n}{n+1} r^x(n) +\frac{1}{n+1}
	r^{x^{+n}}(1).
\end{equation}
Now $r^{x^{+n}}(1)=e_{x(n+1)}$ and so 
$r^x(n+1) = \frac{n}{n+1} r^x(n) +\frac{1}{n+1} e_{x(n+1)}$. When $n$ is
large $\frac{n}{n+1}\approx 1$ and $\frac{1}{n+1}$ is small, and so
this says that appending $x(n+1)$ to the length $n$ segment $x([n])$ moves
the relative frequency $r^x$ from $r^x(n)$ in the direction of $e_{x(n+1)}$
by a small amount. Observe that the \emph{only} directions which the point
$r^x(n)$ can be moved is towards one of the vertices of the $k$-simplex,
$e_1,\ldots, e_k$.  Thus if we had, for a fixed $n$ that
$r^x(n)=\pold$ and we wished to append $m$ additional elements
$\Def{z}$ to $x$ to produce $xz$ such that $r^{xz}(n+m)=\pnew$, we need to
figure out a way of heading in the direction
$d=\pnew-\pold$ when at each step we are constrained to move a small amount
to one of the vertices. The solution is to approximate the
direction $d$ by a quantized choice that can be obtained by an integer
number of elements of $[k]$.

Given arbitrary $\pold\ne\pnew\in\relint\Delta^k$, we define $p^*$ to be
the intercept by $\partial\Delta^k$ of the line segment starting at $\pold$
and passing through $\pnew$. 
(If $p_\mathrm{new}\in\partial\Delta^k$ set $p^*=p_\mathrm{new}$.)
The intercept on the boundary of $\Delta^k$ is denoted $p^*$ and is given by 
\begin{equation}\label{eq:p-star-def}
	\Def{p^*\coloneqq\gamma (\pnew-\pold)+\pold}
\end{equation}
for some $\gamma>0$. We can determine $\gamma$ as follows. The choice of
$\gamma$ can not take $p^*$ outside the simplex. Thus let $\gamma_i$
($i\in[k]$) satisfy $\gamma_i(\pnew_i-\pold_i)+\pold_i=0$. Thus 
$\Def{\gamma_i=\frac{\pold_i}{\pold_i-\pnew_i}}$.  Any $\gamma_i<0$ points in the
wrong direction and so we choose 
\begin{equation}
\label{eq:gamma-def}
	\Def{\gamma\coloneqq\min\{\gamma_i\colon i\in[k]\ \mbox{and}\ \gamma_i >0\}.}
\end{equation}
Such a choice of $\gamma$
guarantees that $p^*\in\partial\Delta^k$.
Observe that the requirement that $\gamma_i>0$ means the denominator in the
definition of $\gamma_i$ is positive and less than the numerator, and thus
all $\gamma_i$ which are positive exceed $1$, and consequently $\gamma>1 $.

We can now take $p^*$ to be the
direction we would like to move $r^x(n)$ towards. However our only control
action is to choose a sequence $z\in[k]^m$.  To that end we suppose we
quantize the vector $p^*$  so that it has rational components with denominator
$T\in\naturals$ (which will be strategically chosen henceforth). As we will
shortly show, this will allow us to move (approximately) towards $p^*$.
Thus let $\Def{\iota_i\coloneqq\lfloor Tp_i^*\rceil}$ for $i\in[k]$ and set
\begin{equation}
\label{eq:p-tilde-star-def}
	\Def{\tilde{p}^*\coloneqq
	\left(\frac{\iota_1}{T},\ldots,\frac{\iota_k}{T}\right).}
\end{equation}
Observe that $\tilde{p}^*$ is not guaranteed to be in $\Delta^k$ because there is no
guarantee that $\sum_{i=1}^k \tilde{p}^*_i=1$.  

We have 
\begin{lemma} 
	\label{lem:tilde-p-star-error}
 Let $p^*$  and $\tilde{p}^*$ be defined by (\ref{eq:p-star-def}) and 
 (\ref{eq:p-tilde-star-def}) respectively.
	Then for all $i\in [k]$, $\left|\tilde{p}_i^*-p_i^*\right|\le \frac{1}{2T}$.
\end{lemma}
\begin{proof}
	$
		\left|\tilde{p}_i^*-p_i^*\right|=\textstyle\frac{1}{T}\left|T\tilde{p}_i^*-Tp_i^*\right|
		= \textstyle\frac{1}{T}\left|\lfloor T{p}_i^*\rceil -T p_i^*\right|
		\le\textstyle\frac{1}{2T},
	$
by definition of the rounding operator $\lfloor\cdot\rceil$.
\end{proof}

\subsection{Determining the Number of Steps to Take}
Observe that (\ref{eq:r-decomposition}) can be written as
\begin{equation}
	r^x(n+m) = (1-\alpha) r^x(n) +\alpha r^{x^{+n}}(m).
\end{equation}
where $\alpha=\frac{m}{n+m}$ and thus $1-\alpha=\frac{n}{n+m}$. This
suggests that we can engineer the construction of $x$ by requiring a
suitable $\alpha$ such that 
\[
	(1-\alpha)\pold +\alpha p^*=\pnew.
\]
Recall we want the sequence $r^x$ to move from $\pold$ to $\pnew$ which can be achieved by a taking a suitable comvex combination of $\pold$ and $p^*$, which corresponds to appending a suitable number of copies of $z$ to $x$, where $z$ is chosen to move $r^x(n)$ in the direction of $p^*$.
If we substitute the definition of $p^*$ from (\ref{eq:p-star-def}) we
obtain the problem:
\begin{align*}
	& \mbox{\ Find\ } \alpha\mbox{\ such that\ }
	(1-\alpha)\pold +\alpha[\gamma(\pnew-\pold)+\pold]=\pnew\\
\Leftrightarrow &  \mbox{\  Find\ } \alpha\mbox{\ such that\ } (1-\alpha)\pold +\alpha\gamma\pnew -\alpha\gamma\pold
	+\alpha\pold -\pnew =0_k\in\reals^k\\
\Leftrightarrow & \mbox{\ Find\ } \alpha\mbox{\ such that\ } \pold(1-\alpha\gamma) +\pnew(\alpha\gamma-1)=0_k,\\
\end{align*}
which is only true when either $\pold=\pnew$ (which is a trivial case) or
when $1-\alpha\gamma=0$ and thus $\Def{\alpha\coloneqq 1/\gamma}$, which we take as a definition.
Since $\gamma>1$ this implies $\alpha<1$, which is consistent with our original motivation 
for taking convex combinations.

We will append the segment $z$ to $x$, where $z=y^{[\ell]}$ and 
$y=\langle 1^{[\iota_1]},\ldots,k^{[\iota_k]}\rangle$.
Now if each $y$  is of length $T$ and we
notionally made $\ell$ repetitions, we would have $m=\ell T$ 
From the
definition of $\alpha$ this means 
\begin{equation}
\label{eq:alpha-def}
{\alpha=\frac{\ell T}{n+\ell T}}.
\end{equation}
We presume $n$ is given (at a particular stage of construction)
and $T\in\naturals$ is a
fixed design parameter. We can thus solve for $\ell$ to obtain
\[
	\ell= \frac{\alpha n}{T(1-\alpha)}=\frac{(1/\gamma)
	n}{T(1-1/\gamma)}=\frac{n}{T(\gamma-1)}.
\]
Observe that $\ell$ is not guaranteed to be an integer, a complication we will deal with later.
If it was an integer, we would
create $z$ by concatenating $\ell$ copies of $y$ which is of length
$T$. The vector $y$ moves $\pold$ towards $p^*$. By appending $\ell$
copies we should move $r^x$ to $\pnew$ as desired.

However we do not head exactly in the direction of $p^*$, since we worked
with a quantized version $\tilde{p}^*$ instead, and we can not always take 
$\ell$ copies because $\ell$ is not guaranteed to be an integer; instead we will take 
\[
	\Def{\tilde{\ell}\coloneqq
	\left\lceil\frac{\alpha n}{T(1-\alpha)}\right\rceil=
	\left\lceil\frac{n}{T(\gamma -1)}\right\rceil}
\]
copies of $y$ which will move $r^x(n)$ towards $\hat{p}^*$ instead of $p^*$.
We now proceed to analyse the effect of these approximations on our
construction.

\subsection{Analyzing the Effect of Approximations}

The parameter $T\in\naturals$ is a design variable.  
Our construction will utilize 
\begin{equation}
\label{eq:tilde-T-def}
\Def{\tilde{T}\coloneqq\sum_{j=1}^k\iota_j},
\end{equation}
 where $\iota_j=\lfloor Tp_j^*\rceil$ for $j\in[k]$.
Then $\tilde{T}\approx T$, a claim which we quantify below.
\begin{lemma}\label{lem:T-tilde-T}
Suppose $T\in\naturals$, 
and $\tilde{T}$ is defined as above. Then
	 $T-\frac{k}{2}\le \tilde{T}\le T+\frac{k}{2}$.
\end{lemma}
\begin{proof}
By definition of
the rounding operator $\lfloor\cdot\rceil$, we have that 
\[
	\left|\iota_j - \lfloor Tp_j^*\rceil\right| \le \textstyle\frac{1}{2} \ \ \
	\forall j\in[k].
\]
Thus
\begin{align*}
	& Tp_j^* -\frac{1}{2} \le \iota_j \le T p_j^*
	+\frac{1}{2}\
	\ \ \ \forall j\in[k]\\
\Rightarrow\ \   & \sum_{j=1}^k\left(Tp_j^*-\frac{1}{2}\right) \le
	\sum_{j=1}^k \iota_j \le \sum_{j=1}^k
		\left(Tp_j^*+\frac{1}{2}\right)\\
\Rightarrow\ \  & T-\frac{k}{2}\le \tilde{T}\le T+\frac{k}{2}.
\end{align*}
\end{proof}

Ideally we move $r^x(n)$ to 
\begin{equation}
	\label{eq:pnew}
	\pnew=(1-\alpha)\pold+\alpha p^*.
\end{equation}
by appending $z$ (i.e.~ we hope that $r^{xz}(n+m)=\pnew$).
But in fact the segment $z$ which we will append  to $x$ will move $r^x$ from
$\pold$ instead to
\begin{equation}
	\label{eq:phatnew}
	\Def{\hpnew\coloneqq(1-\tilde{\alpha})\pold+\tilde{\alpha}\hat{p}^*},
\end{equation}
where
\begin{equation}
\label{eq:tilde-alpha-def}
    \Def{\tilde{\alpha}\coloneqq \frac{\tilde{\ell}T}{\tilde{\ell}T+m}.}
\end{equation}
and
\begin{equation}
\label{eq:hat-p-star-def}
   \Def{\hat{p}^*\coloneqq\left(\frac{\iota_1}{\tilde{T}},\ldots,\frac{\iota_k}{\tilde{T}}\right).}
\end{equation}
We now determine the error incurred from these approximations. 
We first need the following Lemma
\begin{lemma}\label{lem:tilde-alpha-minus-alpha}
Suppose $n,T\in\naturals$, $\alpha$ is defined by (\ref{eq:alpha-def})
and $\tilde{\alpha}$ is defined by
(\ref{eq:tilde-alpha-def}). Then
	\[
		\left|\alpha - \tilde{\alpha}\right| \le \frac{T}{n}.
	\]
\end{lemma}
\begin{proof}
	By definition of $\tilde{\alpha}$ we have 
	\begin{align*}
		|\alpha-\tilde{\alpha}| & = \left|
			\frac{\tilde{\ell}T}{\tilde{\ell} T +n}   
		         - \frac{\ell T}{\ell T +n} \right|\\
	        &= \left|\frac{\tilde{\ell} T (\ell T+n) -\ell
		     T(\tilde{\ell} T+n)}{(\tilde{\ell} T +n)(\ell T+n)} \right|\\
	        &=\left|\frac{\tilde{\ell}Tn - \ell Tn}{(\tilde{\ell} T +n)(\ell T+n)}
		        \right|.\\
			\intertext{Since $(\tilde{\ell} T +n)(\ell T+n)>0$
			and $\tilde{\ell}=\lceil\ell\rceil\ge\ell$, we have}
		(\tilde{\ell} T +n)(\ell T+n) & \ge (\ell T+n)(\ell T+n)\ge n^2\\
		\intertext{and since $|\tilde{\ell}-\ell|\le 1$,}
		|\alpha-\tilde{\alpha}| & = \frac{|\tilde{\ell}-\ell|\cdot Tn}{n^2}
		\le \frac{T}{n}.
	\end{align*}
\end{proof}
Our construction does not move $r^x(n)$ towards $\tilde{p}^*$ but an
approximation of it, namely
$\hat{p}^*$  defined in (\ref{eq:hat-p-star-def}).
We exploit the fact that repeating a segment does not change its relative frequencies, 
which we state formally as
\begin{lemma}
	Let $z=\left\langle 1^{[\iota_1]},\ldots,k^{[\iota_k]}\right\rangle$. Then	
	$r^z(\tilde{T})=\hat{p}^*=r^{z^{[\tilde{\ell}]}}(\tilde{\ell}\tilde{T})$.
\end{lemma}
\begin{proof}
	For any $i\in[k]$ we have
	$r_i^z(\tilde{T})=\frac{1}{\tilde{T}}\left|\left\{j\in[\tilde{T}]\colon
z_j=i\right\}\right|=\frac{1}{\tilde{T}} \iota_i$. 
The first equality is immediate. For
the second, similarly we have 
$r_i^{z^{[\tilde{\ell}]}}(\tilde{\ell}\tilde{T})=
\frac{1}{\tilde{\ell}\tilde{T}}\left|\left\{j\in[\tilde{\ell}\tilde{T}]\colon
z_j^{[\tilde{\ell}]}=i\right\}\right|=\frac{1}{\tilde{\ell}\tilde{T}}\cdot 
\tilde{\ell} \iota_i = \frac{1}{\tilde{T}}\iota_i$ by definition of
$z^{[\tilde{\ell}]}$.
\end{proof}
Since
$\tilde{p}^*=\left(\frac{\iota_1}{{T}},\ldots,\frac{\iota_k}{{T}}\right)$
we have that $\tilde{p}^*=\frac{\tilde{T}}{T}\hat{p}^*$. This allows us to
show:
\begin{lemma}\label{lem:hat-p-tilde-p}
Suppose $T\in\naturals$ and $\hat{p}^*$ is defined via 
(\ref{eq:hat-p-star-def}).  Then
	$\left\|\hat{p}^*-\tilde{p}^*\right\|\le \frac{k}{2T}$.
\end{lemma}
\begin{proof}
	We have
	$
		\left\|\hat{p}^*-\tilde{p}^*\right\|  =\left\|\hat{p}^* -
		\frac{\tilde{T}}{T} \hat{p}^*\right\|
		=\left|1-\frac{\tilde{T}}{T}\right|\cdot\|\hat{p}^*\|
		\le \left|1-\frac{\tilde{T}}{T}\right|.
	$
	Suppose $\tilde{T}<T$, then $1-\frac{\tilde{T}}{T}>0$ and
	$
		\left|1-\frac{\tilde{T}}{T}\right|=1-\frac{\tilde{T}}{T}\le
		1-\frac{T-k/2}{T} = \frac{k}{2T}
	$
	by Lemma \ref{lem:T-tilde-T}.  Similarly if $\tilde{T}>T$, then 
	$1-\frac{\tilde{T}}{T}<0$ and
	$
		\left|1-\frac{\tilde{T}}{T}\right|=\frac{\tilde{T}}{T}-1\le
		\frac{T+k/2}{T}-1 = \frac{k}{2T}
	$
	completing the proof.
\end{proof}

\begin{lemma}
	\label{lem:hpnew-pnew}
 Suppose $k,T\in\naturals$ and $\pnew$ and $\hpnew$ are defined as above. Then
	\begin{equation}
		\|\hpnew-\pnew\|\le \frac{4 T}{n}
		+\frac{k}{T}.
		\label{eq:hpnew-pnew}
	\end{equation}
\end{lemma}
\begin{proof}
	From (\ref{eq:pnew}) and (\ref{eq:phatnew}) we have
	\begin{align}
		\|\hpnew-\pnew\| &= \|(1-\tilde{\alpha})\pold
		+\tilde{\alpha}\hat{p}^* -(1-\alpha)\pold -\alpha p^*\|\nonumber\\
		&=\|[(1-\tilde{\alpha})-(1-\alpha)]\pold+ 
		     (\tilde{\alpha}\hat{p}^*-\alpha p^*)\|\nonumber\\
		& \le \|(\alpha-\tilde{\alpha})\pold\| +
			\|\tilde{\alpha}\hat{p}^* -\alpha p^*\|\nonumber\\
			&\le \sqrt{2} |\tilde{\alpha}-\alpha| +
				\|\tilde{\alpha}\hat{p}^* -\alpha p^*\|.\label{eq:two-term-expression}
	\end{align}
	The second term in (\ref{eq:two-term-expression}) can be bounded as follows:
	\begin{align}
	     \|\tilde{\alpha}\hat{p}^* -\alpha p^*\| &=
		     \|(\tilde{\alpha}-\alpha+\alpha)\hat{p}^* 
		     -\alpha p^*\|\nonumber\\
	     &=\|(\tilde{\alpha}-\alpha)\hat{p}^* +
		     (\alpha\hat{p}^* -\alpha p^*)\|\nonumber\\
	     &\le \|(\tilde{\alpha}-\alpha)\hat{p}^*\| +
			     \|\alpha\hat{p}^* -\alpha p^*\|\nonumber\\
	     &=|\tilde{\alpha}-\alpha|\cdot\|\hat{p}^*\|
			     +\alpha\|\hat{p}^*-p^*\|\nonumber\\
	     &=|\tilde{\alpha}-\alpha|\cdot\|\hat{p}^*\|
		     +\alpha\|(\hat{p}^* - \tilde{p}^*)+(\tilde{p}^* -p^*)\|\nonumber\\
	     &\le|\tilde{\alpha}-\alpha|\sqrt{2}
		     +\alpha\|\hat{p}^* - \tilde{p}^*\|+\alpha\|\tilde{p}^* -p^*\|\nonumber\\
	     &\le  \sqrt{2}|\tilde{\alpha}-\alpha|  
		             +\frac{k}{2T} +\alpha\left(
			     \textstyle\sum_{i=1}^k \left(\tilde{p}_i^*-p_i^*\right)^2
			     \right)^{1/2},\nonumber\\
	     \intertext{by Lemma \ref{lem:hat-p-tilde-p} and the fact that
				     $\|\hat{p}^*\|\le 1$,}
	     &\le  \sqrt{2}|\tilde{\alpha}-\alpha| + \frac{k}{2T} +
			     \alpha \left(\textstyle\sum_{i=1}^k
			     \left(\frac{1}{2T}\right)^2\right)^{1/2}\nonumber\\
	     &= \sqrt{2}|\tilde{\alpha}-\alpha| +\frac{k}{2T}+ \alpha
			     \frac{\sqrt{k}}{2T}, \label{eq:intermediate-bound}
	\end{align}
	where we used Lemma \ref{lem:tilde-p-star-error} in the
	penultimate step.  Since $\alpha\le 1$, combining 
    (\ref{eq:two-term-expression}) and (\ref{eq:intermediate-bound})  we have
	\[
		\|\hpnew-\pnew\|\le 2\sqrt{2}\left|\tilde{\alpha}-\alpha\right| +
		\frac{k}{2T}+
		\frac{\sqrt{k}}{2T} \le  4\left|\tilde{\alpha}-\alpha\right| +
		\frac{k}{T}.
	\]
	Appealing to Lemma \ref{lem:tilde-alpha-minus-alpha} gives us
	(\ref{eq:hpnew-pnew}).
\end{proof}

The above arguments control the errors at the end of a \emph{piece} (and thus
in a \emph{segment}). But for later purposes we need control at each \emph{step}.  This
follows immediately by the fact that we make small steps:
\begin{lemma}\label{lem:r-within-segment}
	For $n\in\naturals$ and $m\in[\tilde{T}]$ and any
	$x\colon\naturals\rightarrow[k]$,
	\[
		\left\|r^x(n)-r^x(n+m)\right\|\le \frac{2T+k}{n}.
	\]
\end{lemma}
\begin{proof}
	By Lemma \ref{lemma:r-decomposition},
	\begin{align*}
	\left\|r^x(n)-r^x(n+m)\right\| &=\left\|r^x(n)-\textstyle\frac{n}{n+m}r^x(n)
		-\frac{m}{n+m}r^{x^{+n}}\!(m)\right\|\\
		&=\left\|(1-\textstyle\frac{n}{n+m}) r^x(n) -
		\frac{m}{n+m}
		r^{x^{+n}}\!(m)\right\|\\
		&=\frac{m}{n+m}\left\|r^x(n)-r^{x^{+n}}\!(m)\right\|\\
		&\le \frac{2m}{n+m}\\
		&\le \frac{2m}{n}\\
		&\le \frac{2T+k}{n},
	\end{align*}
	where the first inequality holds since
	$\|r^x(n)\|,\|r^{x^{+n}}\!(m)\|\le 1$ and 
	the last step follows from Lemma \ref{lem:T-tilde-T}.
\end{proof}
\subsection{Completing the Proof of Theorem \ref{th:CP-r-C}}

The remaining piece of the argument concerns the piecewise linear
approximation of $C$ by $C^g$. 
For $A,B\subseteq\Delta^k$ and $a\in\Delta^k$ define
$\Def{d(a,B)\coloneqq\min\{\|a-b\|\colon b\in B\}}$ and the \Def{Hausdorff distance}
\begin{equation}
\label{eq:hausrdorff-distance-def}
\Def{d(A,B)\coloneqq\max\{d(a,B)\colon a\in A\}}. 
\end{equation}
Let us write $V^g$ (the number of vertices in the piecewise 
linear approximation at generation $g$) in functional form as
$\Def{V(g)}$. Let $\Def{n(g)}$ denote the length of the segment 
of $x$ that has been constructed at the beginning of 
generation $g$, and let 
$\Def{g(n)\coloneqq \inf\{g\in\naturals\colon n\le n(g)\}}$
denote its quasi-inverse. Clearly $n(g)$ is strictly 
increasing in $g$ and $g(n)$ is increasing, but often constant.
With these definitions, we have $V=V(g)=V(g(n))$.

Denote by $\Def{\tilde{C}(V)}$ the best piecewise linear approximation 
of  $C$ with $V$ vertices, in the sense of minimizing 
$\Def{\psi_C(V)\coloneqq d(C,\tilde{C}(V))}$.
Since every rectifiable curve  $C$ has a Lipschitz continuous 
parametrisation, we have that $V\mapsto \psi_C(V)$ is decreasing in $V$
and $\lim_{V\rightarrow\infty} \psi_C(V)=0$.  Thus 
$\lim_{g\rightarrow\infty} \psi_C(V(g))=0$ and 
$\lim_{n\rightarrow\infty} \psi_C(V(g(n)))=0$, although the convergence could be 
very slow (in $n$) and its speed will depend on the choice of $C$.
Denote by  $\Def{\bar{C}(n)\coloneqq \tilde{C}(V(g(n)))}$
the sequence of best possible piecewise linear approximations  
of $C$ indexed by $n$, and let $\Def{\bar{\psi}_C(n)\coloneqq d(C,\bar{C}(n))}$.
We have thus shown:
\begin{lemma}\label{lem:bar-psi}
Let $C\subset\Delta^k$ be a rectifiable curve.  Then
\[
    \lim_{n\rightarrow\infty}  \bar{\psi}_C(n)=0.
\]
\end{lemma}

We summarize what we know so far.
\begin{enumerate}
	\item For all generations $g$, 
		$\|\hpnew-\pnew\|\le \frac{4T}{n}+\frac{k}{T}$, where
		$\pnew=p_v$ for $v\in[V^g]$ (Lemma \ref{lem:hpnew-pnew}). This means the
		following. Suppose at the beginning of segment $v$ in
		generation $g$ we have $n=\operatorname{length}(x)$.  By
		definition, we have $r^x(n)=\pold$ and
		$r^{xz^{[\tilde{\ell}]}}(n+\tilde{\ell}\tilde{T})=\hpnew$.
		Furthermore, for $m=i\tilde{T}$, $i\in[\tilde{\ell}]$  we
		have $r^{xz^{[i]}}(n+i\tilde{T})\in l(\pold,\hpnew)$.
	\item  Furthermore, (by Lemma \ref{lem:r-within-segment}) 
		for all $j\in[\tilde{\ell}\tilde{T}]$, 
		$d\left(r^{xz^{[\tilde{\ell}]}}(j), l(\pold,\hpnew)\right)\le
		\frac{2T+k}{n}$ --- the relative frequencies for all
		points in the segment are close to the line segment
		$l(\pold,\hpnew)$.
\end{enumerate}
Combining these facts, and appealing to the triangle inequality, we conclude
that  the sequence $x$ constructed by Algorithm 1 satisfies
\begin{equation}
	\label{eq:overall-error-abstract}
	d(r^x(n),C)\le \frac{4T}{n} +\frac{k}{T}+\frac{2T+k}{n}
	+\bar{\psi}_C(n) \ \ \ \forall n\in\naturals.
\end{equation}
We now choose $T=T(n)$ and $V=V(g(n))$ appropriately. One choice is to choose
$T(n)=\sqrt{n}$.
Equation \ref{eq:overall-error-abstract} then implies
\[
	d(r^x(n),C)\le \frac{4\sqrt{n}}{n}
	+\frac{k}{\sqrt{n}}+\frac{2\sqrt{n}+k}{n}
	  + \bar{\psi}_C(n) \ \ \ \forall n\in\naturals
\]
which implies that for all $n\in\naturals$ with $\sqrt{n}> k$,
\[
	d(r^x(n),C)\le \frac{4}{\sqrt{n}}
	+\frac{k}{\sqrt{n}}+\frac{4}{\sqrt{n}}
	+\bar{\psi}_C(n)=
 \frac{8+k}{\sqrt{n}} +\bar{\psi}_C(n),
\]
and thus Lemma \ref{lem:bar-psi} implies $\lim_{n\rightarrow\infty} d(r^x(n),C)=0$.
Furthermore, by the generational nature of our construction,
$r^x(n)$ repeatedly approximately follows the curve $C$. In each generation
it gets closer and closer. Thus for any $\epsilon$, for any point $y\in C$,
the sequence $r^x$ enters an $\epsilon$-ball of $y$ infinitely often. Thus
$\CP(r)\supseteq C$.

Finally, since in each generation the bounds above constrain $r^x$ more and
more tightly, there cannot exist cluster points that are not in $C$; that
is, $\CP(r)\subseteq C$.  We have thus proved Theorem \ref{th:CP-r-C}.

\subsection{Remarks on the Construction}
We make a few remarks on the construction.
\label{app:remarksonconstruction}
\begin{enumerate}[nolistsep]
	\item By the definition of $\tilde{\ell}$ we are guaranteed that
		$\tilde{\ell}\ge 1$ for each segment and each generation.
		For the construction to approximate well,  we need
		$\tilde{\ell}\gg 1$, which it will inevitably be when $n$
		gets large enough.
	\item Recall $n(g)$ is the length of $x$ at the beginning of
		generation $g$.  
      Let $L_C^g\coloneqq\operatorname{length}(C^g)>0$. Since each step
      which $r^x$ moves is of size less than $1/n(g)$ (see Equation 
      \ref{eq:r-one-step}), we require at least $L_C^g\cdot n(g)$ steps 
      for $r^x$ to traverse the whole of $C$ in generation $g$.  
      Thus at the end of generation $g$ and the beginning of generation $g+1$ 
      we have \[ n(g+1)\ge n(g)+ L_C^g\cdot n(g) =\lambda_C^g\cdot n(g).\]
      Furthermore, $L_C^g$ is increasing in $g$ and approaches 
      $\operatorname{length}(C)$.  Thus for all $g$, $\lambda_C^g>1$
      (and is in fact increasing in $g$).
      Hence $n(g)$ grows exponentially with $g$.
	\item The growth of the length of $x$ is controlled in a complex
		manner by the nature of the curve $C$. In particular if $C$
		is very complex, then $\bar{\psi}_C$ must decay slowly.
		Furthermore, if $C$ has parts close to $\partial\Delta^k$,
		in particular if some vertices $v$ of the piecewise linear
		approximation $C^g$ are close to $\partial\Delta^k$, then
		$\tilde{\ell}$ can end up very large for that segment,
		meaning that the length of the sequence $x$ grows more rapidly. See
		Lemma \ref{lem:extreme-case} for an illustration of this
		observation.
\end{enumerate}

Finally note that since (by Lemma \ref{lem:connected}) $\CP(r^x$) must always
be connected, we have in Theorem \ref{th:CP-r-C} what appears to be
the most general result possible
(under the restriction that $x$ takes values only in a finite set $[k]$). We
do not know what the appropriate generalization is to sequences $x$ that can
take values in an infinite (or uncountable) set\footnote{
Although we do not pursue this in any detail, we remark that one could design
an algorithm to construct $x$ such that $\CP(r^x)$ is \emph{any} subset 
$D\subseteq\Delta^k$  by using  our algorithm as a subroutine. The idea would 
be to construct a space filling curve that fills $D$, each generation of 
which is a rectifiable curve. One would appeal to our algorithm for each 
generation, and then once within a suitable tolerance, change the target 
to be the next generation of the space filling curve.  
A suitable method would be to simply intersect an extant families of 
closed space filling curves for $k$ dimensional cubes with the 
$(k-1)$-simplex (e.g. generalisations of the Moore curve),
attaching joins on $\partial D$ where necessary.   Since (see the main body
of the paper) it ends up being only the convex hull of $\CP(r^x)$ that matters, 
such an exotic construction is of little direct interest.}.

\subsection{From Boundary to Curves}
\label{app:fromboundarytocurves}
\begin{proof}[Proof of Corollary~\ref{cor:CP-convexbody}]

\citet{bronshteyn1975approximation} show  that if
$D$ is a convex set contained in the unit ball (w.r.t.\@\xspace the Euclidean norm) in $\reals^n$ and 
$\epsilon<10^{-3}$, then there exists a set of at most $K_\epsilon \coloneqq 3\sqrt{n}(9/\epsilon)^{(n-1)/2} $ points 
whose convex hull is at most $\epsilon$ away from $D$.   
Thus for any convex $D\subseteq \Delta^k$, and any $\epsilon<10^{-3}$ there is a 
polyhedron $D_\epsilon$ comprising the convex 
hull of $K_\epsilon<\infty$ points $Q_\epsilon\coloneqq \{q_i\in\Delta^k\colon i\in[K_\epsilon]\}$  
such that $d(D,D_\epsilon)<\epsilon$, where 
$d$ is the Hausdorff distance (\ref{eq:hausrdorff-distance-def}).  Hence for 
any $\epsilon<10^{-3}$ 
there exists a closed rectifiable curve $C_\epsilon$ (constructed by linearly 
connecting  successive points in $Q_\epsilon$) such that $d(\co C_\epsilon,D)<\epsilon$. 
One can then construct a sequence $x$ such that $\CP(r^x)=\partial D$ as follows. 
Start with $\epsilon_0=10^{-3}$. Pick and construct $C_{\epsilon_0}$ as above. 
Construct $x$ according to the previous procedure (Appendix~\ref{app:logicofconstruction} - \ref{app:remarksonconstruction}) to get an entire generation
of $r^x$ within $\epsilon_0$ of $C_{\epsilon_0}$.  Then let $\epsilon_{1}=\epsilon_0/2$ 
and repeat the procedure, appending the constructed sequence. 
Continue iterating (dividing the $\epsilon$ in half each phase) and one 
achieves that in the limit $\CP(r^x)=\partial D$.
\end{proof}

\subsection{Illustration}
\label{subsec:illustration}

We illustrate our construction for $k=3$ for a two (identical) generations
and thus a single piecewise linear (in fact polygonal) approximation $C^g$ of $C$. 
We take for $C$ (a scaled version of) the lemniscate of Bernoulli 
\citep[Chapter 12]{lockwood1961},
\citep[Section 5.3]{lawrence1972} mapped onto the
2-simplex as  the space curve
$\{(z_1(t),z_2(t),z_3(t)\colon t\in[0,2\pi]\}$, where
\begin{align*}
	z_1(t)&=\frac{1}{3}+\frac{1}{12}
	\frac{2\cos(t)}{1+\sin^2(t)}\\
	z_1(t)&=\frac{1}{3}+\frac{1}{12}
	\frac{2\sin(t)\cos(t)}{1+\sin^2(t)}\\
	z_3(t)&=1-z_1(t)-z_2(t).
\end{align*}

\begin{figure}
	\centering
	\begin{subfigure}[b]{0.40\textwidth}
		\centering
		\includegraphics[width=0.99\textwidth]{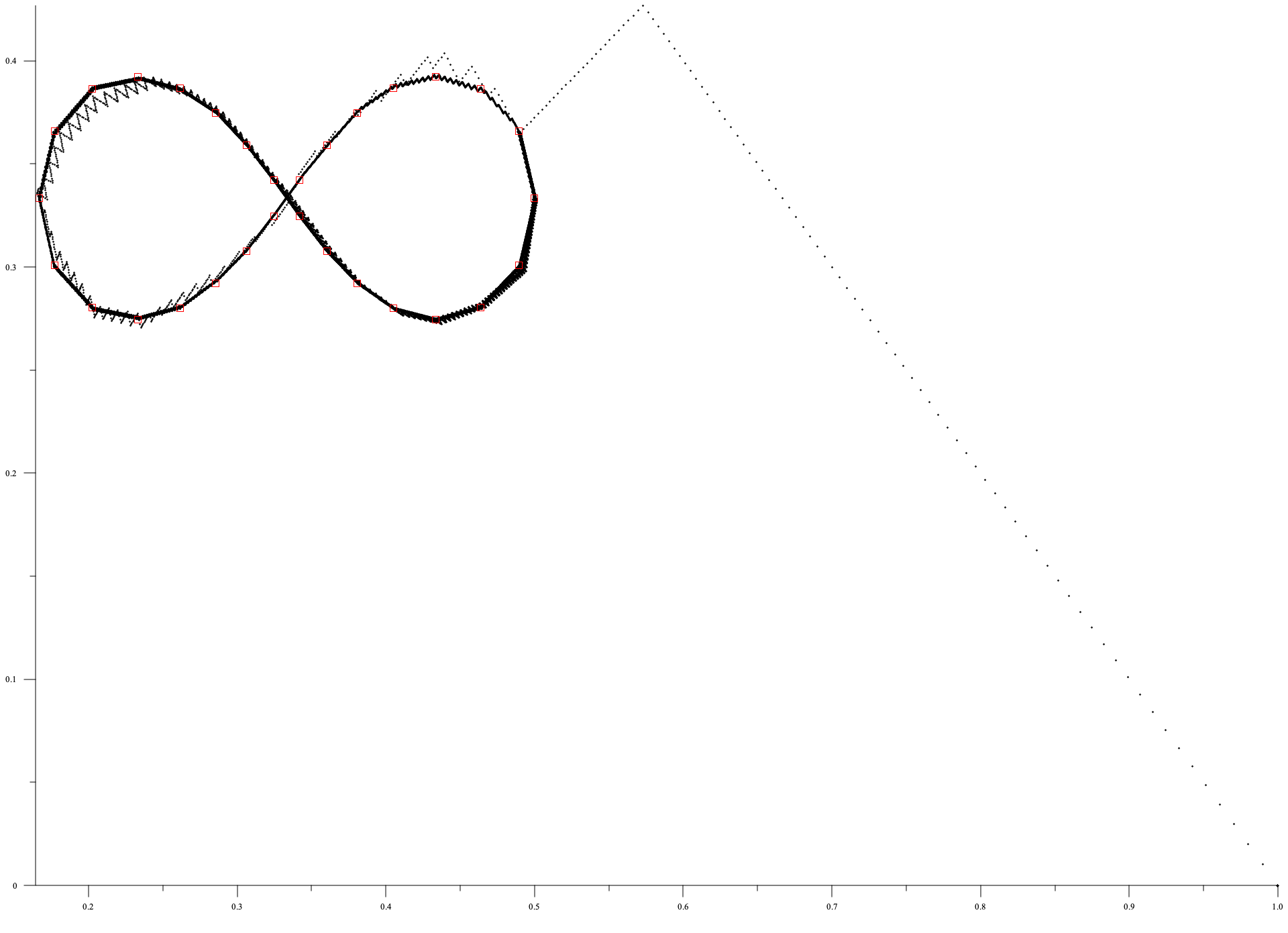}
		\caption{Overall view, showing
		initialisation.}
		\label{fig:lem-0}
	\end{subfigure}
	\begin{subfigure}[b]{0.40\textwidth}
		\centering
		\includegraphics[width=.99\textwidth]{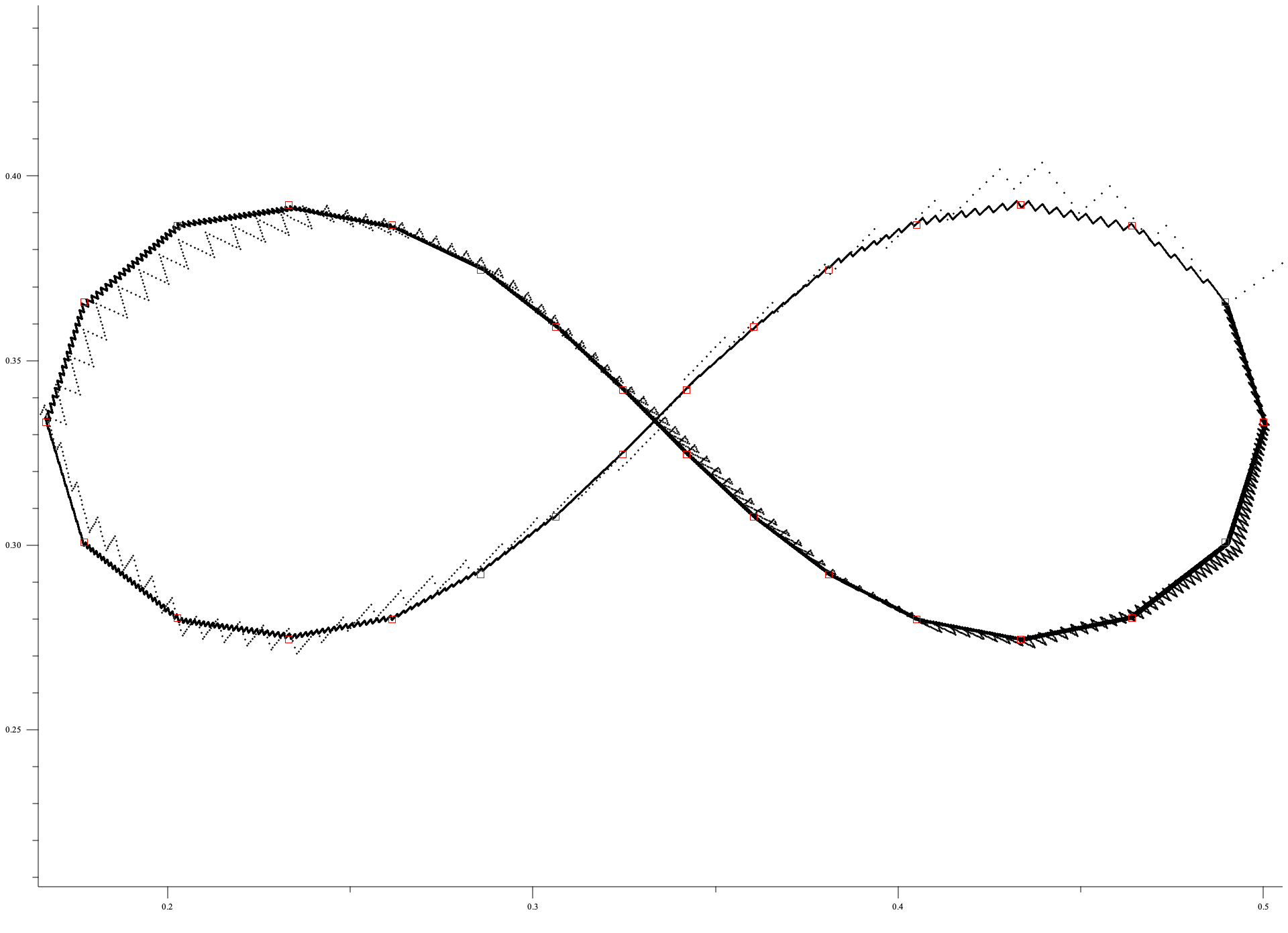}
		\caption{Closer view of just the region containing $C$.}
		\label{fig:lem-1}
	\end{subfigure}
	\begin{subfigure}[b]{0.40\textwidth}
		\centering
		\includegraphics[width=.99\textwidth]{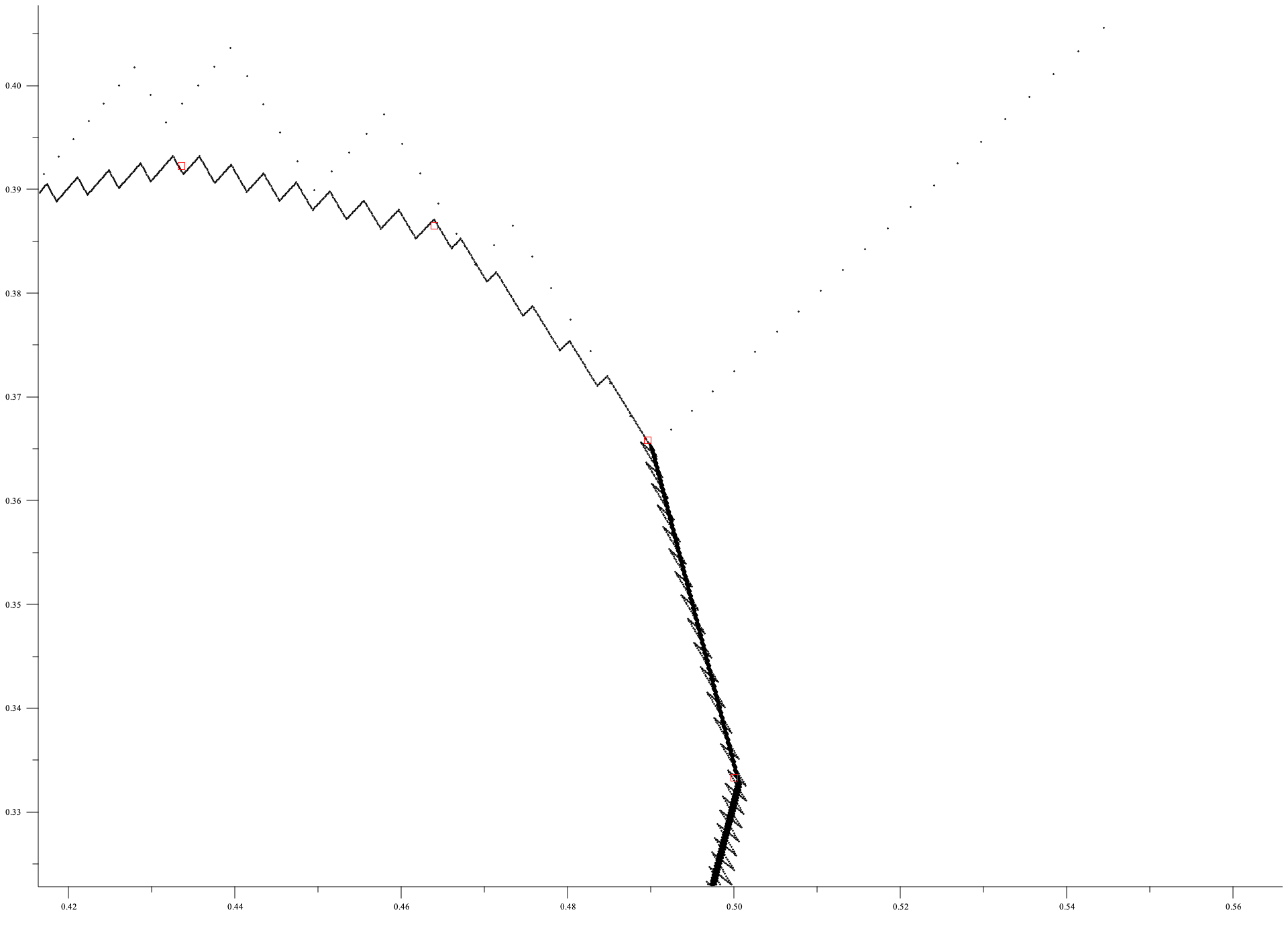}
		\caption{Zoom of upper right corner of Figure
			\ref{fig:lem-1}.}
		\label{fig:lem-2}
	\end{subfigure}
	\begin{subfigure}[b]{0.40\textwidth}
		\centering
		\includegraphics[width=.99\textwidth]{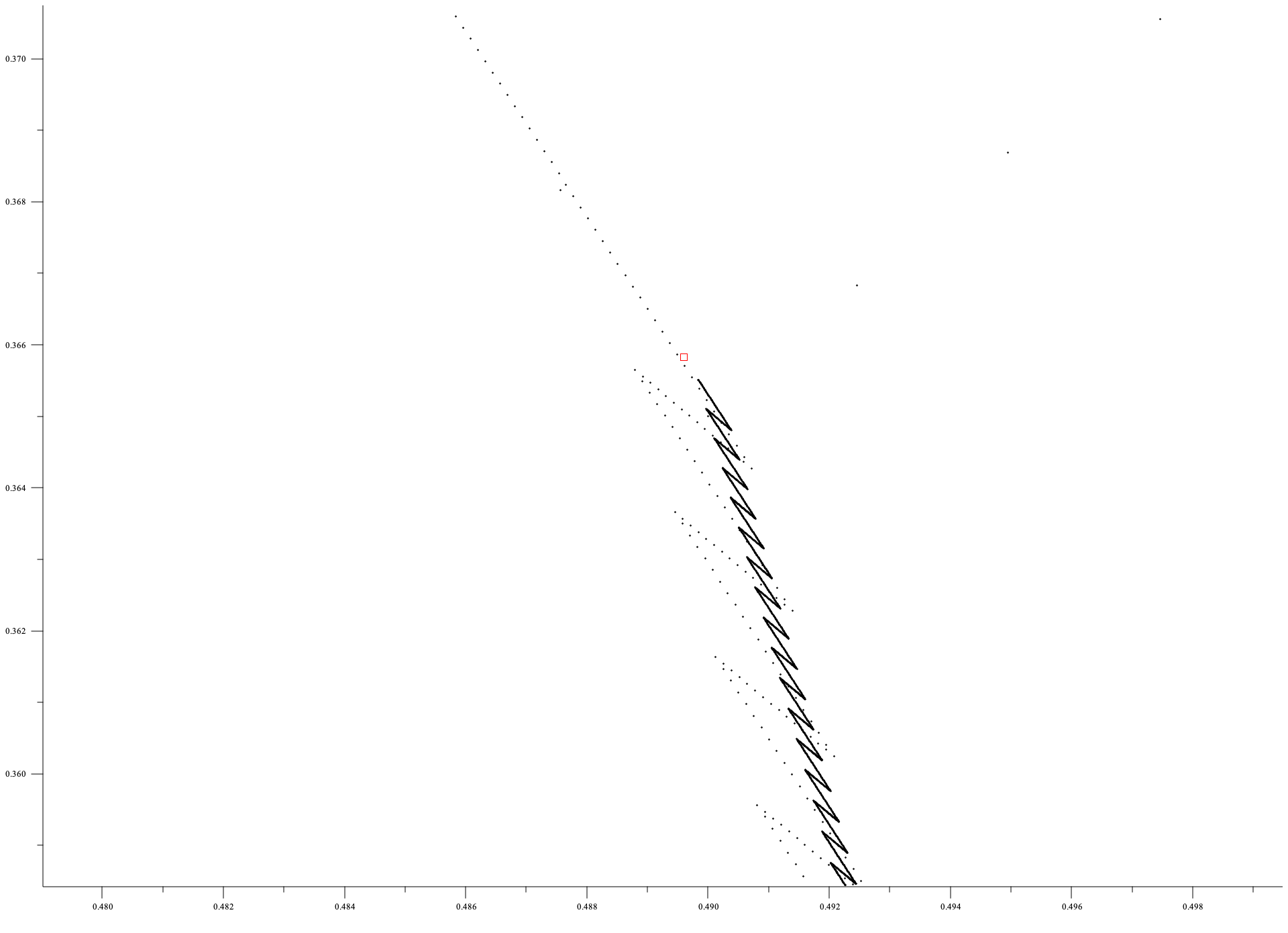}
		\caption{Closer view; two generations are
		apparent.}
		\label{fig:lem-3}
	\end{subfigure}
	\begin{subfigure}[b]{0.40\textwidth}
		\centering
		\includegraphics[width=.99\textwidth]{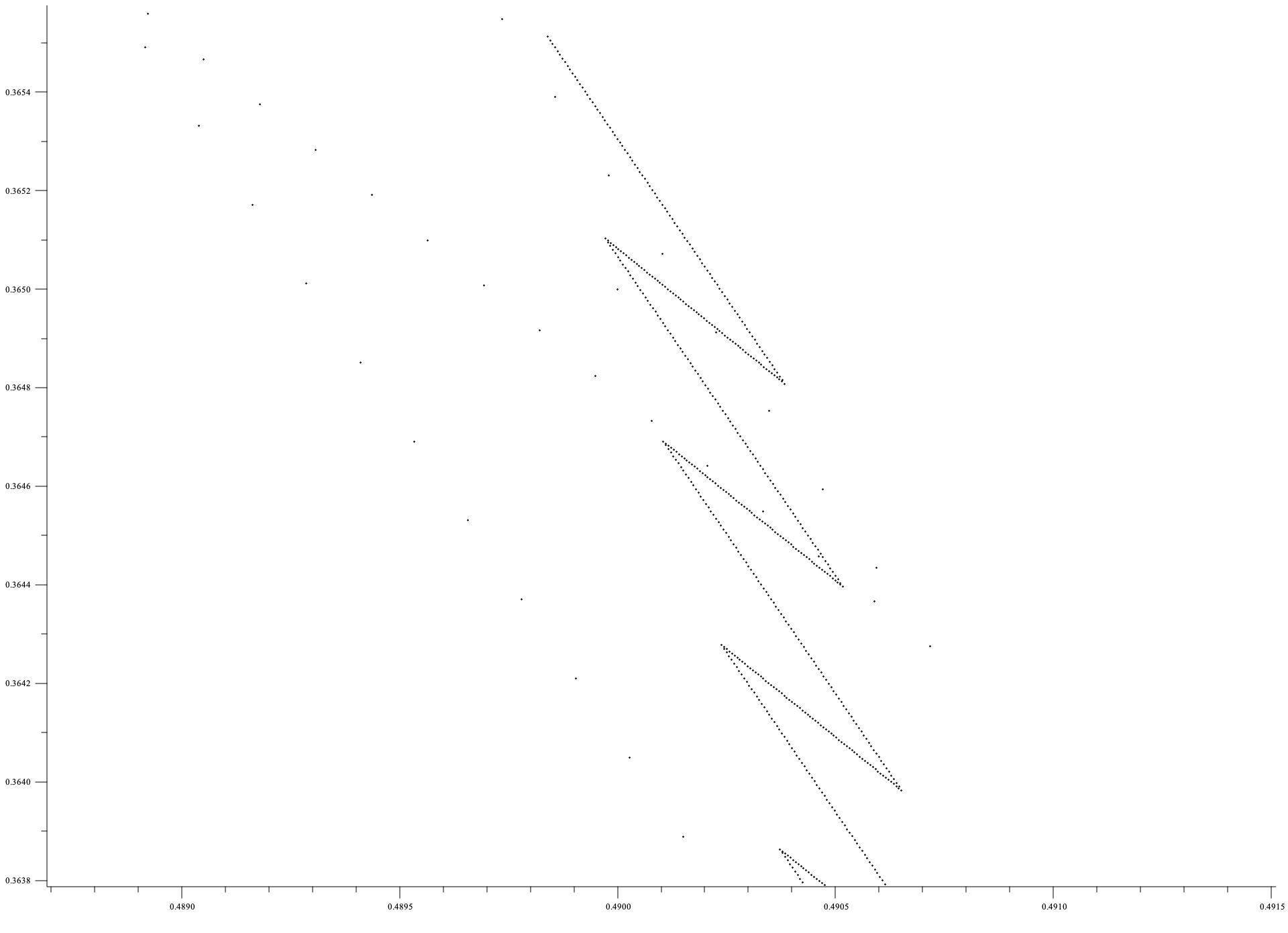}
		\caption{Closer again zoom of Figure \ref{fig:lem-3}.}
		\label{fig:lem-4}
	\end{subfigure}
	\begin{subfigure}[b]{0.40\textwidth}
		\centering
		\includegraphics[width=.99\textwidth]{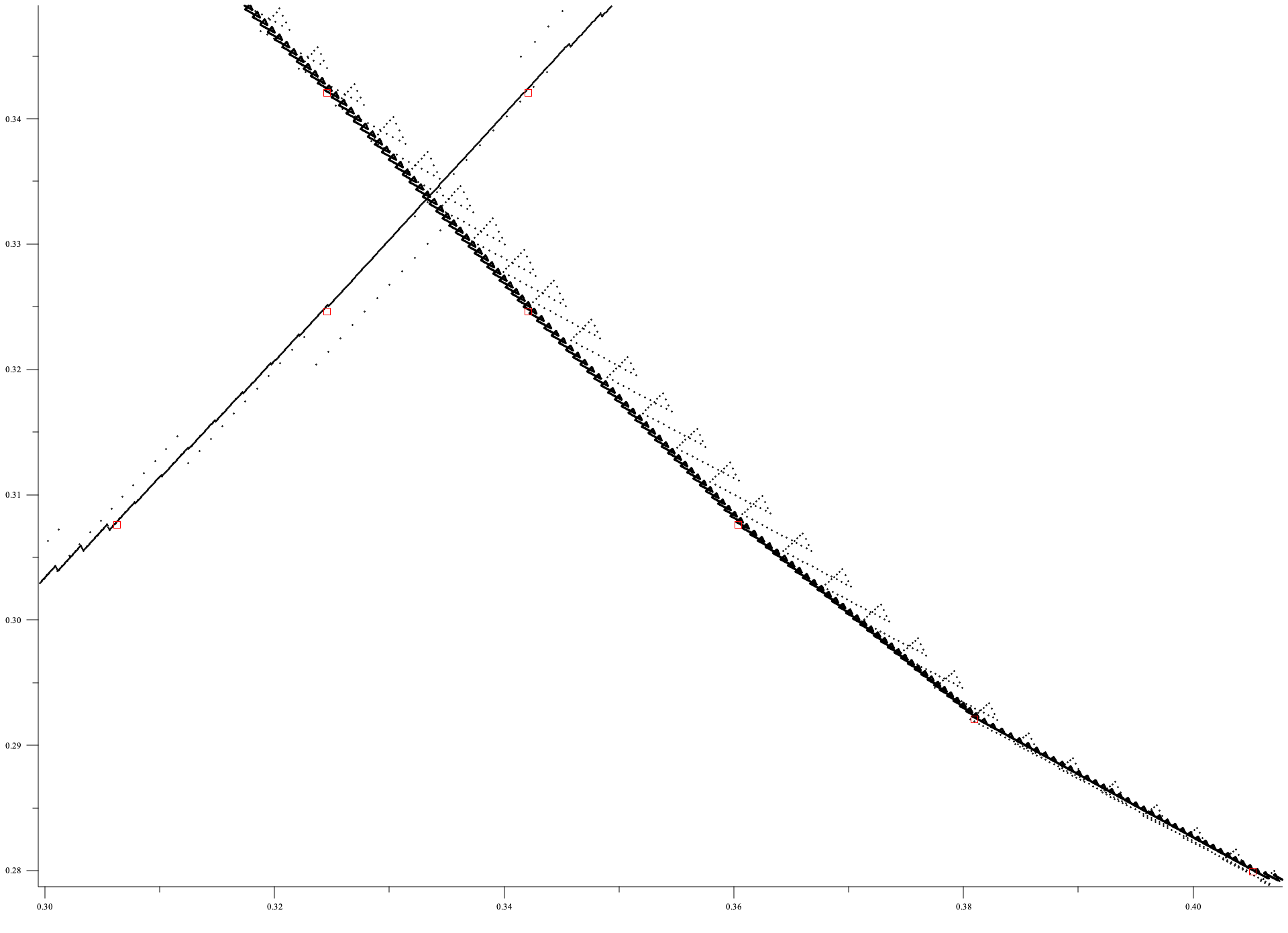}
		\caption{Near the centre of Figure \ref{fig:lem-1}.}
		\label{fig:lem-5}
	\end{subfigure}
	\begin{subfigure}[b]{0.40\textwidth}
		\centering
		\includegraphics[width=.99\textwidth]{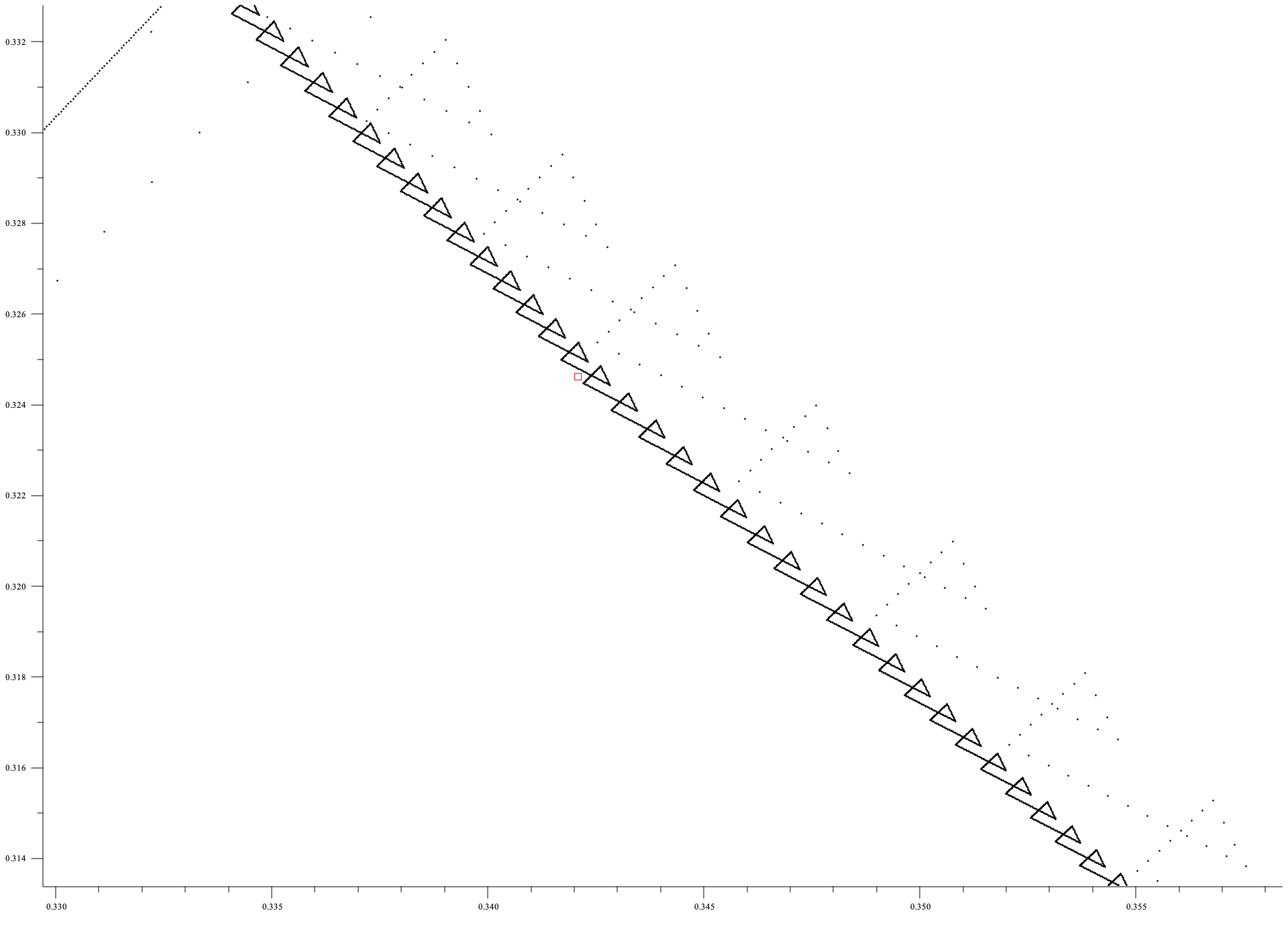}
		\caption{Further zoom of Figure \ref{fig:lem-5}.  }
		\label{fig:lem-6}
	\end{subfigure}
	\begin{subfigure}[b]{0.40\textwidth}
		\centering
		\includegraphics[width=.99\textwidth]{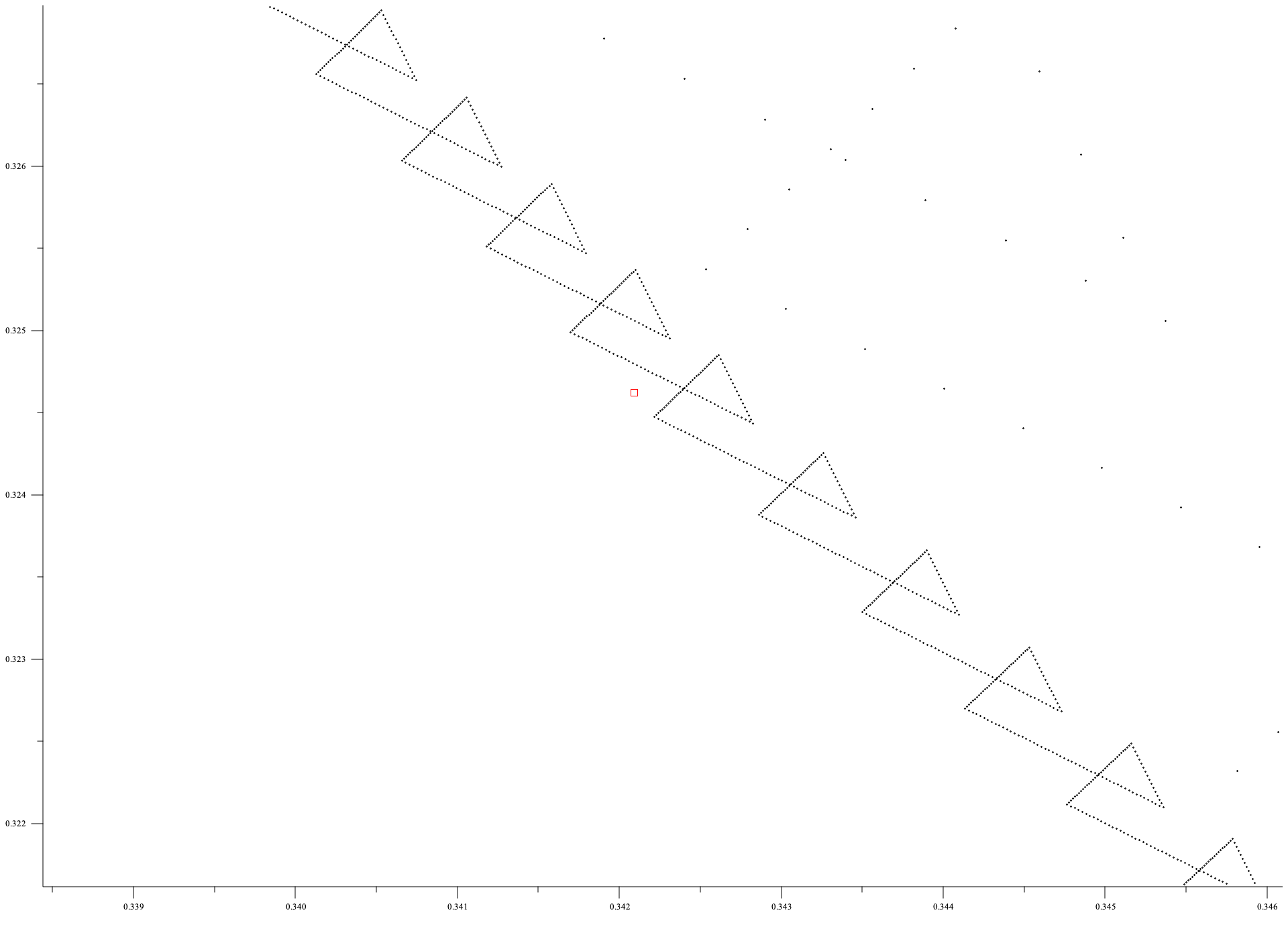}
		\caption{An even closer view of Figure \ref{fig:lem-6}.}
		\label{fig:s4}
	\end{subfigure}
     \caption{Illustration of approximation of the polygonal curve $C^g$
	     by the relative frequencies of the sequence constructed
	     according to Algorithm 1. Two generations were used. Red
	     squares are vertices of $C^g$. See Section \ref{subsec:illustration} 
      for details.\label{fig:lem}}\vspace*{-6pt}
\end{figure}

We set $T=12$ and $V=30$ (number of nodes in the polygonal approximation
$C^g$ for both $g=1,2$ to make the figure less cluttered) and we iterated
long enough to go around the lemniscate twice, which resulted in a sequence
of length 85677.  As the construction proceeded, $\tilde{T}$ went from 9 up
to 164 and $\tilde{\ell}$ went from 1 or 2 for the first few segments up to
69 for the last (with the largest value being 104).  The results can be
seen in Figure \ref{fig:lem}  which plots the achieved relative frequencies
at different zoom levels. The small red squares are the vertices of $C^g$.

\subsection{Construction of \texorpdfstring{$x$}{x} such that \texorpdfstring{$\co\CP(r^x)=\Delta^k$}{the Convex Hull of Cluster Points is the Simplex}}
\label{subsec:maximally-nonstochastic}
How much of the simplex can we cover with $\CP(r^x)$? This question
is poorly posed as (it seems) we can only ever construct one dimensional
sets that are the set of cluster points of $r^x$. However, for inducing
an upper prevision, all that matters is the convex hull of $\CP(r^x)$. 
We now show the convex hull of $\CP(r^x)$ can be made
as large as conceivable with simpler and more explicit construction: 
\begin{lemma}
	\label{lem:extreme-case}
	Suppose $k\in\naturals$. There exists a sequence
	$x\colon\naturals\rightarrow [k]$ such that $\co(\CP(r^x))=\Delta^k$.
\end{lemma}
\begin{proof}
	As before, our proof is constructive. 
     Recall $e_1,\ldots,e_k$ are the vertices (and extreme points)
     of the $(k-1)$-simplex,  and $\co\{e_1,\ldots,e_k\}=\Delta^k$.
     We will construct a sequence
	$x$ such that for all $\epsilon>0$, and all $i\in[k]$, $r^x$ visits
	$N_{\epsilon}(e_i)$ infinitely often. Since 
     $\lim_{\epsilon\rightarrow 0} \co \bigcup_{i\in[k]} N_\epsilon(e_i)
     = \Delta^k$ we will
     have achieved the desired result. 

	We again make use of (\ref{eq:r-decomposition}). We will construct
	a sequence $x$ by adding segments ($s$) such that for each
	successive $s$ we drive $r^x$ closer and closer towards one of the
	vertices $e_i$ ($i\in[k]$). In order to do this, at each $s$
	we append $m$ copies of $i$ to the current $x$. Specifically, we
	construct $x$ as follows:
	\begin{algorithm}[h]
	\begin{algorithmic}[1]
	\Require $\phi\colon\naturals\rightarrow\naturals$
	    \State  $x\gets \langle \rangle$
	     \State $s\gets1$
	     \While{true}
	           \State $i \gets s \mod k$   
		   \Comment{Cycle around the vertices of the simplex}
		   \State $m\gets\phi(s+1)-\phi(s)$
		   \State $x\gets x\, i^m $  
		   \Comment{Append $m$ copies of $i$ to $x$}
		   \State $s\gets s+1$
	     \EndWhile
	\end{algorithmic}
	\end{algorithm}

	We need to make $m$
	large enough so that the convex combination coefficient
	$\frac{m}{n+m}$ approaches $1$.  To that end, consider an
	increasing function $\phi\colon\naturals\rightarrow\naturals$ which
	will further restrict later. The role of $\phi$ is to control $n$
	as a function of segment number $s$; that is $n=\phi(s)$ and thus
	$m=\phi(s+1)-\phi(s)$. With this choice, we have
	\[
		\frac{n}{n+m}=\frac{\phi(s)}{\phi(s+1)}\ \ \ \mbox{and}\ \
		\ \ \frac{m}{n+m}=1-\frac{\phi(s)}{\phi(s+1)}.
	\]
	and thus for all $s\in\naturals$
	\begin{equation}
		\label{eq:r-phi-s}
		r^x(\phi(s+1))= \frac{\phi(s)}{\phi(s+1)} r^x(\phi(s)) +
		\left(1- \frac{\phi(s)}{\phi(s+1)}\right)
		r^{x^{+\phi(s)}}\left(\phi(s+1)-\phi(s)\right) .
	\end{equation}
	We demand that
	$\lim_{s\rightarrow\infty}\frac{\phi(s)}{\phi(s+1)}=0$ so that as
	$s$ increases, the second term in (\ref{eq:r-phi-s}) dominates.
	For any $s\in\naturals$,  $r^x(\phi(s))\in\Delta^k$ and thus
	\begin{align*}
            \|e_{s\operatorname{mod} k} - r^x(\phi(s+1))\| & = 
	    \left\|e_{s \operatorname{mod} k} - \frac{\phi(s)}{\phi(s+1)} r^x(\phi(s)) -
		\left(1- \frac{\phi(s)}{\phi(s+1)}\right)
		r^{x^{+\phi(s)}}\!\left(\phi(s+1)-\phi(s)\right)\right\|\\
		&=\left\|\frac{\phi(s)}{\phi(s+1)} e_{s\operatorname{mod} k} -
		\frac{\phi(s)}{\phi(s+1)} r^x(\phi(s))\right\|\\
		&=\frac{\phi(s)}{\phi(s+1)} \left\|e_{s \operatorname{mod} k}-r^x(\phi(s))\right\|\\
		&\le \frac{\phi(s)}{\phi(s+1)} \cdot 2,
	\end{align*}
	where the second line follows from the fact that we constructed
	$x$ such that
	$r^{x^{+\phi(s)}}\!\left(\phi(s+1)-\phi(s)\right)=e_{s\operatorname{mod} k}$. But
	by assumption,
	$\lim_{s\rightarrow\infty}\frac{\phi(s)}{\phi(s+1)}=0$ and hence
	for any $\epsilon>0$ there exists $s_\epsilon$ such that for all
	$i\in[k]$,
	\[
		\left|\left\{s\in\naturals\colon s>s_\epsilon,\ s\operatorname{mod} k=i,\
		\left\|e_i-r^x(\phi(s+1))\right\|\le\epsilon\right\}\right|=\aleph_0.
	\]
	That is, for each $i\in[k]$, each $\epsilon$-neighbourhood of $e_i$
	is visited infinitely often by the sequence $r^x$ and hence
	$\{e_1,\ldots,e_k\}\subseteq\CP(r^x)$.  But since
	$r^x(n)\in\Delta^k$ for all $n\in\naturals$ we conclude that indeed
	$\co(\CP(r^x))=\co(\{e_1,\ldots,e_k\})=\Delta^k$ as required. 
\end{proof}
A suitable choice of $\phi$ is $\phi(s)=\left\lceil\exp(s^\alpha)\right\rceil$ for
some $\alpha>1$, in which case $\frac{\phi(s)}{\phi(s+1)}\approx
\exp\left(-\alpha s^{\alpha-1}\right)$.  An argument as in the proof of Theorem \ref{th:CP-r-C}
would show that $S_k\coloneqq\bigcup_{i\in[k]}
l\left(e_i,e_{(i+1)\operatorname{mod} k}\right)\subseteq\CP(r^x)$.  Observe that when  $k=3$,
$S_k=\partial\Delta^k$, but for $k>4$, that is not true, even though $\co
(S_k)=\partial\Delta^k$.

\section{Ivanenko's Sampling Nets}
\label{app:ivanenkonets}
Ivanenko seeks to abstract away from sequences and hence defines a the notion of a \textit{sampling net}.

First, we recall the standard definition of a net in topology, which generalizes sequences in an important way. A directed set $(\Lambda, \geq)$ consists of an arbitrary set $\Lambda$ and a direction $\leq$ on it, which satisfies the following properties:
\begin{enumerate}[label=\textbf{DIR\arabic*.}, ref=DIR\arabic*]
\item If $\lambda \in \Lambda$, then $\lambda \leq \lambda$ \quad (reflexivity).
\item If $\lambda_1, \lambda_2, \lambda_3 \in \Lambda$ and $\lambda_1 \leq \lambda_2$ and $\lambda_2 \leq \lambda_3$, then $\lambda_1 \leq \lambda_3$ \quad (transitivity).
\item If $\lambda_1,\lambda_2 \in \Lambda$, then $\exists \lambda_3 \in \Lambda$ such that $\lambda_1 \leq \lambda_3$ and $\lambda_2 \leq \lambda_3$ \quad (upper bound).
\end{enumerate}
That is, $\leq$ is a pre-order and two any two elements there exists a common upper bound.

A \textit{net}\footnote{Ivanenko calls this ``directedness''.} is a function $\varphi\colon \Lambda \rightarrow \mathcal{Y}$, where the domain is a directed set $(\Lambda,\leq)$.

Fix a topology $\tau$ on $\mathcal{T}$. We say that a point $T \in \mathcal{T}$ is a cluster point of the net $\varphi$, $T \in \CP(\varphi)$, if:
\[
    \forall N \text{, where N is any neighbourhood of } T \text{ with respect to } \mathcal{T} \text{, } \forall \lambda_0 \in \Lambda: \exists \lambda_1: \lambda_0 \leq \lambda_1: \varphi(\lambda_1) \in N.
\]

As an example, a sequence is a net, where $\Lambda = \mathbb{N}$ and $\leq$ is the familiar order on the natural numbers.

Define the \textit{space of samples from $\Omega$} as:
\[
\Omega^{(\infty)} \coloneqq \bigcup_{n=1}^\infty \Omega^n, \quad \Omega^n \coloneqq \underbrace{\Omega \times \cdots \times \Omega}_{n \text{ times}}.
\]
\citet{ivanenkobook} then calls a net $\varphi\colon \Lambda \rightarrow \Omega^{(\infty)}$, which takes values in the space of samples from $\Omega$ a \textit{sampling directedness} or \textit{sampling net} (\eg in \citet{ivanenkoonregularities}). To such a net, Ivanenko associates a net $n_\varphi \colon (\Lambda,\leq) \rightarrow \mathbb{N}$ of ``frequency counts'':
\[
n_\varphi(\lambda) \coloneqq n \text{ such that } \varphi(\lambda) \in \Omega^n.
\]

Furthermore, a corresponding net of relative frequencies $P_\varphi\colon (\Lambda,\leq) \rightarrow \PF(\Omega)$ can be defined as follows:
\[
P_\varphi(\lambda) \coloneqq A \mapsto \frac{1}{n_\varphi(\lambda)} \sum_{i=1}^{n_\varphi(\lambda)} \chi_A(\varphi_i(\lambda)),
\]
where $A \subseteq \Omega$ and $\varphi(\lambda)=(\varphi_1(\lambda),\ldots,\varphi_n(\lambda)) \in \Omega^{n_\varphi(\lambda)}$.

The non-empty closed set of limit points of $P_\varphi$ Ivanenko calls the \textit{statistical regularity} of the sampling net $\varphi$. The main result in \citep{ivanenkobook} is then the following. Call any non-empty weak* closed subset $\mathcal{P} \subseteq \PF(\Omega)$ a \textit{regularity}.
\begin{proposition}[\protect{\citet[Theorem 4.2]{ivanenkobook}}]
\label{prop:ivanenkomaintheorem}
    Any sampling net has a regularity, and any regularity is
the regularity of some sampling net [..].
\end{proposition}

Thus, the concept of a sampling net is in a satisfying one-to-one correspondence to that of a weak* closed set of linear previsions (in \citep{ivanenkobook}, finitely additive probabilities). Nonetheless, we remain skeptical about the utility of this concept and raise the question: what is the \textit{meaning} of a sampling net?
\citet{ivanenko2017expected} give an intuition: take for example $\Lambda = \mR^+$ with the familiar order $\leq$. Then, we could interpret $P_\varphi(t)(A)$ as the frequency of the number of hits in $A$ of the observations $(\omega_1,..,\omega_{n_\varphi(t)})$ that are performed at time $t$. Importantly, at any time $t$ we could record a totally different number of observations, since that number is itself given by the net $n_\varphi$. And to obtain the ``relative frequencies'' at time $t$, we consider only data which was observed at time $t$ and completely neglect the past. Contrast this with the case of a sequence: at each time step, we make exactly one observation, and to compute the relative frequencies at time $n$, we use \textit{the complete past}.

Moreover, in the above example, the directed set $\mR^+$ with the familiar order $\leq$ was easily intuited. However, Ivanenko's proof for the direction ``to any regularity there exists a corresponding sampling net'' is non-constructive in the sense that he uses the exotic directed set
\[
\Lambda \coloneqq \mR^+ \times (\linfty)^{(\infty)} \times P,
\]
where $P \subseteq \PF(\Omega)$ is the regularity in question and
\[
(\linfty)^{(\infty)} \coloneqq \bigcup_{n=1}^\infty (\linfty)^n, \quad (\linfty)^n \coloneqq \underbrace{\linfty \times \cdots \times \linfty}_{n \text{ times}}.
\]
It is not clear to us what a realistic interpretation for 
such a sampling net $\varphi \colon \Lambda \rightarrow \Omega$
would look like.
\section{Pathological or Normal?}
\label{app:pathornormal}
\vspace*{-8mm}
\hfill\begin{minipage}{0.55\textwidth}
    \footnotesize{\it
Much of the confusion about probability arises because the true depth of
the law of large numbers as an extremely hard analytical assertion is not
appreciated at all.} --- Detlef D\"{u}rr and Stefan Teufel
\citeyearpar[p.\@\xspace 62]{durr2009}
\end{minipage}

\label{sec:pathologies-or-norm}

When one looks at \emph{finite} sequences $x\colon[n]\rightarrow[2]$, there
is a simple counting argument using the binomial theorem that illustrates
that the vast majority of the $2^n$ possible sequences have roughly equal
numbers of elements with values of 1 and 2. If one \emph{assumes} that an
infinite sequence $x\colon\naturals\rightarrow[2]$ is generated i.i.d. then
this argument can be used to prove the law of large numbers, which  ensures
``most'' sequences have relative frequencies which converge. 

Hence the construction, as illustrated in the present paper, of sequences
$x\colon\naturals\rightarrow [k]$ with \emph{divergent} relative
frequencies naturally raises the question of how contrived they are. That
is, are we examining a rare pathology, or something ``normal'' that we
might actually encounter in the world?  We will refer to sequences whose
relative frequencies converge as ``\Def{stochastic sequences}'' and
sequences whose relative 
frequencies do not converge as  ``\Def{non-stochastic sequences}''\footnote{
	\label{footnote:rivas}
    This dichotomy would appear clear-cut, but there is a subtlety: there 
    exist sequences  such that 
    $(x_i)_{i\in\naturals}$ and $(y_i)_{i\in\naturals}$ are both
    stochastic, but  the
    \emph{joint} sequence $((x_i,y_i))_{i\in\naturals}$ is 
    \emph{non}-stochastic \citep{rivas2019role}; that is, the marginal relative
    frequencies converge, but the joint relative frequencies do not!
}.

The classical law of large numbers suggests 
that indeed ``almost all'' sequences are stochastic,
and therefore, by such reasoning, the non-stochastic sequences with which
we have concerned ourselves in the present paper are indeed pathological
exceptions.  In this appendix we will argue:
\begin{enumerate}[nolistsep]
    \item This very much depends upon what one means by ``rare'' or ``almost
	    all'' and there are many choices, and the only real
	    argument in favour of the usual ones (which declare
	    non-stochastic sequences rare) is familiarity --- different
	    notions of ``typicality'' (for that is what is at issue) lead
	    to very different conclusions. Specifically, there are choices
	    (arguably just as ``natural'' as the familiar ones) which imply
	    that rather than non-stochastic sequences being rare, they are
	    in fact the norm in a very strong sense.
    \item Nevertheless, none  of the mathematical nuances of the previous
	    point allow one to conclude \emph{anything} about the empirical
	    prevalence of stochastic or non-stochastic sequences  in the
	    world.  Indeed, no purely mathematical reasoning allows one to
	    draw such conclusions, unless one wishes to appeal to some
	    conception of a Kantian ``synthetic a priori.''
\end{enumerate}
We will first explore what can said from a purely mathematical perspective,
illustrating that there is a surprising amount of freedom of choice in 
precisely posing the problem, and that the choices are consequential. Then in 
Subsection \ref{subsec:real-measured-sequences}
we examine the question of prevalence of non-stochastic sequences
actually in the world.  

\subsection{The Mathematical Argument ---  The Choices to be Made}

The classical Law of Large Numbers says ``almost all sequences'' are
stochastic.  But the ``almost all'' claim  comes from the mapping of
sequences to real numbers in $[0,1]$  and then making a claim that ``almost
all'' numbers correspond to stochastic sequences.  Thus there are at least
three choices being made here:
\begin{description}[nolistsep]
    \item[Mapping from Sequences to Real Numbers]  
	The choice of mapping from sequences to real numbers, to enable to
	use of some notion of typicality on $[0,1]$ to gauge how common
	stochastic sequences are.
   \item[Notion of Typicality] The notion of typicality to be used (e.g. 
         Cardinality, Hausdorff dimension, Category or Measure).
   \item[Specific Index of Typicality] 
       Within the above choice of notion of typicality, the particular 
       choice of typicality index,  \eg the measure or 
       topology that underpins the notion of typicality.
\end{description}
The choices for the classical law of large numbers are 
1) $k$-ary positional representation; 
2) a $\sigma$-additive measure on $[0,1]$; 
3) The Lebesgue measure.  As we shall summarize
below, each of these three choices substantially affects the theoretical
preponderance of non-stochastic sequences.

That there are alternate choices that lead to the unusual conclusion that
non-stochastic sequences  are ``typical'' has been known for some time:
``This result may be interpreted to mean that the category analogue of the
strong law of large numbers is false''  \citep[p.\@\xspace 85]{oxtoby1980measure};
see also \citep{mendez1981law}.  The significance of this fact has been
stressed recently \citep{whenlargeisalso,cisewski2018standards}.  And it
has been observed that the introduction of alternate topologies can change
whether sequences are stochastic \citep{khrennikov2013p}. However, the
strongest results arise in number theory, motivated by the notion of a
``normal number.'' 

\subsection{Notions of Typicality --- 
Cardinality, Dimension, Comeagreness, and Measure}

Let $\Def{\mathcal{S}}$ (resp.~$\Def{\mathcal{N}}$) denote the set of
stochastic (resp.~non-stochastic) sequences $\naturals\rightarrow [k]$.
That is, $\Def{\Scal\coloneqq\{x\colon\naturals\rightarrow[k]\colon 
	\lim_{n\rightarrow \infty} r^x(n)\text{ exists}\}}$ and 
	$\Def{\Ncal\coloneqq [k]^\naturals\setminus\Scal}$.
(For simplicity, and alignment with Appendix~\ref{app:construction}, we restrict
ourselves to sequences whose domain is $[k]$.)

In order to make a claim regarding the relative preponderance of stochastic
versus non-stochastic sequences, they are often mapped onto the unit
interval.\footnote{The only attempts to judge the relative
	sizes of $\mathcal{S}$ and $\mathcal{N}$ which do not rely on such
	a mapping are described in the first and third of the cases
	listed below, and rely upon imposing a topology directly on 
	the set of sequences $[k]^\naturals$.}
In such cases, the question of relative preponderance of classes of
sequences is reduced to that of a question concerning the relative
preponderance of classes of subsets of $[0,1]$.  The question then arises
of how to measure the size of such subsets. Unlike in the finite case
mentioned above, merely counting (i.e. determining the cardinality of the
respective subsets) is hardly adequate, as it is easy to argue that
$|\Scal|=|\Ncal|=\aleph_1$.  There are three notions that have
been used to compare the size of $\Scal$ and $\Ncal$:
\begin{description}[nolistsep]
	\item[Measure] A countably additive measure, usually the Lebesgue
		measure on $[0,1]$.
	\item[Meagre / Comeagre] 
		A subset $S$ of a topological space
		$X$ is \Def{meagre} if it is a countable union of nowehre
		dense sets (i.e. sets whose closure has empty interior). 
		A set $S$ is \Def{comeagre} (residual) if $X\setminus S$ is
		meagre.
	\item[Dimension]  A variety of fractal dimensions, such as the
		Hausdorff dimension, have also been used to judge the size
		of non-stochastic sequences (and the  numbers they induce);
		however for space considerations we omit discussion of
		these results\footnote{
			See for example
			\citep{eggleston1949fractional,
				olsen2004extremely,
				gu2011effective,
				bishop2017fractals,
				albeverio2017non}.  
		}. 
\end{description}

Some of the results obtained in the literature are summarized below. The
object is not to state them in an entirely formal manner, or even to
describe them in their full generality. Rather we simply wish to show the
diversity of conclusions available by tweaking the three choices enumerated
above. If no representation is mentioned, the usual $k$-ary positional
representation is used, whereby $\tilde{x}\in[0,1]$ is constructed from
$x\colon\naturals\rightarrow [k]$ via
$\Def{\tilde{x}\coloneqq\sum_{i\in\naturals} (x_i-1) k^{-i}}$ (the
$(x_i-1)$ term is required because our sequences map to
$[k]=\{1,\ldots,k\}$).  Obviously every $x\in[k]^\naturals$ maps to some
$\tilde{x}\in [0,1]$; and every $z\in [0,1]$ corresponds to at least one
$x\in [k]^\naturals$ (recalling we have to handle the situation that, when
$k=10$ for example, $0.4\overline{9}=0.5\overline{0}$, where $\overline{i}$
means that $i$ is repeated infinitely, and thus there are two sequences
$x_1,x_2\in[k]^\naturals$ such that $1/2=\tilde{x}_1=\tilde{x}_2$).
\footnote{This non-uniqueness of the representation will not affect
	the results below because 
	$\{\tilde{x}_1=\tilde{x}_2\in[0,1]\colon x_1\ne
	x_2\}=\mathbb{Q}\cap[0,1]$ and is of cardinality $\aleph_0$,
	whereas $|[k]^\naturals|=|[0,1]|=\aleph_1$.}  
Let  $\Def{\Scalt\coloneqq\{\tilde{x}\in [0,1]\colon x\in \Scal\}}$ and
$\Def{\Ncalt\coloneqq\{\tilde{x}\in [0,1]\colon x\in \Ncal\}}$.

\begin{description}[nolistsep]
	\item[Most (Lebesgue measure) sequences are stochastic]  
		  This is the classical strong law of large numbers. If
		  $\mu_{\mathrm{leb}}$ denotes the Lebesgue measure on
		  $[0,1]$, then the claim is that
		  $\mu_{\mathrm{leb}}(\tilde{\Scal})=1$.
	\item[Most (comeagre) sequences are non-stochastic] 
		 Let $X_i=[2]$ and $X=\bigtimes_{i\in\naturals} X_i$
		 equipped with the product topology. As a set $X\cong
		 [2]^\naturals$. Then $\Ncal\subset X$
		is comeagre \citep{oxtoby1980measure}.
        \item[Most (comeagre) sequences are stochastic] 
		 With different choices of topology, the opposite
		 conclusion holds --- there are topologies such that
		 $\Scal$ is comeagre \citep{calude2003topological}.
	\item[Most (comeagre) sequences are extremely non-stochastic]
		Let $\Def{\tilde{\Ncal}^\star}$ denote the subset 
		of $[0,1]$ of $\tilde{x}$ corresponding to
		$x\in[k]^\naturals$ which satisfy
		$
			\forall i \in [k],\ \liminf_{n\rightarrow\infty} 
			r_i^x(n)=0 
			\text{ and } 
			\limsup_{n\rightarrow\infty}r_i^x(n)=1.
		$
		These sequences are (justifiably) called \Def{extremely
		non-stochastic}; the sequence constructed in Subsection
		\ref{subsec:maximally-nonstochastic} is an
		example.
		Then the set $\tilde{\Ncal}^\star$ is comeagre in the usual 
		topology of real numbers \citep{calude1999most};
		confer \citep[section 7.3]{calude2002information}.
	\item[Most (comeagre) sequences are perversely non-stochastic]
		Denote the set of \Def{perversely
		nonstochastic sequences}
		$\Def{\Ncal^{\star\star}\coloneqq\{x\in[k]^\naturals\colon 
		\CP(r^x)=\Delta^k\}}$. Observe
		$\Ncal^{\star\star}\subset\Ncal^\star$.  
		Let
		$\tilde{\Ncal}^{\star\star}\coloneqq\{\tilde{x}\in [0,1]
		\colon x\in\Ncal^{\star\star}\}$ (what 
		\citet{olsen2004extremely} calls ``extremely non-normal
		numbers'', but we use ``extremely''  for the larger set
		$\tilde{\Ncal}^\star$).   Then $\tilde{\Ncal}^{\star\star}$
		is comeagre (in the usual topology of real numbers) 
		\citep{aveni2022most,olsen2004extremely}.
		An even stronger result holds.
		Let $A$ denote a (not necessarily uniform)
		finite averaging operator and let
		$\Def{\Ncal^{\star\star\star}\coloneqq 
		\{x\in[k]^\naturals\colon 
		\CP(A(r^x))=\Delta^k\}}$ and
		$\Def{\tilde{\Ncal}^{\star\star\star}\coloneqq\{
		\tilde{x}\in [0,1] \colon x\in\Ncal^{\star\star\star}\}}$.
		Observe $\Ncal^{\star\star\star}\subset\Ncal^{\star\star}$
		and $\tilde{\Ncal}^{\star\star\star}\subset 
		\tilde{\Ncal}^{\star\star}$.
		Then $\tilde{\Ncal}^{\star\star\star}$ is also comeagre
		\citep{stylianou2020typical}!
	\item[Most (Lebesgue measure) sequences are non-stochastic]
		There exist a range of representations of real numbers
		called  $Q^*$-representations ($Q^*$ is a
		$k\times\infty$ matrix valued parameter of the
		representation); see \citep[Section
		4]{albeverio2005topological} for details.  Let
		$\Def{\tilde{x}^{Q^*}}$ denote the $Q^*$-representation of
		a sequence $x\in[k]^\naturals$, and
		$\Def{\tilde{\Scal}^{Q^*}\coloneqq\left\{\tilde{x}^{Q^*}\colon
		x\in\Scal\right\}}$ and
		$\Def{\tilde{\Ncal}^{Q^*}\coloneqq\left\{\tilde{x}^{Q^*}\colon
		x\in\Ncal\right\}}$.
		Then there exist $Q^*$ such that
		$\mu_\mathrm{leb}\left(\tilde{\Ncal}^{Q^*}\right)=1$.
		\citep[p.\@\xspace 627]{albeverio2005topological}.  Thus if the
		size of $\Ncal$ is judged via certain $Q^*$
		representations, Lebesgue almost all sequences are
		non-stochastic!  
\end{description}

An obvious conclusion to draw from the above examples is that in answering
the question of the preponderance of non-stochastic sequences, one can get
essentially whatever answer one wants by choosing a range of different
precise formulations of the question. At the very least, this should make
us skeptical of any purely mathematical attempts to reason whether one
might expect to encounter non-stochastic sequences in practice --- the
topic to which we now turn.

\subsection{Typical Real Sequences}
\label{subsec:real-measured-sequences}
\vspace*{-8mm}
\hfill\begin{minipage}{0.45\textwidth}
    \footnotesize{\it
	The laws of large numbers cannot be applied for describing the
	statistical stabilization of frequencies in sampling experiments.}
	\hfill --- Andrei Khrennikov 
	\citeyearpar[p.\@\xspace 20]{khrennikov2009interpretations}
\end{minipage}

What do the above points imply about the likelihood one will encounter
stochastic or non-stochastic sequences when performing real measurements?

\emph{Nothing}. 

This is not to say that in actuality we will often encounter non-stochastic
sequences. Rather our point is that no amount of purely theoretical
reasoning will be able to tell us in advance how ``likely'' it is to do so.
What is at issue is whether stochastic sequences are in fact ``typical'' in
our world. 

Perhaps the most surprising thing about the mathematical results summarized
above is the extent to which different notions of typicality affect the
conclusions. This raises the question of whether some notions of typicality
are more justified when wishing to consider real sequences that have been
measured in the world.  In the study of physics (especially aspects of
physics that are apparently intrinsically statistical) such questions have
been raised, and below we briefly summarize what is known. 

Traditionally, ``probability'' is considered as a primitive, and notions of
typicality are derived from that in terms of their ``probability'' of
occurring.  And the above examples illustrate that attempts to argue for the Lebesgue measure having a privileged role as the ``right'' notion of
typicality are barking up the wrong tree; confer
\citep{pitowsky2012typicality}.  But this will not do for our question.
Typicality is a more fundamental notion
\citep{galvan2006bohmian,galvan2007typicality} --- arguably the ``mother of
all'' notions of probability  \citep{goldstein2012typicality}.  Typicality
is at the core of questions of non-stochastic randomness in physics, thus
(consistent with the perspective of the present paper)  leading to
non-additive measures of typicality \citep{galvan2022non} (essentially
defining a measure of typicality inspired by a coherent upper probability)
which allows the extension to notions of mutual typicality necessary to
reason about situations such as that referred to in footnote
\ref{footnote:rivas}.

In fact, typicality plays an even stronger role than answering questions
regarding the preponderance of non-stochastic sequences. As 
\citet[p.\@\xspace 36]{durr2021typicality} observe ``the notion of typicality is
necessary to understand what the statistical predictions of a physical
theory really mean.'' They note that the usual appeal to the law of large
numbers misses the point because while its conclusion is true (convergence
of relative frequencies) \emph{if} one sees typical sequences, but 
	``What needs to be explained is why we only see typical sequences!
	That’s actually the deep question underlying the meaning of
	probability theory from its very beginning \ldots'' 
	\citep[p.\@\xspace 37]{durr2021typicality}.
In classical mechanics, appeal to Liouville's theorem suggests an
``invariant measure'' as being a natural choice; in the quantum realm,
there is an analogous choice (invariant to Bohmian flow) 
\citep[p.\@\xspace 41]{durr2021typicality}. But these situations are rather
special from the perspective of a statistician.
The situation is well summarized by  \citet[p.\@\xspace 130]{durr2001bohmian}:
	``What is typicality? It is a notion for defining the smallness of
	sets of (mathematically inevitable) exceptions and thus
	permitting the formulation of law of large numbers type statements.
	Smallness is usually defined in terms of a measure. What
	determines the measure? In physics, the physical theory.''
Confer \citep[Chapter 4]{durr2009} who observe that from a
\emph{scientific} perspective (where one wants to make claims about the
world) establishing the pre-conditions for the law of large numbers to
hold is ``exceedingly difficult''\footnote{The example they give is for the
	Galton board, or quincunx, a device often appealed to in order to teach the
	reality of the central limit theorem --- an even stronger claim than the
	law of large numbers. The irony is that it is rarely checked empirically.
	And when it has been, it has been found to be untrue!
	\citep[Figure 8]{bagnold1983nature}.}.

Very well one might say, but the arguments in favor of typicality of
non-stochastic sequences given above all rely on topological arguments, or
unusual encodings of sequences to numbers. What is the justification for
topological notions of typicality when considering sequences of measurements
obtained from the world?  \citet[p.\@\xspace 270]{sklar2000topology} has actually
argued that the topological perspective might offer a foundational
perspective with \emph{fewer} opportunities for claims of arbitrariness
than measure theoretical approaches.  See also \citep[p.\@\xspace
185]{sklar1995physics} and the discussion in \citep[Chapter
4]{guttmann1999concept} which reframes the problem away from typicality to
viewing the whole question from an approximation perspective where the
notion of smallness of sets is naturally one of meagreness.  Our point is
that even within the restricted realm of physics, there are compelling
arguments at least not to take the measure-based notion of typicality for
granted. Once that is accepted, non-stochastic sequences seem less unusual.

\subsection{Violations of the Law of Large Numbers}

A typical universe is in equilibrium; but ``our universe is atypical or in
non-equilibrium'' \citep[p.\@\xspace 81]{durr2009} and ``what renders knowledge at
all possible is nonequilibrium'' \citep[p.\@\xspace 886]{durr1992quantum} so we
should not be surprised if it is not ``typical''. And indeed that is what
we see as long as we look: ``The so-called law of large numbers is also
invalid for social systems with finite elements during transition''
\citep{chen1991nonequilibrium}.
\citet{gorban2011statistical,gorban2017statistical,gorban2018randomness}
has documented many examples of real phenomena failing to be statistically
stable.  Such failures are held to explain departures from ``normal''
distributions \citep{philip1987some}. But more importantly, they mean we
should not expect even convergence of relative frequencies in
non-equilibrium  situations. 

Such was the conclusion of Prigogine in his ground-breaking studies of
non-equilibrium thermodynamics  where he spoke of a ``breakdown of the `law
of large numbers'{}'' \citep[p.\@\xspace 9 and 228]{Nicolos1977}; see also
\citep[p.\@\xspace 781]{prigogine1978time}, \citep[p.\@\xspace 180]{prigogine1984order}
and \citep[p.\@\xspace 131]{prigogine1980}.  And more recently, studies of the use
of machine learning systems ``in the wild'' have recognized that
non-stochasticity is not so exotic after all
\citep{katsikopoulos2021classification}.  Thus perhaps its time to
downgrade this ``law'' of nature.

\subsection{Repeal of the Law of Large Numbers}
\vspace*{-8mm}
\hfill\begin{minipage}{0.45\textwidth}
    \footnotesize{\it
	A typical universe shows statistical regularities as we
	perceive them in a long run of coin tosses. It looks as if
	objective chance is at work, while in truth it is not.
	There is no chance. That is the basis of mathematical
	probability theory.} \hfill ---  
	Detlef D\"{u}rr and Stefan Teufel \citeyearpar[p.\@\xspace 64]{durr2009}. 
\end{minipage}

 \citet{desrosieres1998politics} in his history of statistical reasoning
 has observed the awe with which stable frequencies were viewed when they
 were first encountered; the effect been interpreted as a hidden divine
 order\footnote{``I know of scarcely anything so apt 
	to impress the imagination as
	the wonderful form of cosmic order expressed by the `Law of
	Frequency of Error.' The law would have been personified by the
	Greeks and deified, if they had known of it. It reigns with
	serenity and in complete self-effacement amidst the wildest
	confusion. The huger the mob, and the greater the apparent anarchy,
	the more perfect is its sway. It is the supreme law of Unreason.
	Whenever a large sample of chaotic elements are taken in hand and
	marshalled in the order of their magnitude, an unsuspected and most
	beautiful form of regularity proves to have been latent all
	along.''
	\citep[p.\@\xspace 66]{galton1889natural}. See also  \citep{rose2016end}
	for a recent discussion on statistical normality.
}.
And indeed in many practical situations, stable frequencies  do 
arise.  But that does not mean we should take such situations as the only
ones that can occur. We may well legitimately call them ``normal.'' But we
can better understand the normal by studying the pathological \citep[p.\@\xspace
19--20]{canguilhem1978normal}. Ironically in his attempt to clarify the
notion of ``normal'' \citet[p.\@\xspace 103]{canguilhem1978normal}  
considered whether ``normal'' was simply ``average''
and concluded ``the concepts of norm and average must be considered as two
different concepts''. As we have seen, averages can indeed be far from
normal, and potentially quite often.

Perhaps we have been misled by the strange name given to the famous theorem
we are considering: by calling it a ``law'' we are inheriting a lot of
baggage as to what we mean by that, baggage that has been traced to notions
of divine origin \citep{zilsel1942genesis}\footnote{The contrary views
	regarding the historical origin of the notion of a scientific law 
	\citep{milton1981origin, ruby1986origins,weinert1995laws} do not 
	contradict our point.}
of lawfulness. And we hanker after lawfulness:
\begin{quote}
	We \ldots  naturally hope that the world is orderly. We like it
	that way...  All of us \ldots  find this idea sustaining. It
	controls confusion, it makes the world seem more intelligible. But
	suppose the world should happen in fact to be not very
	intelligible? Or suppose merely that we do not know it to be so?
	Might it not then be our duty to admit these distressing facts?
	\citep[p.\@\xspace 199]{midgley2013science}
\end{quote}
Perhaps the theory of imprecise probabilities presented in this
paper which we have grounded in the instability of relative
frequencies may help us to admit this ``distressing fact.''  It does
suggest to us that the law of large numbers, while a fine and true theorem,
as a ``law'' might be in need of  repealing. 



\end{document}